\documentclass[a4paper,12pt,leqno]{article}
\usepackage[latin1]{inputenc} 
\usepackage[T1]{fontenc}   
\usepackage[french, english]{babel}  
\usepackage{typearea}
\usepackage{color}

\usepackage{textcomp}
\usepackage{amsmath}
\usepackage{amsfonts}
\usepackage{amssymb}
\usepackage{tikz}
\usepackage{fancyhdr}
\usepackage{url}
\usepackage{mathrsfs}
\usepackage{dsfont}
\usepackage{cancel}
\usepackage[justification=centering]{caption}

\usepackage{hyperref}

\usepackage[all]{xy}

\usepackage[a4paper,top=3cm,bottom=3cm,left=2cm,right=2cm]
{geometry}

\usepackage{tikz}
 \usetikzlibrary{patterns}
 
 \numberwithin{equation}{section}

\newtheorem{df}{Definition}[section]

\newtheorem{prop}[df]{Proposition}
\newtheorem{thm}[df]{Theorem}

\newtheorem{lem}[df]{Lemma}
\newtheorem{crl}[df]{Corollary}
\newcommand{\prf}{\textit{Proof}}
\newtheorem{rmk}[df]{Remark}

\newcommand{\dbar}{\overline{\partial}}
\newcommand{\ddbar}{\partial\overline{\partial}}

\newcommand{\R}{\mathbb{R}}

\newcommand{\C}{\mathbb{C}}
\newcommand{\Z}{\mathbb{Z}}
\newcommand{\N}{\mathbb{N}}
\newcommand{\D}{\mathbb{D}}

\newcommand{\ric}{\operatorname{Ric}}

\newcommand{\vareps}{\varepsilon}

\newcommand{\supp}{\operatorname{supp}}

\newcommand{\bb}{\boldsymbol{b}}
\newcommand{\bL}{\boldsymbol{L}}
\newcommand*{\longhookrightarrow}
{\ensuremath{\lhook\joinrel\relbar\joinrel\rightarrow}}

\newcommand*{\LargerCdot}{\raisebox{-0.25ex}
{\scalebox{1.5}{$\cdot$}}}
\DeclareRobustCommand*{\point}{\LargerCdot}  

\newcommand{\cqfd}{ \hfill $\square$ }

\newcommand{\comment}[1]{}

\title{Bergman kernels on punctured Riemann surfaces}
\author{\normalsize\textsc{Hugues AUVRAY
\footnote{71 rue Caulaincourt, 75018 Paris, France. 
E-mail: hugues.auvray@gmail.com} and 
Xiaonan MA\footnote{
Universit\'e de Paris, CNRS, 
Institut de Math\'ematiques de Jussieu-Paris Rive Gauche, 
F-75013 Paris, France.
E-mail: xiaonan.ma@imj-prg.fr} 
and George MARINESCU\footnote{
Universit{\"a}t zu K{\"o}ln,  
Mathematisches Institut, Weyertal 86-90, 50931 K{\"o}ln, 
Germany. E-mail: gmarines@math.uni-koeln.de}}}
\date{}

\makeatletter
\def\blfootnote{\xdef\@thefnmark{}\@footnotetext}
\makeatother

\begin{document} 
\blfootnote{First author is partially supported by ANR contract 
ANR-14-CE25-0010. 
Second author is partially supported by ANR contract 
ANR-14-CE25-0012-01, NNSFC No.11528103, No.11829102
and funded through the Institutional Strategy of the University of Cologne
within the German Excellence Initiative. 
Third author is partially supported by CRC TRR 191.
            \vspace{1pt}
            \footnoterule
            \vspace{-3pt}}
\makeatletter
\renewcommand\section{\@startsection {section}{1}{\z@}%
                                   {-3.5ex \@plus -1ex \@minus -.2ex}%
                                   {2.3ex \@plus.2ex}%
                                   {\centering\sc\normalsize}}
\renewcommand\subsection{\@startsection{subsection}{2}{\z@}%
                                     {-3.25ex\@plus -1ex \@minus -.2ex}%
                                     {1.5ex \@plus .2ex}%
                                     {\normalsize\sf}}
\renewcommand\subsubsection{\@startsection{subsubsection}{3}
{\z@}%
                                     {-3.25ex\@plus -1ex \@minus -.2ex}%
                                     {1.5ex \@plus .2ex}%
                                     {\normalsize\it}}

\makeatother
 
\maketitle
\abstract{In this paper we consider a punctured Riemann surface 
endowed 
with a Hermitian metric that equals the Poincar\'e metric near
the punctures, and a holomorphic
line bundle that polarizes the metric.
We show that the Bergman kernel can be localized around 
the singularities
and its local model is the Bergman kernel of the punctured unit disc 
endowed with the standard Poincar\'e metric.
One of the technical tools is a new weighted elliptic estimate near 
the punctures, which is uniform with respect to the tensor power.
As a consequence, we obtain an optimal uniform estimate of the 
supremum 
norm of the Bergman kernel, involving a fractional growth order 
of the tensor power. This holds in particular for the Bergman kernel 
of cusp forms of high weight of non-cocompact geometrically finite
Fuchsian groups of first kind without elliptic elements. }

\comment{ \selectlanguage{french}
\abstract{On consid\`ere dans cet article une surface de Riemann
\'epoint\'ee munie d'une m\'etrique hermitienne, 
\'egale \`a la m\'etrique de Poincar\'e pr\`es des points exceptionnels, 
et un fibr\'e en droites holomorphe polarisant la m\'etrique. 
On d\'emontre que le noyau de Bergman peut \^etre localis\'e
autour des singularit\'es, 
et compar\'e au mod\`ele local donn\'e par le noyau de Bergman du 
disque unit\'e \'epoint\'e \'equip\'e de la m\'etrique de 
Po\'encar\'e standard. 
On utilise entre autres une estim\'ee elliptique pr\`es des 
 singularit\'es originale, uniforme en la puissance tensorielle. 
 On obtient en cons\'equence une estim\'ee uniforme optimale 
sur la borne sup\'erieure du noyau de Bergman, 
o\`u la puissance tensorielle appara\^{\i}t avec exposant non entier. 
 Ceci s'applique en particulier dans le cas du noyau de Bergman 
de formes paraboliques de grand poids 
 des groupes fuchsiens de premi\`ere esp\`ece non cocompacts, 
 g\'eom\'etriquement finis, et sans \'el\'ements elliptiques. 
   } }
\selectlanguage{english}

\tableofcontents

\section{Introduction}
 
 In this paper we study the Bergman kernels of a singular Hermitian 
 line bundle
 over a Riemann surface under the assumption that the curvature 
 has singularities of Poincar\'e type at a finite set. 
 Our first result shows that the Bergman kernel can be localized 
 around the singularities
 and its local model is the Bergman kernel of the punctured disc 
endowed with the standard Poincar\'e metric.
The proof follows the principle that the spectral gap of the Kodaira
Laplacian and uniformly elliptic estimates near the singularities
imply the localization of the Bergman kernel \cite{mm}.
By a detailed analysis of the local model we deduce a sharp 
uniform estimate of the supremum norm of the Bergman kernels. 
  
Let us describe our setting. Let $\overline\Sigma$ be a compact 
Riemann surface and let
$D=\{a_1,\ldots,a_N\}\subset\overline\Sigma$ be a finite set.
We consider the punctured Riemann surface 
$\Sigma = \overline{\Sigma}\smallsetminus D$ and a Hermitian form
$\omega_{\Sigma}$ on $\Sigma$.
Let $L$ be a holomorphic line bundle on $\overline{\Sigma}$, 
and let $h$ be a singular Hermitian metric on $L$ such that:
   \begin{itemize}
    \item[($\alpha$)] $h$ is smooth over $\Sigma$, 
                     and for all $j=1,\ldots,N$, 
there is a trivialization of $L$
 in the complex neighborhood $\overline{V_j}$ of $a_j$ in 
 $\overline{\Sigma}$, with associated coordinate $z_j$
 such that $|1|_{h}^2(z_{j})= \big|\!\log(|z_j|^2)\big|$.
    \item[($\beta$)] There exists $\varepsilon>0$ such that the
    (smooth) curvature $R^L$ of $h$ satisfies
   $iR^L\geq\varepsilon\omega_{\Sigma}$ over $\Sigma$ 
   and moreover, $iR^L=\omega_{\Sigma}$ 
	on $V_j:=\overline{V_j}\smallsetminus\{a_j\}$;
                     in particular, $\omega_{\Sigma} = \omega_{\D^*}$ 
                     in the local coordinate $z_j$
                     on $V_j$ 
                     and $(\Sigma, \omega_{\Sigma})$ is complete. 
   \end{itemize}
Here $\omega_{\D^*}$ denotes the Poincar\'{e} metric on the 
punctured unit disc $\D^*$, normalized as follows:
 \begin{equation}   \label{eqn_omegaPcr}
  \omega_{\D^*} := \frac{idz\wedge d\overline{z}}
  {|z|^2\log^2(|z|^2)}\,\cdot
 \end{equation}

 For $p\geq1$, let $h^p:=h^{\otimes p}$ be the metric induced 
 by $h$ on $L^p\vert_{\Sigma}$, 
where $L^p:=L^{\otimes p}$. We denote by 
$H^0_{(2)}(\Sigma,L^p)$ 
the space of ${\bL}^2$-holomorphic sections of $L^p$ 
relative to the metrics 
$h^p$ and $\omega_\Sigma$, 
\begin{equation}\label{e:bs}
H^0_{(2)}(\Sigma,L^p)=\left\{S\in H^0(\Sigma,L^p):\,
\|S\|_{{\bL}^2}^2:=\int_{\Sigma}|S|^2_{h^p}\,
\omega_\Sigma<\infty\right\},
\end{equation}
endowed with the obvious inner product. The sections from 
$H^0_{(2)}(\Sigma,L^p)$
extend to holomorphic sections of $L^p$ over $\overline\Sigma$, 
i.\,e.,  (see \cite[(6.2.17)]{mm})
\begin{equation}\label{e:bs1}
    H^0_{(2)}(\Sigma,L^p)\subset 
H^0\big(\overline\Sigma,L^p\big).
\end{equation}
 In particular, 
the dimension $d_p$ of $H^0_{(2)}(\Sigma,L^p)$ is finite.
 
We denote by respectively by $B_p(\LargerCdot,\LargerCdot)$ and 
by $B_p(\LargerCdot)$ 
the (Schwartz-)Bergman kernel and the Bergman kernel function
of the space $H^0_{(2)}(\Sigma,L^p)$, defined as follows: 
if $\{S_\ell^p\}_{\ell\geq1}$ is an orthonormal 
basis of $H^0_{(2)}(\Sigma,L^p)$, then 
\begin{equation}\label{e:BFS1}
B_p(x,y) = \sum_{\ell=1}^{d_p}S^p_\ell(x)\otimes(S_p(y))^{*} 
\quad\text{and}\quad 
B_p(x)=\sum_{\ell=1}^{d_p}|S^p_\ell(x)|_{h^p}^2\,.
\end{equation} 
Note that these are 
independent of the choice of basis (see 
\cite[(6.1.10)]{mm} or \cite[Lemma 3.1]{CM11}). 
Similarly, let $B_p^{\D^*}(x,y)$ and $B_p^{\D^*}(x)$ be 
the Bergman kernel and Bergman kernel function of
$\big(\D^*, \omega_{\D^*}, 
\C, h_{\D^{*}}^{p} \big)$, and 
$h_{\D^{*}}=\big|\!\log(|z|^2)\big| h_{0}$
with $h_{0}$ the flat Hermitian metric on the trivial line bundle 
$\C$.
 
 The main result of this paper is a weighted estimate in
 the $C^m$-norm near the punctures for the 
 global Bergman kernel $B_p$ compared to the Bergman kernel 
 $B_p^{\D^*}$ of the punctured disc, uniformly
in the tensor powers of the given bundle.
 
Note that for $k\in \N$, the $C^{k}$-norm at $x\in \Sigma$
is defined for $\sigma\in C^\infty(\Sigma, L^p)$ as 
 \begin{equation}   \label{eq:2.13c}\begin{split}
 &|\sigma |_{C^k(h^p)}(x)= \big(  |\sigma|_{h^p}
 +\big|\nabla^{p,\Sigma}\sigma
\big|_{h^p,\omega_{\Sigma}}+\ldots+\big|(\nabla^{p,\Sigma})^k
\sigma\big|_{h^p,\omega_{\Sigma}}\big)(x),
   \end{split}  \end{equation}
with $\nabla^{p,\Sigma}$ is the connection on
$(T\Sigma)^{\otimes\ell}\otimes L^p$ induced by the Levi-Civita 
connection on $(T\Sigma, \omega_{\Sigma})$ and the Chern connection
on $(L^{p},h^p)$,
and the pointwise norm $|\,\LargerCdot\,|_{h^p,\omega_{\Sigma}}$
is induced by $\omega_{\Sigma}$ and $h^{p}$. 
In the same way, for $f\in  C^\infty(\Sigma,\C)$, its $C^{k}$-norm 
$|f|_{C^{k}}(x)$ at $x\in \Sigma$ is defined by using 
the Levi-Civita connection on $(T\Sigma, \omega_{\Sigma})$.

Let us fix a point $a\in D$ and work 
in coordinates centered at $a$. 
Let $\mathfrak{e}$ be the holomorphic frame of $L$ near $a$ 
corresponding to the trivialization in the condition ($\alpha$).
Let  $0<r<e^{-1}$ and let $\D_{r}^*$ be the punctured disc of radius $r$ 
centered at $a$.
  By the assumptions ($\alpha$), ($\beta$), under our trivialization
  $\mathfrak{e}$ of $L$ on the coordinate $z$ on 
  $\D^*_{r}$, we have the identification of
  the geometric data
$(\Sigma,\omega_{\Sigma}, L,h)|_{\D^*_{r}}
= (\D^*,\omega_{\D^*}, \C, h_{\D^*})|_{\D^*_{r}}$.

  \begin{thm}\label{t6.1}
  Assume that $(\Sigma, \omega_{\Sigma}, L, h)$ fulfill conditions
  ($\alpha$) and ($\beta$). Let $a\in D$, and $0<r<e^{-1}$
  as above. 
    Then for any $k\in \N,  \ell>0$, $\alpha\geq 0$, there exists 
    a constant $C_{k, \ell,\alpha}$ such that for $p\gg1$, 
     we have on $\D_{r/2}^*\times \D_{r/2}^*$
     \begin{equation}\label{eq:6.6}
      \Big|{B}_p^{\D^*}(x,y)  - {B}_p(x,y)\Big|_{C^k(h^p)} 
            \leq C_{k,\ell,\alpha}p^{-\ell}
\big|\!\log(|x|^2)\big|^{-\alpha}\big|\!\log(|y|^2)\big|^{-\alpha}, 
     \end{equation}
with norms computed with help of (double copies of) 
$\omega_{\Sigma}$, $h$, and the associated Levi-Civita
and Chern connections on $\D_{r/2}^*\times \D_{r/2}^*$.
   \end{thm}
As a consequence of the proof of Theorem \ref{t6.1},
we obtain in Corollary \ref{crl6.2} uniform
estimates for the Bergman kernel $B_{p}(x,y)$
away from the diagonal. Another consequence of Theorem \ref{t6.1} is the
weighted diagonal expansion of the Bergman kernel:
\begin{thm}   \label{thm_MainThm}
  Assume that $(\Sigma, \omega_{\Sigma}, L, h)$ fulfill conditions
  ($\alpha$) and ($\beta$). 
  Then the following estimate holds: 
  for any $\ell, m\in \N$, and every $\delta>0$, 
  there exists a constant $C=C(\ell, m, \delta)$ such that
  for all $p\in\N^*$, and $z\in V_1\cup\ldots\cup V_N$
  with the local coordinate $z_{j}$,
   \begin{equation}    \label{eqn_MainThm}
     \Big | B_p - B_p^{\D^*}\Big |_{C^m} (z_{j})\leq Cp^{-\ell}
     \, \big|\!\log(|z_{j}|^2)\big|^{-\delta},    
   \end{equation}
  with norms computed with help of $\omega_{\Sigma}$ 
  and the associated Levi-Civita connection on $\D_{r/2}^*$.
    \end{thm}


   \begin{rmk}   \label{t0.2} 
Theorems \ref{t6.1} and \ref{thm_MainThm} admits a generalization to
orbifold Riemann surfaces.   
Assume that $\overline\Sigma$ is a compact orbifold Riemann
surface,  and the finite set 
$D=\{a_1,\ldots,a_N\}\subset\overline\Sigma$
does not meet the (orbifold) singular set of $\overline\Sigma$.
Assume moreover that $L$ is a holomorphic orbifold line bundle on 
$\overline\Sigma$. Let $\omega_{\Sigma}$ be an orbifold Hermitian 
form on $\Sigma$ and $h$ an orbifold Hermitian metric on $L$
in the sense of \cite[\S 5.4]{mm}. 
The proof of Theorems \ref{t6.1} and \ref{thm_MainThm}
can be modified to show: If conditions $(\alpha),  (\beta)$ 
hold in this context, then (\ref{eqn_MainThm}) holds. In fact,
the elliptic estimate \cite[(4.14)]{DLM06}
and the finite propagation speed of wave operators hold on 
orbifolds as observed by \cite[\S 6]{M05},
so the arguments
used in this paper go through for orbifolds to get the conclusion.
       \end{rmk}
       
By \cite[Theorems\,6.1.1,\, 6.2.3]{mm}, for any compact set 
$K\subset\Sigma$ we have the following expansion
  on $K$ in any $C^m$-topology (see Theorem \ref{T:bke}), 
  \begin{equation}    \label{e:bke}
   \frac{1}{p}  B_p(x)= \frac{1}{2\pi}
     +\sum_{j=1}^\infty \bb_j(x) p^{-j}\qquad\text{as }p\to +\infty.
   \end{equation}
By contrast, Theorem \ref{thm_MainThm} gives a precise description
of $B_p$ up to the punctures, in terms of the Bergman kernel function 
of the Poincar\'e metric on the
local model of the punctured unit disc in $\C$. Note that in the 
case of \textit{smooth} Hermitian metrics with
positive curvature the Bergman kernel can be localized and
its local model is the Euclidean
space endowed with a trivial bundle of positive curvature, 
see \cite[Sections 4.1.2--3]{mm}.
This kind of localization is inspired from the analytic 
localization technique of Bismut-Lebeau \cite{BL}
in local index theory.
For the problem at hand we have to overcome 
difficulties linked to the presence of \textit{singularities}. 
For this purpose we prove a new weighted elliptic estimate
near the punctures, which is uniform in the tensor powers $p$.
Then, combining with a revisited, singularity-centered localization principle
and a weighted Sobolev embedding, we can reach 
\eqref{eqn_MainThm} for some \textit{negative} $\delta$ however. 
Using the fact that square integrable holomorphic sections of $L^p$
vanish on $D$, we finally establish \eqref{eqn_MainThm} 
for any \textit{positive} $\delta$.

 From a study of the model Bergman kernel functions 
 $B_p^{\D^*}$ on the punctured unit disc, 
 we get the following ratio estimate as a corollary
 of Theorem \ref{thm_MainThm} and Corollary \ref{crl_fp}: 
  \begin{crl}    \label{crl_CrlIntro}
   Let $(\Sigma, \omega_{\Sigma}, L, h)$ be as in 
   Theorem \ref{thm_MainThm}. Then 
    \begin{equation}   \label{eqn_CrlIntro}
     \sup_{x\in \Sigma}B_p(x)=
     \sup_{\substack{x\in \Sigma\\\sigma\in H^0_{(2)}(\Sigma,L^p)
	 \setminus\{0\}}}
     \frac{|\sigma(x)|^2_{h^p}}{\|\sigma\|^2_{{\bL}^2}} 
      =  \Big( \frac{p}{2\pi}\Big)^{3/2}+ \mathcal{O}(p)
      \qquad \text{as }p\to +\infty. 
    \end{equation}
  \end{crl}
 
 ~
  
 This collection of results represents, to our knowledge, the first example 
 of a \textit{uniform} $\bL^{\infty}$ asymptotic description of 
 the Bergman kernel function of a singular polarization.
 This is of particular interest in arithmetic situations.
 Note that the work of Burgos Gil et al.\ \cite{BBK07, BKK05} 
 developed the arithmetic intersection theory
 for log-singular Hermitian metrics, showing in particular
 that Arakelov heights can be defined, and applied 
 successfully the theory for the Hilbert modular surfaces.
 Our results provide some possible applications 
 in this direction. For example, the classical arithmetic Hilbert-Samuel
 theorem \cite{GS92} for positive Hermitian line bundles
 is usually used to produce global integral sections with small 
 sup-norm;  a combination of the recent work \cite{bf} with 
 the distortion estimate of our Corollary \ref{crl_CrlIntro} 
 should give some
 interesting arithmetic consequence for cusp forms 
 on arithmetic surfaces and Hilbert modular surfaces.

Corollary \ref{crl_CrlIntro} is also quite striking from a K\"{a}hler 
geometry point of view, 
as the supremum of the Bergman kernel is equivalent to 
$\big(\frac{p}{2\pi}\big)^n$ on 
compact polarized manifolds of complex dimension $n$ 
(cf.\ Corollary \ref{C:bksup}), compared to the non-integer 
exponent $\frac32$ in \eqref{eqn_CrlIntro}. 
In K\"{a}hler geometry moreover, a central problem is 
the relation between the existence of special complete/singular metrics
and the stability of the pair $(X, D)$ where $D$ is
a smooth divisor of a compact K\"{a}hler manifold $X$; 
see e.g. the suggestions of \cite[\S 3.1.2]{sze} for the case of 
``asymptotically hyperbolic K\"ahler metrics'', 
which naturally generalize to higher dimensions the complete metrics 
$\omega_{\Sigma}$ studied here. 
In this respect, Theorem \ref{thm_MainThm} is an initial
step towards the application of Bergman kernels
to this problem and the first instance in
which the behavior of the Bergman kernel at infinity 
is fully understood. 
Moreover, the technique developed 
here can be extended to the higher dimensional situation; 
more precisely, in the case of Poincar{\'e} type Kähler metrics 
with reasonably fine asymptotics 
on complement of divisors, 
see the construction of \cite[\S1.1]{auv} and \cite[Theorem 4]{auv2}, 
we expect uniform elliptic estimates and localization techniques to 
transpose in a straightforward way
\footnote{As for the higher-dimensional local reference metric, 
which \textit{has to} be taken as a \textit{perturbation} of 
the product of the one-dimensional Poincar{\'e} cusp metric with 
some smooth metric on the divisor, its Bergman kernel can surely 
be properly understood as well, but with other approaches than 
the explicit description of this paper.}.

Besides, we remark that the Bergman kernel of smooth approximations of
a singular metric can change dramatically near the limit, cf. \cite{DMN16}. 
In the same circle of ideas,   
note that the Bergman kernel $B_p$ provides information 
on holomorphic sections of $L^p$ over the whole $\overline{\Sigma}$ 
that vanish at order $\geq 1$ at the punctures $a_1,\ldots,a_N$, 
whereas, for some fixed $\alpha\in (0,1)$ smooth metrics 
and \og partial Bergman kernels\fg~can be used to derive information 
on  holomorphic sections of $L^p$ on 
$\overline{\Sigma}$ with higher vanishing order 
$\geq \lfloor \alpha p\rfloor$ at the $a_j$'s. (see \cite{CM15}). 
With this in mind, one might roughly observe three different regimes
in the asymptotics of $B_p^{\D^*}$, hence of $B_p$ (Poincar{\'e} type case),
corresponding to sections with vanishing order at least 1 
and $\ll p^{1/2}$, $\sim p^{1/2}$, and $\gg p^{1/2}$, respectively; 
see Section \ref{S:bkepd}, in particular Figure \ref{fig1} 
\footnote{Conversely, one might imagine the parameter $\alpha$ 
mentioned above (smooth case) go to 0 at speed $1/p$, possibly 
together with the approximation 
of some Poincar{\'e} type metric by carefully chosen smooth metrics,  
to try and understand how the partial Bergman kernels approximate $B_p$;  
for sure, some subtleties will occur along this double limit process,  
of which the wild transition region observed in Figure \ref{fig1} might 
be an artefact.}. 
Note also that the behavior of the Bergman kernel on singular 
Riemann surfaces is relevant for the theory of quantum Hall 
effect \cite{LCCW} and attracted attention recently.

We give an important example where Theorem \ref{thm_MainThm}
applies. Let $\overline\Sigma$ be a compact Riemann surface 
of genus $g$ and consider a finite set 
$D=\{a_1,\ldots,a_N\}\subset\overline\Sigma$. 
We also denote by $D$ the divisor $\sum_{j=1}^Na_j$ and let 
$\mathscr{O}_{\overline{\Sigma}}(D)$ 
be the associated line bundle.
The following conditions are equivalent:
\\[3pt]\indent
(i) $\Sigma=\overline\Sigma\smallsetminus D$ admits a complete 
K\"ahler-Einstein metric $\omega_\Sigma$
with $\operatorname{Ric}_{\omega_\Sigma}=-\omega_\Sigma$,
\\[2pt]\indent
(ii) $2g-2+N>0$, 
\\[2pt]\indent
(iii) the universal cover of $\Sigma$ is the upper-half plane
$\mathbb{H}$,
\\[2pt]\indent
(iv) $L=K_{\overline\Sigma}\otimes
\mathscr{O}_{\overline\Sigma}(D)$ is ample.
\\[3pt]
This follows from the Uniformization Theorem \cite[Chapter IV]{fk} 
and the fact that the Euler characteristic of $\Sigma$ equals 
$\chi(\Sigma)=2-2g-N$ and the degree of $L$ is 
$2g-2+N=-\chi(\Sigma)$. 
If one of these equivalent conditions is satisfied,
the K\"ahler-Einstein metric 
$\omega_{\Sigma}$ is induced by the Poincar\'e metric on 
$\mathbb{H}$;
$(\Sigma,\omega_{\Sigma})$ and
the formal square root of $(L,h)$ satisfy conditions $(\alpha)$ 
and $(\beta)$, see Lemma \ref{L:ab}. 
Theorem \ref{thm_MainThm} hence applies to this context.
Let $\Gamma$ be the Fuchsian group associated with the above 
Riemann surface
$\Sigma$, that is, $\Sigma\cong\Gamma\backslash\mathbb{H}$.
Then $\Gamma$ is a geometrically
finite Fuchsian group of the first kind, without elliptic elements. 
Conversely, if $\Gamma$ is such a group, then 
$\Sigma:=\Gamma\backslash\mathbb{H}$
can by compactified by finitely many points $D=\{a_1,\ldots,a_N\}$
into a compact Riemann surface
$\overline\Sigma$ such that the equivalent conditions (i)-(iv) above
are fulfilled.
Let $\mathcal{S}^\Gamma_{2p}$ be the space of cusp forms 
(Spitzenformen) of weight $2p$ 
of $\Gamma$ endowed with the Petersson inner product. 
We can form the Bergman kernel
function of $\mathcal{S}^\Gamma_{p}$ as in \eqref{e:BFS1},
denoted by $B^\Gamma_{p}$. 
We deduce from Corollary \ref{crl_CrlIntro}:
\begin{crl}    \label{crl_CrlIntro1}
   Let $\Gamma\subset\operatorname{PSL}(2,\R)$ be 
   a geometrically finite Fuchsian group of the first kind 
   without elliptic elements. Let $B^\Gamma_{p}$ be the 
   Bergman kernel function of cusp forms of weight $2p$. 
   If $\Gamma$ is cocompact then uniformly on 
   $\Gamma\backslash\mathbb{H}$,
    \begin{equation}   \label{eqn_CrlIntro1}
     B^\Gamma_{p}(x)
      =  \frac{p}{\pi}+ \mathcal{O}(1), \qquad \text{as } p\to +\infty. 
    \end{equation}
  If $\Gamma$ is not cocompact then
    \begin{equation}   \label{eqn_CrlIntro2}
     \sup_{x\in\Gamma\backslash\mathbb{H}}B^\Gamma_{p}(x)
      =  \Big( \frac{p}{\pi}\Big)^{3/2}+ \mathcal{O}(p), 
      \qquad \text{as } p\to +\infty. 
    \end{equation}
  \end{crl}
   
Uniform estimates for 
$\sup_{x\in\Gamma\backslash\mathbb{H}}B^\Gamma_{p}(x)$ 
are relevant in arithmetic geometry and
were proved in various degrees of generality and sharpness in 
\cite{AbUll95,MiUll98,jk04,fjk}. In \cite{fjk} it is proved that in the
cofinite but non-cocompact case 
$\sup_{x\in\Gamma\backslash\mathbb{H}}B^\Gamma_{p}(x) 
= \mathcal{O}( p^{3/2})$
and the result is optimal, at least up to an additive term in the 
exponent of the form $-\varepsilon$ for any $\varepsilon> 0$. 
Estimate \eqref{eqn_CrlIntro2} gives the precise coefficient 
of the leading term $p^{3/2}$ and is sharp
(by killing the ``$\varepsilon$ from below'' from \cite{fjk}).
Estimate \eqref{eqn_CrlIntro1} is the consequence of the general 
expansion of the Bergman kernel on compact
manifolds \cite{Ti90, bou, Ca99, Z98}
(cf. also \cite{DLM06, mm} and Theorem \ref{T:bke}). 

It turns out that Corollary \ref{crl_CrlIntro1} can be formulated 
so as to underline a certain uniformity in $\Gamma$,  
in the same fashion as in \cite{fjk}:

  \begin{thm}   \label{t0.5}
  Let $\Gamma_{0}\subset \operatorname{PSL}(2,\R)$ 
  be a fixed Fuchsian subgroup of the first kind 
    without elliptic elements and 
  let  $\Gamma\subset \Gamma_{0}$ be any subgroup
of finite index. If $\Gamma_{0}$ is cocompact, then
\begin{equation}   \label{eqn_CrlIntro1a}
     B^\Gamma_{p}(x)
      =  \frac{p}{\pi}+ \mathcal{O}_{\Gamma_{0}}(1), 
      \qquad \text{as } p\to +\infty. 
    \end{equation}
  If $\Gamma_{0}$ is not cocompact then
    \begin{equation}   \label{eqn_CrlIntro2a}
     \sup_{x\in\Gamma\backslash\mathbb{H}}B^\Gamma_{p}(x)
      =  \Big( \frac{p}{\pi}\Big)^{3/2}
      + \mathcal{O}_{\Gamma_{0}}(p), 
      \qquad \text{as } p\to +\infty.
    \end{equation}
 Here the implied constants in 
 $\mathcal{O}_{\Gamma_{0}}(1)$, $\mathcal{O}_{\Gamma_{0}}(p)$
 depend solely on $\Gamma_{0}$.
\end{thm}

Note that (\ref{eqn_CrlIntro1a}) is a special case of a more general
result, which is implied in \cite[\S 6.1.2]{mm} and that we state
as Theorem \ref{t2.5} in Section \ref{S:prelim}.

We consider further an extension of 
Theorem \ref{t0.5} to the case when
the group $\Gamma_0$ has elliptic elements. 
Then the quotients
$\Gamma\backslash\mathbb{H}$ are in general orbifolds. 
By using the result of
Dai-Liu-Ma \cite[(5.25)]{DLM06} on the Bergman kernel
asymptotics on orbifolds and the orbifold version of Theorem
\ref{thm_MainThm} we obtain the following.
 \begin{thm}   \label{t0.6}
  Let $\Gamma_{0}\subset \operatorname{PSL}(2,\R)$ 
  be a fixed Fuchsian subgroup of the first kind.
  Let $\{x_{j}\}_{j=1}^q$  be the orbifold points of 
  $\Gamma_{0}\backslash\mathbb{H}$ and $U_{x_{j}}$ be a small 
  neighborhood of $x_{j}$ in $\Gamma_{0}\backslash\mathbb{H}$.
 Let  $\Gamma\subset \Gamma_{0}$ be any subgroup
of finite index and $\pi_{\Gamma}: \Gamma\backslash\mathbb{H}\to 
\Gamma_{0}\backslash\mathbb{H}$ be the natural projection.
If $\Gamma_{0}$ is cocompact, then
as $p\to +\infty$
\begin{equation}   \label{eq:0.7}
     B^\Gamma_{p}(x)
      =  \frac{p}{\pi}+ \mathcal{O}_{\Gamma_{0}}(1), 
\, \, \text{ uniformly   on } (\Gamma\backslash\mathbb{H}) 
\smallsetminus \bigcup_{j=1}^q \pi_{\Gamma}^{-1}(U_{x_{j}}).
    \end{equation}
On each $\pi_{\Gamma}^{-1}(U_{x_{j}})$ we have as $p\to +\infty$, 
    \begin{equation}   \label{eq:0.8}
     B^\Gamma_{p}(x)=  \Bigg(1+ 
 \sum_{\gamma\in \Gamma_{\!\!x_j^\Gamma} \smallsetminus\{1\}}
 \exp\Big(ip\theta_{\gamma}
 - p(1-e^{i\theta_\gamma })|z|^2\Big)
      \Bigg) 
      \frac{p}{\pi}+ \mathcal{O}_{\Gamma_{0}}(1), 
    \end{equation}
 where $x_{j}^\Gamma\in \pi_{\Gamma}^{-1}(x_{j})$ 
is in the same component of $\pi_{\Gamma}^{-1}(U_{x_{j}})$ 
as $x$, $e^{i\theta_{\gamma}}$ is the action of $\gamma$
on the fiber of $K_{\Gamma\backslash\mathbb{H}}$ at 
$x_{j}^\Gamma$, and $z=z(x)$ is the coordinate of $x$ 
 in normal coordinates $z$ centered at 
 $x_{j}^\Gamma$ in $\mathbb{H}$, and 
 $\Gamma_{y}=\{\gamma\in \Gamma : \gamma y = y\}$
 the stabilizer of $y$. 
    
    \noindent
    In particular, if $q_{0}= {\rm lcm}\{|\Gamma_{0,x_{j}}|: j=1,\ldots, 
    q\}$, $n_{\Gamma}= \max\{ |\Gamma_{y}| :
    y\in \pi_{\Gamma}^{-1}(x_{j}),  j=1,\ldots, 
    q\}$, then 
    \begin{equation}   \label{eq:0.9}
     \sup_{x\in\Gamma\backslash\mathbb{H}}B^\Gamma_{q_{0}p}(x)
      =  n_{\Gamma}\frac{q_{0}p}{\pi} 
      + \mathcal{O}_{\Gamma_{0}}(1).
   \end{equation}
  If $\Gamma_{0}$ is not cocompact then as $p\to +\infty$
    \begin{equation}   \label{eq:0.10}
     \sup_{x\in\Gamma\backslash\mathbb{H}}B^\Gamma_{p}(x)
      =  \Big( \frac{p}{\pi}\Big)^{3/2}
      + \mathcal{O}_{\Gamma_{0}}(p).
    \end{equation}
 Here again the implied constants in 
 $\mathcal{O}_{\Gamma_{0}}(1)$, $\mathcal{O}_{\Gamma_{0}}(p)$
 depend solely on $\Gamma_{0}$.
\end{thm}
Theorems \ref{t0.5}, \ref{t0.6} sharpen
in an optimal way the main result of \cite{fjk} that states
that 
\begin{align}\label{eq:0.12}
    \sup_{x\in\Gamma\backslash\mathbb{H}}B^\Gamma_{p}(x)
=\begin{cases}\mathcal{O}_{\Gamma_{0}}(p) \quad\quad \text{ if 
$\Gamma_{0}$ is cocompact},\\
\mathcal{O}_{\Gamma_{0}}(p^{3/2}) \quad  \text{ if 
$\Gamma_{0}$ is not cocompact}.
\end{cases}\end{align}
We obtain in this way the the precise leading terms in \eqref{eq:0.12}.

\smallskip

The results of this paper were announced in \cite{AMM16}.

This paper is organized as follows. In Section \ref{S:prelim} 
we recall the Bergman kernel expansion of complete K\"ahler 
manifolds
and introduce the functional space we need further.
In Section \ref{S:bkepd}, we study our model situation: 
the Bergman kernel on the punctured unit disc with Poincar\'{e} 
metric. In Section \ref{S:EE}, we establish the basic weighted 
elliptic estimate on the punctured unit disc with Poincar\'{e} 
metric uniformly with respect to the $p$-th power of the trivial 
line bundle with Poincar\'{e} metric. In Section \ref{S:sgl}, we 
develop the spectral gap properties of the Kodaira Laplacian
and give a rough uniform estimate of an approximation 
of the Bergman kernel.
In Section \ref{S:ibke}, by combining the finite propagation speed
of the wave operator and Section \ref{S:sgl},
we establish finally the main results
stated in the Introduction.
In the Appendix \ref{app}, we prove a technical result, 
Lemma \ref{lem_approx_psip}.

\smallskip
\noindent
\textbf{\emph{Acknowledgments.}}\
H.\ A.\ is thankful to the University of Cologne 
where this paper was partly written; he would also like to thank 
Michael Singer for inspiring conversations. 
G.~M.\ acknowledges support from 
Universit\'e Paris Diderot--Paris 7 (now Universit\'e de Paris)
where this paper was partly written and warmly thanks the project
\emph{Analyse Complexe et G\'eom\'etrie} for hospitality 
over many years.
\section{Preliminaries}\label{S:prelim}

In Section \ref{S:prelim1} we recall by following \cite{mm}
the asymptotics of the Bergman kernels
on complete manifolds and prove some results of independent
interest about this expansion on Riemann surfaces with 
locally constant curvature and also about its behavior 
with respect to coverings. In Section \ref{ss:FS} we introduce 
some functional and section spaces
that will be used throughout the paper.

\subsection{Expansion of Bergman kernels on complete manifolds}
\label{S:prelim1}
For a Hermitian
holomorphic line bundle $(L,h)$ on a complex manifold we denote 
by $R^L$ its Chern curvature and by $c_1(L,h)=\frac{i}{2\pi}R^L$
its Chern form.

Let $(M,\omega_M)$ be a complete K\"ahler manifold of dimension 
$n$ and $(L,h)$ be a Hermitian holomorphic line bundle on $M$
and $K_{M}$ be the canonical line bundle on $M$.
Then  the ${\bL}^2$-norm on $C^\infty_{0}(M,L^p)$, the space 
of smooth sections of $L^p$ with compact support,
is defined for any $s\in C^\infty_{0}(M,L^p)$ by
\begin{align}\label{eq:2.1}
\|s\|_{{\bL}^2}^2= \int_{M}|s(x)|_{h^p}^2 
\frac{\omega_{M}^n}{n!}\, \cdot
\end{align}
Let ${\bL}^2(M,L^p)$ be the ${\bL}^2$-completion of 
$(C^\infty_{0}(M,L^p),\|\LargerCdot \|_{{\bL}^2})$. 
We denote by $\langle\LargerCdot,\LargerCdot\rangle$ the inner 
product on ${\bL}^2(M,L^p)$ induced by this ${\bL}^2$-norm.  
Then the Bergman kernel function $B_{p}(x)\in C^\infty(M,\R)$
is still defined by (\ref{e:BFS1}) with $\{S_\ell^p\}_{\ell\geq1}$ 
an orthonormal basis of $H^0_{(2)}(M,L^p)$,
the space of ${\bL}^2$-holomorphic sections of $L^p$ on $M$
with respect to (\ref{eq:2.1}).
The Bergman kernel $B_{p}(x,y)$ is the smooth kernel of the 
orthonormal projection from 
$({\bL}^2(M,L^p),\|\LargerCdot \|_{{\bL}^2})$ 
onto $H^0_{(2)}(M,L^p)$. We have
\begin{align}\label{eq:2.2}
B_{p}(x,y)=\sum_{\ell\geq 1} S_{\ell}^p(x)\otimes (S_{\ell}^p(y))^*
\in L^p_{x}\otimes (L^{p}_{y})^*,\quad 
\text{ and }  B_{p}(x,x)= B_{p}(x).
\end{align}
Here $(S_{\ell}^p(y))^*\in (L^{p}_{y})^*$ is the metric dual 
of $S_{\ell}^p(y)$ with respect to $h^p$.

The Bergman kernel function \eqref{e:BFS1}
has the following variational
characterization (see \cite[Lemma 3.1]{CM11}):
\begin{equation}\label{e:Pvar}
B_p(x)=\max\Big\{|S(x)|^2_{h^p}:
\,S\in H^0_{(2)}(M,L^p),\;\|S\|_{{\bL}^2}=1\Big\}.
\end{equation}
We recall the expansion theorem for the Bergman kernel
on a complete manifold \cite[Theorem 6.1.1]{mm}. 
\begin{thm}\label{T:bke}
Let $(M,\omega_M)$ be a complete K\"ahler manifold of dimension 
$n$ and $(L,h)$ be a Hermitian holomorphic line bundle on $M$.
We assume there exist $\varepsilon>0$, $C>0$ such that 
$iR^L\geq\varepsilon\omega_M$ and 
$\ric_{\omega_M}\geq-C\omega_M$,
where $\ric_{\omega_M}=iR^{K_M^*}$ is the Ricci curvature of 
$\omega_M$.
Then there exist coefficients 
$\bb_j\in C^\infty(M)$,  $j\in\N$, 
such that for any compact set $K\subset M$, any $k,m\in\N$, 
there exists
 $C_ {k,m,K}>0$ such that for $p\in\N^*$,
\begin{equation} \label{ell1} 
\Big\|\frac{1}{p^n}B_p(x)
-\sum^k_{j=0}\bb_j(x)p^{-j}\Big\|_{C^m(K)}
\leqslant C_{k,m,K}\,p^{-k-1},
\end{equation} 
where
\begin{equation}\label{e:coeff}
\bb_0=\frac{c_1(L,h)^n}{\omega_{M}^n},\:\: 
\bb_1= \frac{\bb_0}{8\pi}\,(r_\omega-2\Delta_\omega\log\bb_0),
\end{equation}
and $r_\omega$, $\Delta_\omega$, are the scalar curvature, 
respectively 
the (positive) Laplacian, of the Riemannian metric associated to 
$\omega:=c_1(L,h)$.
\end{thm}
We write \eqref{ell1} shortly as 
\begin{align}\label{eq:2.5}
B_p(x)=\sum^k_{j=0}\bb_j(x)p^{n-j}+\mathcal{O}(p^{n-k-1}).
\end{align}
For compact or certain complete K\"ahler-Einstein manifolds 
the expansion was obtained by Tian \cite{Ti90} for $k=0$ and $m=2$.
For general $k$, $m$ and compact manifolds the 
existence of the expansion 
was first obtained in \cite{Ca99,Z98}.

The proof of \cite[Theorem 6.1.1]{mm} crucially relies on 
the following localization principle for Bergman kernels, that we use
in the proof of Corollary \ref{C:bke} below. Namely, 
as illustrated by the formulas for $\bb_0$ and $\bb_1$ in \eqref{e:coeff}, 
the asymptotics of $B_p(x)$ depend only on the
geometric data in any neighborhood of $x\in M$.
Hence, the Bergman kernel function asymptotics are the same
on two open sets (in two possibly different manifolds) 
over which the geometric data are isometric. 
\begin{thm}\label{T:BKloc}
Let $(M_1,\omega_{M_1})$, $(M_2,\omega_{M_2})$
 be  complete K\"ahler manifolds of dimension $n$ and
  $(L_1,h_1)\to M_1$, $(L_2,h_2)\to M_2$ be Hermitian 
  holomorphic line bundles.
We assume there exist $\varepsilon>0$, $C>0$ such that 
for $j=1,2$ we have $iR^{L_j}\geq\varepsilon\omega_{M_j}$ 
and $\ric_{\omega_{M_j}}\geq-C\omega_{M_j}$.
Assume moreover that there are open sets $U_j\subset M_j$, 
$j=1,2$, and biholomorphic isometries $\Phi:U_1\to U_2$, 
$\Psi:(L_1|_{U_1},h_1)\to \Phi^*((L_2|_{U_2},h_2))$, where
$\Psi$ is also a bundle isomorphism.
Let us denote by $B_{j,p}$ the Bergman kernel functions of
$H^0_{(2)}(M_j,L_j^p)$, $j=1,2$.
Then for any $k,m\in\N$ and any compact set $K\subset U_1$,
we have
\[B_{1,p}-B_{2,p}\circ\Phi=\mathcal{O}(p^{-k})\:\:
\text{in $C^m(K)$  as $p\to +\infty$}.\] 
In particular, if $\bb_{1,j}$ and $\bb_{2,j}$ denote the
 coefficients of the expansion \eqref{eq:2.5} of $B_{1,p}$ 
 and $B_{2,p}$, then $\bb_{1,j}=\bb_{2,j}\circ\Phi$ on
$U_1$ for all $j\in\N$.
\end{thm}
Note that in \cite{HsM11}, Hsiao and Marinescu got also a localization
principe on any non-compact manifold for the kernel for lower energy 
forms. 

We immediately obtain from Theorem \ref{T:bke} uniform sup-norm 
bounds for the Bergman kernel on compact subsets.
\begin{crl}\label{C:bksup}
Under the hypotheses of Theorem \ref{T:bke}, let $K\subset M$
be a compact subset such that
$iR^L=\omega_M$ on $K$. Then uniformly on $K$,
\begin{equation}\label{e:bksup}
B_p(x)=\Big(\frac{p}{2\pi}\Big)^n
+\mathcal{O}(p^{n-1}),\qquad \text{as } p\to +\infty.
\end{equation}
\end{crl}
In the case of dimension one and constant curvature we can
state the following.
\begin{crl}\label{C:bke}
Assume that $M$ in Theorem \ref{T:bke} is a Riemann surface and
there exists an open set $V$ such that $\omega_M$ has scalar 
curvature $-4$ and
$iR^L=\omega_{M}$ on $V$.
Then for any $k,m\in\N$ and any compact set $K\subset V$,
\begin{equation}\label{e:bke1}
     B_p(x)= \frac{1}{2\pi}p-\frac{1}{2\pi}
     +\mathcal{O}(p^{-k})\:\:\text{in $C^m(K)$  as $p\to +\infty$}.
\end{equation} 
\end{crl}
\noindent
\textit{Proof. --- } From \eqref{e:coeff}
follows that $\bb_0=\frac{1}{2\pi}$ and $\bb_1=-\frac{1}{2\pi}$
(note that $r_\omega=-8\pi$), thus
the task is to prove that the coefficients $\bb_j$ of the expansions 
\eqref{e:bke}, \eqref{eq:2.5} vanish on $V$
for $j\geq2$.
We divide the proof in three steps.

Firstly, it is easy to observe that $\bb_j$ are constant
functions on $V$ for all $j\in\N$.
Indeed, by \cite[Theorem 4.1.1]{mm} we know that 
$\bb_j$, $j\in\N$, are polynomials in the
curvatures $R^L$ and $R^{T^{(1,0)}M}$ and their derivatives.
On $V$ we have $iR^L=\omega_M$ and 
$iR^{T^{(1,0)}M}=-2\omega_M$.
Thus all the derivatives alluded to above vanish on $V$, hence 
$\bb_j$ are polynomials just in $R^L$ and $R^{T^{(1,0)}M}$, 
hence constant functions on $V$, for all $j\in\N$.

Secondly, we prove the assertion of the Corollary for a compact 
Riemann surface $\Sigma_{1}$ with genus  $g\geq 2$, such that
$\Sigma_{1}\sim \Gamma_{1}\backslash\mathbb{H}$,
with $\Gamma_{1}$ a cocompact Fuchsian group.
We endow $\Sigma_{1}$ with the metric $\omega_{\Sigma_{1}}$
induced from the Poincar\'e metric of $\mathbb{H}$
with scalar curvature $-4$. We consider the line bundle
$L_{1}= T^{*(1,0)}\Sigma_{1}= K_{\Sigma_{1}}$
endowed with the metric $h_{1}$ induced by 
$\omega_{\Sigma_{1}}$. Thus 
$iR^{L_{1}}= 2\omega_{\Sigma_{1}}$.
Let $B_{1,p}(x)$ be the Bergman kernel function of 
$H^0(\Sigma_1, L_1^p)$. By our observation above, 
the coefficients $\bb_{1,j}$ of the expansion \eqref{eq:2.5} 
are constant functions on $\Sigma_{1}$ for all $j\in\N$. Thus
\begin{align}\label{eq:2.7}
    B_{1,p}(x) = \sum_{j=0}^{k} \bb_{1,j} p^{1-j}
    + \mathcal{O}(p^{-k-1}).
\end{align}
By the Riemann-Roch theorem, for $p>1$,
\begin{align}\label{eq:2.9}
\int_{\Sigma_{1}}B_{1,p}(x)\omega_{\Sigma_{1}}
=\dim H^0(\Sigma_{1},L_{1}^p)
= \int_{\Sigma_{1}} \Big(p-\frac{1}{2}\Big)
c_{1}(K_{\Sigma_{1}}, h_{1}),
\end{align}
and $c_{1}(K_{\Sigma_{1}}, h_{1})
= \frac{1}{\pi}\omega_{\Sigma_{1}}$. 
By plugging the expansion \eqref{eq:2.7} into \eqref{eq:2.9},
identifying the coefficients of the powers of $p$ and taking
into account that $\bb_{1,j}$ are constants we get
\begin{align}\label{eq:2.10}
 \bb_{1,0}=\frac{1}{\pi}, \quad 
    \bb_{1,1}=- \frac{1}{2\pi},
    \quad \bb_{1,j}=0 \quad \text{for } j\geq 2.
\end{align}
Thirdly, we use the localization principle for Bergman kernels
formulated in Theorem \ref{T:BKloc}.
We now identify holomorphically and isometrically  ${L}_2$ 
on a neighborhood of $x\in V$ to $L_{1}$ on an open set of 
$\Sigma_{1}$. Indeed,  
by \cite[Theorem 2.5.17 and Corollary 2.5.18]{Wolf}, 
near $x$, the surface is locally isometric to 
the Poincar\'{e} upper half-plane $\mathbb{H}$,
and the holomorphic structure of the surface is determined
by the conformal structure fixed by the metric, thus 
we obtain 
a holomorphic and isometric identification $\Psi$ of
 a convex neighborhood $U$ of $x\in V$ to an open set of 
$\Sigma_{1}$. Then the curvature of the Chern connection on
 the line bundle ${L}^2\otimes \Psi^* K_{\Sigma_1}^{-1}$
 with the induced metric $h$ is zero on $U$.
 If $\sigma$ is a holomorphic frame of 
 ${L}^2\otimes \Psi^* K_{\Sigma_1}^{-1}$ on $U$, 
 this means that $\partial \overline{\partial} \log |\sigma|^2_h=0$,
 so there is a holomorphic function $f$ on $U$ such that
 $\log |\sigma|^2_h= 2 {\rm Im} f$ (which holds in any dimension. 
 For example, from the Poincar\'e Lemma, there exists 
 $g\in C^\infty(U,\R)$ such that 
 $(\overline{\partial}- \partial) \log |\sigma|^2_h	= i d g$, 
 this implies that $g+ i \log |\sigma|^2_h$ is holomorphic).
 Now $e^{-f} \sigma$ is a holomorphic frame 
 of ${L}^2\otimes \Psi^* K_{\Sigma_1}^{-1}$ such that
 $|e^{-f}\sigma|^2_h=1$ on $U$, 
 and this yields a holomorphic and isometric identification
 of ${L}^2$ to $\Psi^* K_{\Sigma_1}$.
 
 By Theorem \ref{T:BKloc}, we know the asymptotics of 
$B_{p}(x)$ is as same as 
of $B_{1,p}$, thus from (\ref{eq:2.7}) and (\ref{eq:2.10}), 
we get that \eqref{e:bke1} holds uniformly on $K\subset V$.
\cqfd
\\

Observe that if $(\Sigma, \omega_{\Sigma}, L, h)$ fulfill conditions 
($\alpha$) and ($\beta$),
the hypotheses of Corollary \ref{C:bke} are satisfied for 
$V=V_1\cup\ldots\cup V_N$
(note that the scalar curvature of the Poincar\'e metric
\eqref{eqn_omegaPcr} 
equals $-4$), thus
\eqref{e:bke1} holds on any compact set 
$K\subset V_1\cup\ldots\cup V_N$.

The following result is a direct consequence of the proof of
 \cite[Theorem 6.1.4]{mm}, and for completeness, 
we include the proof in Section \ref{S:ibke}.
\begin{thm}   \label{t2.5}
 Let $(M, \omega_{M}, L, h)$ be in Theorem \ref{T:bke}
 and assume moreover that $M$ is compact. Let $\pi_{1}(M)$
 be the fundamental group of $M$and $\widetilde{M}$ be the 
 universal covering of $M$. For
  any  subgroup $\Gamma\subset \pi_{1}(M)$
 with finite index, we define the Bergman kernel
     $B^\Gamma_{p}(x,y)$ on $\widetilde{M}/\Gamma$
     with the pull-back objects from 
     $\pi_{\Gamma}: \widetilde{M}/\Gamma\to M$.
 Then for any $k,m\in\N$, 
there exists
 $C_ {k,m}>0$ such that for any $\Gamma$ as above we have
 \begin{equation} \label{ell1a} 
\Big\|\frac{1}{p^n}B^\Gamma_p(x)
-\sum^k_{j=0}(\pi_{\Gamma}^*\bb_j)(x)p^{-j}
\Big\|_{C^m(\widetilde{M}/\Gamma)}
\leqslant C_{k,m}\,p^{-k-1},
\end{equation}     
where $\bb_j$ are the coefficients of the expansion \eqref{ell1} 
on $M$.
\end{thm}

\subsection{Functional spaces, section spaces}\label{ss:FS}
We define a few functional spaces, 
that will be much helpful in what follows. 

\textbf{(i)} $C^{0}(\D^*, \omega_{\D^*})$ is merely 
 the space of bounded continuous functions on $\D^*$, 
  endowed with the $\sup$ norm; notice that the reference to
  the metric, needed when considering bounds on derivatives,
  is superfluous here.

\textbf{(ii)} Let $U\subset\Sigma$ be an open set. The space
$C^k(U, \omega_{\Sigma})$ is defined as the set
of $C^k$ functions on $U$  bounded up to order $k$ on $U$ 
with respect to the metric $\omega_{\Sigma}$, 
and endowed with the natural norm:
\begin{equation*}
 C^k(U, \omega_{\Sigma})= \big\{f \in C^k(U): \, 
 \|f\|_{C^k(U, \omega_{\Sigma})}<\infty \big\}, 
                   \end{equation*}
                  with 
                   \begin{equation}   \label{eqn_Cmnorm}
  \|f\|_{C^k(U, \omega_{\Sigma})}= \sup_{x\in U} |f|_{C^k}(x),
  \quad |f|_{C^k}(x)=
   \big(|f|+|\nabla^{\Sigma}f|_{\omega_{\Sigma}}+\ldots
+|(\nabla^{\Sigma})^k f|_{\omega_{\Sigma}}\big)(x), 
                   \end{equation}
    $\nabla^{\Sigma}$ being the Levi-Civita connection attached 
    to $\omega_{\Sigma}$.

In the same vein, $C^k(U, \omega_{\Sigma}, L^p, h^p)$ is 
the space of $C^k$ sections of $L^p$ on $U$  
such that the following norm is bounded for 
$\sigma\in C^k(U, \omega_{\Sigma}, L^p, h^p)$:
 \begin{equation}   \label{eq:2.13}\begin{split}
 &|\sigma |_{C^k(h^p)}(x)= \big(  |\sigma|_{h^p}
 +\big|\nabla^{p,\Sigma}\sigma
\big|_{h^p,\omega_{\Sigma}}+\ldots+\big|(\nabla^{p,\Sigma})^k
\sigma\big|_{h^p,\omega_{\Sigma}}\big)(x),\\
 &\|\sigma\|_{C^k(U, \omega_{\Sigma})} :=
  \sup_{x\in U}  |\sigma |_{C^k(h^p)}(x) <\infty,
   \end{split}  \end{equation}
with $\nabla^{p,\Sigma}$ is the connection on
$(T\Sigma)^{\otimes\ell}\otimes L^p$ 
induced by the Levi-Civita connection associated with 
$\omega_{\Sigma}$ and the Chern connection relative to $h^p$. 
 
 \comment{
\textbf{(iii)} For $k\geq1$ we denote by $L^{2,k}
\big(\D^*, \omega_{\D^*}, \C, 
 \big|\!\log(|z|^2)\big|^p h_{0}\big)$ the Sobolev space 
    of sections over $\D^*$ of the trivial line bundle $\C$ 
      \text{endowed with the non-trivial Hermitian norm} 
$\big|\!\log(|z|^2)\big|^p h_{0}$
                  (the trivial Hermitian norm being $h_0$), 
                  which are ${\bL}^2$ up to order $k$, with respect to 
$\omega_{\D^*}$ and $\big|\!\log(|z|^2)\big|^p h_{0}$. 
                  In other words, elements of $L^{2,k}\big(\D^*, 
\omega_{\D^*}, \C, \big|\!\log(|z|^2)\big|^p h_{0}\big)$ 
                  are functions $f$ with $L^{2,k}_{\rm loc}$ regularity 
on $\D^*$, such that:
                  \begin{equation}   \label{eqn_H2pnorm}
       \|f\|_{L^{2,k}_p(\D^*)}^2
           :=\int_{\D^*} \big(|f|^2+|\nabla^{p}f|_{\omega_{\D^*}}^2
+\ldots+\big|(\nabla^{p})^kf\big|_{\omega_{\D^*}}^2\big)
\big|\!\log(|z|^2)\big|^p \omega_{\D^*}
          <\infty,
                  \end{equation}
                  with $\nabla^p$ the Chern connection relative to 
$\big|\!\log(|z|^2)\big|^p h_{0}$ and its extension with respect 
to $\omega_{\D^*}$ 
                  (to be made explicit below). 
For $k=0$ we denote $\|f\|_{L^{2,0}_p(\D^*)}$ by 
$\|f\|_{\bL^{2}_p(\D^*)}$                  .
   The space $L^{2,k}\big(\D^*, \omega_{\D^*}, \C,
\big|\!\log(|z|^2)\big|^p h_{0}\big)$ can moreover be 
characterized as the closure of the space of compactly supported 
functions on $\D^*$ under the norm 
$\|\LargerCdot\|_{L^{2,k}_p(\D^*)}$ defined in 
\eqref{eqn_H2pnorm}. 
 }                 %
                 
\textbf{(iii)} For $k\geq1$, the space
 $\bL^{2,k}\big(\Sigma,\omega_{\Sigma}, L^p, h^p\big)$
 is the Sobolev space of sections of the line bundle $L^p$ 
endowed with the Hermitian metric $h^p$ over $\Sigma$ 
that are ${\bL}^2$ up to order $k$, with respect to 
    $\omega_{\Sigma}$ and $h^p$. 
      This way, elements of 
      $\bL^{2,k}\big(\Sigma,\omega_{\Sigma}, L^p, h^p\big)$ 
                  are sections $\sigma$ of $L^p$ with 
$\bL^{2,k}_{\rm loc}$ regularity on $\Sigma$, 
                  such that:
                  \begin{equation}   \label{eqn_H2pSigmanorm}
                   \|\sigma\|_{\bL^{2,k}_{p}(h)}^2 :=\int_{\Sigma} 
     \big(|\sigma|^2_{h^p}+\big|\nabla^{p,\Sigma}\sigma
\big|_{h^p,\omega_{\Sigma}}^2+\ldots+\big|(\nabla^{p,\Sigma})^k
\sigma\big|_{h^p,\omega_{\Sigma}}^2\big)
                          \omega_{\Sigma}
                       <\infty.
                  \end{equation}
                  Alternatively, 
$\bL^{2,k}\big(\Sigma,\omega_{\Sigma}, L^p, h^p\big)$ is the 
$\|\LargerCdot\|_{\bL^{2,k}_p(h)}$-completion 
of the space of smooth and compactly supported 
sections of $L^p$ over $\Sigma$, 
with $\|\LargerCdot\|_{\bL^{2,k}_p(h)}^2$ defined in 
 \eqref{eqn_H2pSigmanorm}. 
For $k=0$ we simply denote $\|\LargerCdot\|_{\bL^{2,0}_p(h)}$ by
$\|\LargerCdot\|_{\bL^{2}_p(h)}$ and the corresponding inner 
product by $\langle\LargerCdot,\LargerCdot\rangle_p$.

When we apply this definition for the trivial line bundle $\C$ 
                  endowed with the non-trivial Hermitian metric
$\big|\!\log(|z|^2)\big|^p h_{0}$
                  (the trivial Hermitian metric being $h_0$), 
we get the space $\bL^{2,k}\big(\D^*, \omega_{\D^*}, \C, 
 \big|\!\log(|z|^2)\big|^p h_{0}\big)$ and the norm   
 $ \|\LargerCdot\|_{\bL^{2,k}_p(\D^*)}$.

  \textbf{(iv)}                 
 We also need in the localization procedure below
 some 
  \textit{weighted Sobolev spaces on $(\Sigma,\omega_{\Sigma})$ 
  {\rm(}resp.\ on a double copy of 
  $(\Sigma,\omega_{\Sigma})${\rm)}}. 
  We first define the weight function $\rho$ on $\Sigma$ 
  as a smooth function, 
                  equal to $1$ far from the punctures, 
to $\big|\!\log(|z_j|^2)\big|$ near the puncture $a_j$, 
and everywhere $\geq1$. 
    Let now $k\in \N$, and $q\geq 1$; 
    the weighted Sobolev space 
    ${\bL}^{q,k}_{\rm wtd}(\Sigma, \omega_{\Sigma})$ 
is defined as the space of ${\bL}^{q,k}_{\rm loc}$ functions $f$ on 
  $\Sigma$ such that: 
                   \begin{equation}\label{eqn_H2pw}
                    \|f\|^q_{{\bL}^{q,k}_{\rm wtd}}                                                                                             
    := \int_{\Sigma} \rho \big(|f|^q+\ldots 
    +|(\nabla^{\Sigma})^kf|^q_{\omega_{\Sigma}}\big) 
    \omega_{\Sigma} <\infty.  
                   \end{equation}
 Notice moreover that $\Sigma\times\Sigma$ is the 
    complement of a simple normal crossing divisor in a compact 
    K\"ahler manifold, namely $\Sigma\times\Sigma
    =\big(\overline{\Sigma}^2\big)\smallsetminus {\bf D}$ with 
    ${\bf D}=(D\times\Sigma)+(\Sigma\times D)$.  
This way, the natural product metric 
$\omega_{\Sigma\times\Sigma}$ is a K\"ahler metric of Poincar\'{e}
type on $\Sigma\times\Sigma$ (see e.g.\ \cite[Definition 0.1]{auv}). 
 Analogously, 
    ${\bL}^{q,k}_{\rm wtd}(\Sigma\times\Sigma, 
    \omega_{\Sigma\times\Sigma})$ 
                  is the space of ${\bL}^{q,k}_{\rm loc}$ functions $f$ 
		  on $\Sigma\times\Sigma$, 
                  endowed with the product metric 
$\omega_{\Sigma\times\Sigma}(x,y) = \omega_{\Sigma}(x)
+\omega_{\Sigma}(y)$
                  such that 
                   \begin{equation}\label{eqn_H2pw1}
                   \begin{split}
   \|&f\|^q_{{\bL}^{q,k}_{\rm wtd}} 
   := \int\limits_{(x,y) \in\Sigma\times\Sigma}\!\!\rho(x)\rho(y) 
\big(|f(x,y)|^q + \ldots +|(\nabla^{\Sigma\times\Sigma})^k
f(x,y)|^q_{\omega_{\Sigma\times\Sigma}}\big)       
 \omega_{\Sigma}(x)\omega_{\Sigma}(y)  
                                      \end{split}
                   \end{equation}  
      is finite.  
                  \begin{lem}\label{t2.4}
 a) We have 
 $\bL^{1,3}_{\rm wtd}(\Sigma,\omega_\Sigma)
 \longhookrightarrow C^0(\Sigma)$,
 i.\,e., there exists $c_0>0$ such that for all 
 $f\in \bL^{1,3}_{\rm wtd}(\Sigma,\omega_\Sigma)$ we have
  \begin{equation}\label{eq:2.18}
  \|f\|_{C^0(\Sigma,\omega_\Sigma)}
  \leq c_0  \|f\|_{ \bL^{1,3}_{\rm wtd}}.
    \end{equation}

\noindent
        b)    There are continuous embeddings 
                    \begin{equation}\label{eq:2.19}
                     \bL^{2,k}_{\rm wtd}(\Sigma\times\Sigma, 
		     \omega_{\Sigma\times\Sigma}) 
   \longhookrightarrow C^m(\Sigma\times\Sigma, 
   \omega_{\Sigma\times\Sigma})
                    \end{equation}
                   for all $k$, $m$ such that $k>m+2$. 
                  \end{lem}
 \prf. --- a) is from  \cite[\S 4.A and Lemme 4.5]{biq}. 
 For b), after noticing that the proof of \cite[Lemme 4.5]{biq} 
 remains valid close to the divisor 
 $\mathbf{D}\subset\overline{\Sigma}\times\overline{\Sigma}$ 
 but far from the crossings 
 $(a_{j_1},a_{j_2})\in \overline{\Sigma}\times\overline{\Sigma}$, 
 we work around one of these, just as in the proof of 
 \cite[Lemma 4.4]{auv}. 
 More precisely, we choose two small punctured discs 
 $\D_{r_1}^*$ and $\D_{r_2}^*$ around $a_{j_1}$ and $a_{j_2}$
 in each $\Sigma$ respectively, 
 and cover the product $\D_{r_1}^*\times\D_{r_2}^*$ 
 in $\Sigma\times\Sigma$ with help of (self-overlapping) 
 holomorphic polydiscs:
  \begin{equation*}   \label{eq:Phil1l2}
   \begin{aligned}
   \Phi_{\ell_1,\ell_2} : \D_{\epsilon} \times \D_{\epsilon} 
   &\longrightarrow  \D_{r_1}^*\times\D_{r_2}^* \\
     (u, v)   \,\,\,     &\longmapsto      
     \Big(e^{-2^{\ell_1}\tfrac{1+u}{1-u}}, 
  e^{-2^{\ell_2}\tfrac{1+v}{1-v}} \Big), 
   \end{aligned}
  \end{equation*}
 with $\ell_1$, $\ell_2 \geq 0$, and $0<\epsilon<1$ 
 fixed independently of $\ell_1$ and $\ell_2$. This way, 
 $\D_{r_1}^*\times\D_{r_2}^* \subset
 \bigcup_{\ell_1,\ell_2=0}^{\infty}\Phi_{\ell_1,\ell_2}
 (\D_{\epsilon} \times \D_{\epsilon})$, 
 and we can even assume that 
 $\D_{r_1}^*\times\D_{r_2}^* \subset 
 \bigcup_{\ell_1,\ell_2=0}^{\infty}\Phi_{\ell_1,\ell_2}
 (\D_{\epsilon/2} \times \D_{\epsilon/2})$.  
 Moreover, for \textit{any} $(\ell_1,\ell_2)$, 
  \begin{equation}   \label{eq:varpi}
   (\Phi_{\ell_1,\ell_2})^*\omega_{\Sigma\times\Sigma} 
   = \frac{idu\wedge d\overline{u}}{(1-|u|^2)^2}
   +\frac{idv\wedge d\overline{v}}{(1-|v|^2)^2} =: \varpi,
  \end{equation}
 which does not depend on $(\ell_1,\ell_2)$. 
 On the other hand, 
 $(\Phi_{\ell_1,\ell_2})^*\rho(x) 
 = 2^{\ell_1+1}\big|\frac{1+u}{1-u}\big|$, which is of size 
 $2^{\ell_1+1}$ for $|u|\leq \epsilon$, 
 with derivatives (at every order) of the same size for $\varpi$, 
 and similarly for $(\Phi_{\ell_1,\ell_2})^*\rho(y)$
 with $2^{\ell_2+1}$. 
 
 Set $U = \D_{r_1}\times\D_{r_2}\subset 
 \overline{\Sigma}\times\overline{\Sigma}$, 
 so that $\D_{r_1}^*\times\D_{r_2}^*=U\smallsetminus\mathbf{D}$.   
 Take $w\in {\bL}^{q,k}_{\rm wtd}(\Sigma\times\Sigma,
 \omega_{\Sigma\times\Sigma})$, $q\geq 1$, $k\geq 0$, 
 and pick $m\geq 0$, $m<k-\frac{2}{q}$, so that 
 $w\in C^m(U\smallsetminus\mathbf{D})$; 
 what precedes thus yields: 
      \begin{equation}\label{eq:2.20}
      \begin{split}
        \|w\|_{C^m(U\smallsetminus \mathbf{D})}^q &
\leq \sup_{\ell_1,\ell_2 \geq0} 
\|(\Phi_{\ell_1,\ell_2})^*w\|_{ C^m(\D_{\epsilon/2} 
\times \D_{\epsilon/2},\varpi)}^q\\
&\leq \sum_{\ell_1,\ell_2 =0}^{\infty} 
\|(\Phi_{\ell_1,\ell_2})^*w\|_{C^m(\D_{\epsilon/2}\times
\D_{\epsilon/2},\varpi)}^q\\
&= \sum_{\ell_1,\ell_2=0}^\infty
\frac{1}{2^{\ell_1+\ell_2+2}}2^{\ell_1+\ell_2+2}
     \|(\Phi_{\ell_1,\ell_2})^*w\|_{C^m(\D_{\epsilon/2} \times
     \D_{\epsilon/2},\varpi)}^q     \\
        &\simeq \sum_{\ell_1,\ell_2=0}^{\infty}
\frac{1}{2^{\ell_1+\ell_2+2}}\big\|(\Phi_{\ell_1,\ell_2})^*
       \big(\rho(x)^{\frac1q}\rho(y)^{\frac1q}w\big)
\big\|_{C^m(\D_{\epsilon/2} \times \D_{\epsilon/2},\varpi)}^q\\
  &\leq c\sum_{\ell_1,\ell_2=0}^{\infty}\frac{1}{2^{\ell_1+\ell_2}}
    \big\|(\Phi_{\ell_1,\ell_2})^*
\big(\rho(x)^{\frac1q}\rho(y)^{\frac1q}w\big)
\big\|_{{\bL}^{q,k}(\D_{\epsilon} \times \D_{\epsilon},\varpi)}^q
      \end{split}
\end{equation}
 by the \textit{fixed} usual Sobolev embedding (or, more exactly,
 continuous restriction)
 \[{\bL}^{q,k}(\D_{\epsilon} \times \D_{\epsilon},\varpi)
  \longhookrightarrow 
  C^m(\D_{\epsilon/2} \times \D_{\epsilon/2},\varpi)\] 
 applied to \textit{all} the 
 $(\Phi_{\ell_1,\ell_2})^*\big(\rho(x)^{\frac1q}\rho(y)^{\frac1q}
 w\big)$. 
 Now, observe that our choices provide
  \begin{multline}\label{eq:2.22}
   \sum_{\ell_1,\ell_2=0}^{\infty}\frac{1}{2^{\ell_1+\ell_2}}
\big\|(\Phi_{\ell_1,\ell_2})^*\big(\rho(x)^{\frac1q}
\rho(y)^{\frac1q}w\big)\big\|_{{\bL}^{q,k}
(\D_{\epsilon} \times \D_{\epsilon},\varpi)}^q\\
        \leq C(q) \big\|\big(\rho(x)^{\frac1q}\rho(y)^{\frac1q}
w\big)\big\|_{{\bL}^{q,k}(\Sigma\times\Sigma,
\omega_{\Sigma\times\Sigma})}^q, 
  \end{multline}
 since the number of self-overlaps of 
 $\Phi_{\ell_1,\ell_2}(\D_{\epsilon} \times \D_{\epsilon})$ 
 is of order $2^{\ell_1+\ell_2}$, 
 and since each $\Phi_{\ell_1,\ell_2}(\D_{\epsilon} \times 
 \D_{\epsilon})$ overlaps 
 only a finite number of other
 $\Phi_{\ell_1',\ell_2'}(\D_{\epsilon} \times \D_{\epsilon})$, 
 this number being bounded independently of $\ell_1$ and $\ell_2$. 
 Hence 
  \begin{equation*}
   \|w\|_{C^m(U\smallsetminus \mathbf{D})}^q \leq   
   C\big\|\big(\rho(x)^{\frac1q}\rho(y)^{\frac1q}w\big)
   \big\|_{{\bL}^{q,k}(\Sigma\times\Sigma,\omega_{\Sigma\times
   \Sigma})}^q
   \simeq C\| w \|_{{\bL}^{q,k}_{\rm wtd}}^q, 
  \end{equation*}
 and one concludes by specializing to $q=2$, 
 and gathering such estimates around the crossings 
 $(a_{j_1},a_{j_2})$ with analogous estimates along the divisor 
 $\mathbf{D}\subset\overline\Sigma\times\overline\Sigma$ 
 and far from the crossings, 
 and estimates far from the divisor. \cqfd
           
\section{Bergman kernels on the punctured unit disc}\label{S:bkepd} 
 In this section we give a detailed description of the Bergman kernel 
 on the punctured unit disc.
 We first obtain an explicit formula in \S \ref{ss:bkf} and then in 
 \S \ref{ss:bka} we get precise
 asymptotics near the puncture by using a natural rescaling.
 
 \subsection{Expression of the Bergman kernels on the punctured 
 unit disc}\label{ss:bkf}
 Let $p\in \N^*$  and let
 \begin{align}\label{eq:3.1}
H_{(2)}^p(\D^*) := H^0_{(2)}\big(\D^*, \omega_{\D^*}, \C, 
 \big|\!\log(|z|^2)\big|^p h_{0}\big),
\end{align}
be the space of holomorphic functions $S$ on $\D^*$ with finite
 ${\bL}^2$-norm defined in Section \ref{ss:FS} {\bf (iii)} for $k=0$. 
 The purpose here is to study of the Bergman kernel 
 of $H_{(2)}^p(\D^*)$, as $p\to \infty$. 
 
 \begin{lem} \label{t3.1} For $p\geq 2$, the set
     \begin{equation}\label{eq:onb}
     \Big\{\Big(\dfrac{\ell^{p-1}}{2\pi (p-2)!}\Big)^{1/2}z^\ell:
     \ell\in\N,\,\ell\geq 1\Big\}
      \end{equation}
forms an orthonormal basis of $H_{(2)}^p(\D^*) $.
\end{lem}
\textit{Proof. --- }  
Let $H^0(\D,\C)$ be the space of holomorphic functions on $\D$.
By  \cite[(6.2.17)]{mm}, we know
    \begin{align}\label{eq:3.2}
H_{(2)}^p(\D^*) \subset H^0(\D,\C).
\end{align}
Note that for $p\geq 2$, $\ell\geq 1$,
    \begin{equation}\label{eq:3.4}
    \begin{split}
   \int_{\D^*} \big|\!\log(|z|^2)\big|^p \omega_{\D^*} 
   &= \int_{\D^*} \big|\!\log(|z|^2)\big|^{p-2} 
   \frac{idz\wedge d\overline{z}}{|z|^2} \\
    &= \int_{\mathbb{S}^1 }2^{p-1} \, d\theta
       \int_0^{1}|\log r|^{p-2}\frac{dr}{r} 
                                                       = \infty, 
     \end{split}                                                  
  \end{equation}
 and 
  \begin{equation}\label{eq:3.6}
  \begin{split}
   \int_{\D^*} |z^{\ell}|^2\big|\!\log(|z|^2)\big|^p \omega_{\D^*} 
   &= \int_{\D^*} \big|\!\log(|z|^2)\big|^{p-2}  |z|^{2\ell}
   \frac{idz\wedge d\overline{z}}{|z|^2} \\
    &= \int_{\mathbb{S}^1 }2^{p-1}\, d\theta \int_0^{1}r^{2\ell-1}
    |\log r|^{p-2}\, dr   \\
     &= 2^p\pi\cdot (2\ell)^{1-p}\cdot\Gamma(p-1) 
            = \frac{2\pi(p-2)!}{\ell^{p-1}}<\infty.
  \end{split}
  \end{equation}
  By (\ref{eq:3.4}), (\ref{eq:3.6}) and 
  the circle invariance of $\omega_{\D^*}$ 
 and $\big|\!\log(|z|^2)\big|^p h_{0}$, the set \eqref{eq:onb} 
 forms an orthonormal basis of $H_{(2)}^p(\D^*) $.
\cqfd

 \comment{
 As it turns out, for fixed $p\geq 2$, 
 such a Bergman kernel admits a rather explicit expression. 
 Indeed, up to normalization, the family $(z^{\ell})_{\ell\geq1}$
 is a Hilbert base of the space of holomorphic functions 
 $H_{(2)}^p(\D^*)$; 
 this can be seen using the circle invariances of $\omega_{\D^*}$ 
 and $\big|\!\log(|z|^2)\big|^p h_{0}$, 
 and the elementary calculations
  as soon as $p\geq 2$ and $\ell\geq1$;   
  }

 \begin{rmk}\label{t3.2}
Notice that a similar computation shows that 
 the elements of $H_{(2)}^0(\Sigma, L^p)$ are, for $p\geq 2$, 
 exactly the sections of $L^p$ over the whole $\overline{\Sigma}$
 vanishing on the puncture divisor $D=\{a_1,\ldots,a_N\}$. 
\end{rmk}

 Back to $\D^*$ and according to Lemma \ref{t3.1},
 the Bergman kernel of $H_{(2)}^p(\D^*)$, for any $p \geq 2$,
 is thus 
  \begin{equation}\label{eq:3.7}
   B_p^{\D^*}(x,y) =  
          \frac{
	  \big|\!\log(|y|^2)\big|^{p}}{2\pi(p-2)!}
     \sum_{\ell=1}^{\infty}\ell^{p-1} x^{\ell}\overline{y}^{\ell}.
  \end{equation}
 Here the metric dual of the canonical section $1$
 with respect to $h_{0}$ is identified to $1$, hence the metric 
 dual of $1$ with respect to $\big|\!\log(|z|^2)\big|^p h_{0}$
 is $1^*(z) = \big|\!\log(|z|^2)\big|^p 1$.
 Specializing to the diagonal, 
 we get in particular the Bergman kernel function of
 $H_{(2)}^p(\D^*)$ for all $p \geq 2$,
  \begin{equation}   \label{eqn_BpD}
   B_p^{\D^*}(z) =  \frac{\big|\!\log(|z|^2)\big|^{p}}{2\pi(p-2)!}
   \sum_{\ell=1}^{\infty}\ell^{p-1} |z|^{2\ell}.
  \end{equation}

 ~
 
 This readily provides the behavior of $B_p^{\D^*}$ far from 
 $0\in \D$.
 \begin{prop}\label{t3.3}
For any  $0<a<1$ and any $m\geq 0$,
 there exists $c=c(a)>0$ such that 
  \begin{align}\label{eq:3.10}
   \Big\|B_p^{\D^*}(z) - \frac{p-1}{2\pi}
   \Big\|_{C^m(\{a\leq |z|<1\}, \omega_{\D^*})} 
   = \mathcal{O}(e^{-cp})
   \:\: \text{ as } p\to +\infty\,.
  \end{align} 
  More generally,  for any  $0<a<1$ and $0<\gamma<\frac{1}{2}$, 
  there exists $c = c(a,\gamma)>0$ such that
   \begin{align}\label{eq:3.10a}
   \Big\|B_p^{\D^*}(z) - \frac{p-1}{2\pi}
   \Big\|_{C^m(\{ae^{-p^\gamma}\leq |z|<1\}, \omega_{\D^*})} 
   = \mathcal{O}\big(e^{-cp^{1-2\gamma}}\big)
   \:\: \text{ as } p\to +\infty\,,
  \end{align} 
  and if $b\in (0,1)$, 
   \begin{align}\label{eq:3.10b}
   \Big\|B_p^{\D^*}(z) - \frac{p-1}{2\pi}
   \Big\|_{C^m(\{b \, e^{-p^{1/2}(\log p)^{-1}}\leq |z|\leq 
a\, e^{-p^{\gamma}}\}, \omega_{\D^*})} 
   = \mathcal{O}\big(p^{-\infty}\big)
   \:\: \text{ as } p\to +\infty\,.
  \end{align} 
\end{prop}
\textit{Proof. --- } 
 Let us recall the celebrated formula from complex analysis: 
  \begin{equation*}
   \frac{1}{\sin^2 w} = \sum_{k\in \Z} \frac{1}{(w-k\pi)^2} 
   \qquad\text{on }\C\smallsetminus\pi\Z. 
  \end{equation*}
  Thus for $t>0$,
  \begin{align}\label{eq:3.8}
\sum_{k\in \Z} \frac{1}{(2ki \pi+t)^2} = (e^{t/2}- e^{-t/2})^{-2}
= \sum_{\ell=1}^{\infty} \ell e^{- \ell t}.
\end{align}
 Combining (\ref{eqn_BpD}), (\ref{eq:3.8}) with an easy
  induction on $p\geq2$, one gets the identity
  \begin{equation}\label{eq:3.9}
  \begin{split}
   B_p^{\D^*}(z) &= \frac{(p-1)}{2\pi}\sum_{k\in \Z} 
   \frac{\big|\!\log(|z|^2)\big|^{p}}{(2ik\pi+|\!\log(|z|^2)|)^p} \\
                 &= \frac{(p-1)}{2\pi}\Bigg(1+\sum_{k\in \Z, \,k\neq 0} 
\frac{\big|\!\log(|z|^2)\big|^{p}}{(2ik\pi+|\!\log(|z|^2)|)^p}\Bigg).
\end{split}
  \end{equation}
  To obtain (\ref{eq:3.10}) for $m=0$,  from (\ref{eq:3.9}),
  for $p\geq 2$, we use 
  \begin{align}\label{eq:3.11}\begin{split}
&\sum_{k\in \Z, \,k\neq 0} 
\frac{\big|\!\log(|z|^2)\big|^{p}}{|2ik\pi+|\!\log(|z|^2)||^p}\\
&\hspace{10mm} < 2 \left( 1+ \frac{(2\pi)^2}{\big|\!\log(|z|^2)\big|^2}
\right)^{-\frac{p-2}{2}}
\sum_{k=1}^{\infty}
\frac{\big|\!\log(|z|^2)\big|^{2}}{(2k\pi)^2+|\!\log(|z|^2)|^2}
\quad \text{for } 0<|z|\leq e^{-1/2},\\
&\sum_{k\in \Z, \,k\neq 0} 
\frac{\big|\!\log(|z|^2)\big|^{p}}{|2ik\pi+|\!\log(|z|^2)||^p}
\leq  2 \big|\log(|z|^2)\big|^{p}\sum_{k=1}^{\infty}(2k\pi)^{-p}
\quad  \text{for }\, e^{-1/2}\leq |z| <1.
\end{split}\end{align}
For $m\geq 1$, by considering separately $a<|z|\leq e^{-1/2}$
and $e^{-1/2}\leq |z| <1$ as above, we get  again
 (\ref{eq:3.10}) from (\ref{eq:3.9}).
 
 For $ae^{-p^{\gamma}}\leq |z|\leq e^{-1/2}$, from
  \begin{align*}
	  \begin{split}
      \log \left( 1+ \frac{(2\pi)^2}{\big|\!\log(|z|^2)\big|^2}
\right)\geq  \frac{C}{\! |\log(|z|^2)|^2 }
\geq  C p^{-2\gamma},\\
\sum_{k=1}^{\infty} 
\frac{\big|\!\log(|z|^2)\big|^{2}}{(2k\pi)^2+|\!\log(|z|^2)|^2}
\leq \sum_{k=1}^{\infty} 
\frac{Cp^{2\gamma}}{(2k\pi)^2+ 1},
\end{split}\end{align*}
and (\ref{eq:3.11}), we get
also (\ref{eq:3.10a}).

Note that $\log(1+x)\geq  x - \frac{x^{2}}{2}$ for $x\in [0,1]$.
From this we know that for $b e^{-p^{1/2}(\log p)^{-1}}\leq |z|\leq 
ae^{-p^{\gamma}}$, the analogue of the above equation is
 \begin{align*}
	 \begin{split}
        \log \left( 1+ \frac{(2\pi)^2}{\big|\!\log(|z|^2)\big|^2}
\right)\geq  \frac{C}{\! |\log(|z|^2)|^2 }
\geq  \frac{C}{p}(\log p)^{2},\\
\sum_{k=1}^{\infty} 
\frac{\big|\!\log(|z|^2)\big|^{2}}{(2k\pi)^2+|\!\log(|z|^2)|^2}
\leq \sum_{k=1}^{\infty} 
\frac{C\, p\, (\log p)^{-2}}{(2k\pi)^2+ 1}.
\end{split}\end{align*}
From the above equation and  (\ref{eq:3.11}) we get (\ref{eq:3.10b}).
\cqfd  \\

 Observe that the expected behavior for an Einstein metric
 of scalar curvature $-4$ such as $\omega_{\D^*}$, 
 at least on compact subsets of $\D^*$, according to 
 \eqref{e:bke},  Theorem \ref{T:bke} and Corollary \ref{C:bke},
 is $B_p^{\D^*}(z) - \frac{p-1}{2\pi}=\mathcal{O}(p^{-\infty})$. 
 From our explicit description of $B_p^{\D^*}$, 
 we hence benefit an improvement, namely 
 the existence of regions of \emph{exponential decay} of the remainder, 
 such regions extending
 \textit{up to the exterior boundary} $\partial\D$ of $\D^*$ on one side, 
 and \textit{exponentially close to the singularity $0\in \D$} 
 on the other side; we also get a first estimate on how close 
 to 0 the usual $\mathcal{O}(p^{-\infty})$ remains valid.

 \subsection{Asymptotics of the density functions near
 the puncture}\label{ss:bka}
 We are also interested in a global description of 
 $B_p^{\D^*}$ \textit{up to the singularity} $0\in \D$, 
 especially in the geometric context of Theorem \ref{thm_MainThm}, 
 and such a description requires another angle of attack. 
 Let us simplify notations: for $p\in \N^*$, set
 \begin{align}\label{eq:3.14}\begin{split}
 b_p(y) &= \frac{\big|\!\log y\big|^{p+1}}{2\pi(p-1)!}
   \sum_{\ell=1}^{\infty}\ell^{p} y^{\ell} \qquad\text{for }
   y\in (0,1),\\
   \varphi(\xi )&=e \, \xi\,  |\!\log\xi | \qquad \text{for } 
    \xi\in(0,1),\\
\nu(p)&= (2\pi p)^{1/2} p^p e^{-p}( p!)^{-1} -1.
\end{split}\end{align}
Note that by Stirling's formula and (\ref{eq:3.14}), 
\begin{align}\label{eq:3.15}
\nu(p)= \mathcal{O}(p^{-1}) \quad\text{ as } p\to +\infty.
\end{align}
By (\ref{eqn_BpD}) and (\ref{eq:3.14}), we have
\begin{align}\label{eq:3.13}
B_{p+1}^{\D^*}(z) = b_p(|z|^2) \quad \text{for } z\in \D^*.
\end{align}

 Motivated by 
 the observation that at fixed $y$, the index of 
 the largest term of the sum 
 $\sum_{\ell=1}^{\infty}\ell^{p} y^{\ell}$
 is determined by $y^{1/p}$, 
 we further proceed to the change of variable $x=y^{1/p}$, 
 and focus on the function $f_p:(0,1)\to \R$ given by:
  \begin{equation}   \label{eqn_fp}\begin{split}
   f_p(x)&:=b_p(x^p) =\frac{\big|\!\log(x^p)\big|^{p+1}}{2\pi(p-1)!}
  \sum_{\ell=1}^{\infty}\ell^{p} x^{p\ell}\\
   &= \frac{p^{p+2}e^{-p}}{2\pi p!}|\!\log x| 
   \sum_{\ell=1}^{\infty}\big(e x^{\ell}\big|\!\log(x^{\ell})
   \big|\big)^p\\
   &=\Big(\frac{p}{2\pi}\Big)^{3/2}
   \big(1+\nu(p)\big)|\!\log x| \sum_{\ell=1}^{\infty}
   \big(\varphi(x^{\ell})\big)^p.   
\end{split}  \end{equation}
 
 The smooth function $\varphi$ maps $(0,1)$ to $(0,1]$, with 
 $\varphi(\xi)=1$ iff $\xi=e^{-1}$. 
 Thus, for $\ell$ fixed, $x\mapsto \big( \varphi(x^{\ell})\big)^p$ 
 heuristically converges to a thinner and thinner 
 Gaussian-shaped bump of height 1 centered 
 at $e^{-1/\ell}$, 
 and $|\!\log x| \sum_{\ell=1}^{\infty}\big(\varphi(x^{\ell})\big)^p$ 
 can thus be thought of as a series of these bumps centered at 
 $e^{-1}$, $e^{-1/2}$, $e^{-1/3}$, 
 of respective heights $1$, $\frac{1}{2}$, $\frac{1}{3}$ 
 (because of the factor $|\!\log x|$) and so on;  
 this actually holds for $x$ in \og low regime
 \fg~($x\leq e^{-p^{-\delta}}$, $\delta>1/2$, say), 
 the \og tail\fg~($x\geq e^{-p^{-\delta}}$, $0<\delta<1/2$) 
 consisting in an agglomeration 
 of such bumps mixing up with one another to follow an almost 
 constant behavior near $x=1$: 
 see Figure \ref{fig1} below. 
   \smallskip  
   \begin{figure}[!ht] 
   \begin{centering}
    \includegraphics[height=6cm]{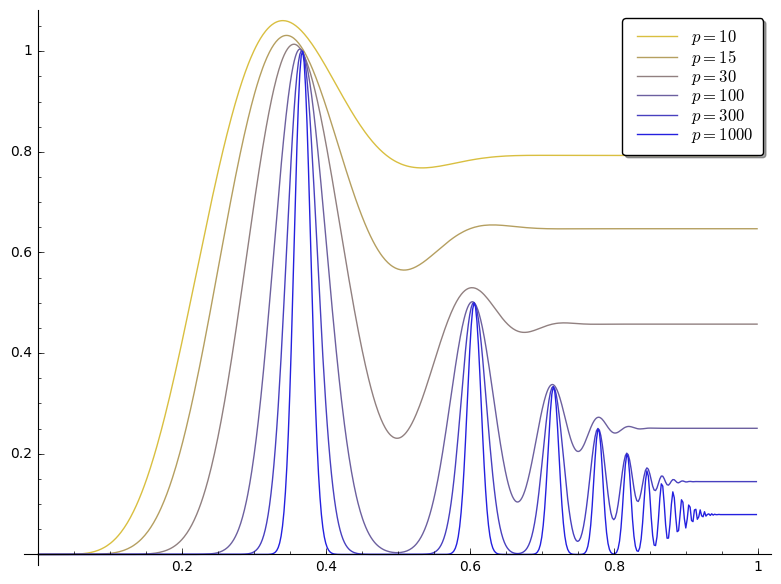}
   \caption{The scaled functions $\big(\frac{2\pi}{p}\big)^{3/2}f_p$ 
   on $(0,1)$}
      \label{fig1}
   \end{centering}
   \end{figure} 
   
  We develop in the following lines some elementary analysis that 
  justifies these heuristic considerations. 
  First, set
   \begin{align}   \label{eqn_psi}\begin{split}
 \psi_p(\zeta) &= \big(\varphi({e}^{-\zeta})\big)^p 
  = {e}^{p(1-\zeta+\log\zeta)}\quad \text{for  } \zeta>0,\\
    \mathcal{G}_0(\eta) &= {e}^{-\eta ^2/2},
    \qquad   \mathcal{G}_{1}(\eta) = \eta ^3 {e}^{-\eta ^2/2}  
    \quad \text{for  } \eta \in \R.
 \end{split} \end{align}
  We prove in the appendix \ref{app} the following estimate, 
  linking $\psi_p$ to the Gaussian-type functions $\mathcal{G}_0$ 
  and $\mathcal{G}_1$: 
   \begin{lem}    \label{lem_approx_psip}
    There exists $C>0$ such that for all $\zeta>0$ and all 
    $p\geq 1$, 
     \begin{equation*}
      \Big|\psi_p(\zeta) - \mathcal{G}_0\big(\sqrt{p}(1-\zeta)\big)
      +  \frac{1}{3\sqrt{p}}\mathcal{G}_1\big(\sqrt{p}(1-\zeta)\big)
      \Big| 
       \leq \frac{C}{p(1+p(1-\zeta)^2)}.
     \end{equation*}
   \end{lem}
   For $p\geq 1$ and $x\in (0,1)$, set 
     \begin{equation}\label{eq:3.17}
      \mathbb{G}_p(x)=|\!\log x|\Big(\sum_{\ell=1}^{\infty} 
      \mathcal{G}_0\big(\sqrt{p}[1+\log(x^{\ell})]\big) 
         -\frac{1}{3\sqrt{p}}\sum_{\ell=1}^{\infty} \mathcal{G}_1
	 \big(\sqrt{p}[1+\log(x^{\ell})]\big)\Big).
     \end{equation}
  Remembering that we are looking for an approximation of 
  $|\!\log x| \sum_{\ell=1}^{\infty}\big(\varphi(x^{\ell})\big)^p$ 
  and keeping in mind the relation \eqref{eqn_psi} between 
  $\varphi$, $p$ and $\psi_p$, 
  we state:
   \begin{prop}   \label{prop_gssn}
    There exists $C>0$ such that for all $p\geq1$
    and $x\in (0,1)$, 
     \begin{equation}   \label{eqn_gssn}
      \Big| |\!\log x| \sum_{\ell=1}^{\infty}
      \big(\varphi(x^{\ell})\big)^p - \mathbb{G}_p(x) \Big| 
      \leq \frac{C}{p +p|\!\log x|}. 
     \end{equation}
   \end{prop}
  
  \begin{crl}    \label{crl_fp}
   There exists $C>0$ such that for all $p\geq1$ and 
   $z\in \D^*$, 
    \begin{equation}\label{eq:3.19}
     \Big|\Big(\frac{2\pi}{p}\Big)^{3/2}\big(1+\nu(p)\big)^{-1}
    B_{p+1}^{\D^*}(z)  - \mathbb{G}_p(|z|^{2/p}) \Big| 
    \leq \frac{C}{p +2 |\!\log |z||}.
    \end{equation}
   In particular, 
   \begin{align}\label{eq:3.20}
\sup_{z\in \D^*} B_{p}^{\D^*}(z)
= \big(\frac{p}{2\pi}\big)^{3/2} 
   + \mathcal{O}(p).
\end{align}
  \end{crl}

\noindent
\textit{Proof of Proposition \ref{prop_gssn}. --- }
 We set:
  \begin{equation}\label{eq:3.21}
   \delta_p(\zeta) = \psi_p(\zeta) - 
   \mathcal{G}_0\big(\sqrt{p}(1-\zeta)\big) 
   +  \frac{1}{3\sqrt{p}}\mathcal{G}_1\big(\sqrt{p}(1-\zeta)\big).
  \end{equation}
 For all $p\geq 1$ and $x\in (0,1)$, 
 by (\ref{eqn_psi}), (\ref{eq:3.17}) and (\ref{eq:3.21}),
  \begin{equation}\label{eq:3.22}\begin{split}
 \Big| |\!\log x| \sum_{\ell=1}^{\infty}
      \big(\varphi(x^{\ell})\big)^p - \mathbb{G}_p(x) \Big| 
=  |\!\log x|\Big|\sum_{\ell=1}^{\infty} \delta_p\big(-\log(x^{\ell})
   \big) \Big|\\
   \leq 
   |\!\log x|\sum_{\ell=1,\, \ell\neq |\!\log x|^{-1}}^{\infty} 
   \big|\delta_p\big(-\log(x^{\ell})\big) \big|; 
 \end{split} \end{equation}    
 this takes into account the vanishing of $\delta_p(\zeta)$ 
 at $\zeta=1$. By Lemma \ref{lem_approx_psip}, 
 the latter is bounded above by 
  \begin{equation}\label{eq:3.23}
\frac{C}{p}|\!\log x|\sum_{\ell=1,  \ell\neq |\!\log x|^{-1}}^{\infty} 
\frac{1}{1+p(\ell\log x + 1)^2};   
  \end{equation}
 we can thus conclude if we bound this quantity above by an 
 expression of type $\frac{C}{p(1+|\!\log x|)}$.
 
 If $0<|\!\log x|\leq 2$, 
 by bounding the terms associated with 
 $\ell= \lfloor -(\log x)^{-1}\rfloor,  \lfloor -(\log x)^{-1}\rfloor 
 +1$ by $1$,  
 where we note $\lfloor u\rfloor$ the integer part of $u\in \R$,
 we get 
  \begin{align*}
   \sum_{\substack{\ell=1\\\ell\neq |\!\log x|^{-1}}}^{\infty} 
   \frac{1}{1+p(\ell\log x + 1)^2} 
        &\leq 2 + \frac{1}{|\log x|}\int_0^{\infty} 
 \frac{d\alpha}{1+p(\alpha-1)^2}     \\
       &\leq 2 + \frac{1}{|\log x|}\int_{-\infty}^{\infty} 
       \frac{d\alpha}{1+p(\alpha-1)^2}                           \\
                           &= 2 + \frac{\pi}{\sqrt{p}|\log x|}, 
  \end{align*}

 Thus by (\ref{eq:3.22}) and (\ref{eq:3.23}),
  for $0<|\!\log x|\leq 2$,
  \begin{equation}\label{eq:3.24}
       \Big| |\!\log x| \sum_{\ell=1}^{\infty}
      \big(\varphi(x^{\ell})\big)^p - \mathbb{G}_p(x) \Big| 
      \leq \frac{C|\!\log x|}{p} + \frac{C}{p^{3/2}} 
			\leq \frac{C}{p}, 
  \end{equation}
 and this yields the upper bound $\dfrac{C}{p(1+|\!\log x|)}$ as 
 $\dfrac{1}{1+|\!\log x|}\geq \dfrac{1}{3}$. 
 This way, the estimate \eqref{eqn_gssn} is proved on the region 
 $\{0<|\!\log x|\leq 2\}$. 
 
 Let us assume now that $|\!\log x|\geq 2$. 
 Then for all $\ell\geq 1$, 
 $\ell|\!\log x| - 1 \geq \dfrac{\ell|\!\log x|}{2}$, 
 thus \[(\ell\log x + 1)^2 = (\ell|\!\log x| - 1)^2
 \geq \dfrac{\ell^2|\!\log x|^2}{4},\] 
 and 
  \begin{equation}\label{eq:3.25}
   \sum_{\substack{\ell=1\\  \ell\neq |\!\log x|^{-1}}}^{\infty} 
   \frac{|\!\log x|^2}{1+p(\ell\log x + 1)^2} 
      \leq \sum_{\ell=1}^{\infty} \frac{|\!\log x|^2}{p\big(
      \frac{\ell^2|\!\log x|^2}{4}\big)} 
       = \frac{2\pi^2}{3p}. 
  \end{equation}
 In other words, by (\ref{eq:3.22}), (\ref{eq:3.23})
 and (\ref{eq:3.25}), 
 on the region $\{|\!\log x|\geq 2\}$,
  \begin{equation}\label{eq:3.26}
         \Big| |\!\log x| \sum_{\ell=1}^{\infty}
      \big(\varphi(x^{\ell})\big)^p - \mathbb{G}_p(x) \Big| 
       \leq \frac{C}{p^2 |\!\log x|},  
  \end{equation}
 and this upper bound yields here again an upper bound 
 $\dfrac{C}{p(1+|\!\log x|)}$, since $\dfrac{1}{|\!\log x|}
 \leq \dfrac{3}{2(1+|\!\log x|)}$ 
 when $|\!\log x|\geq 2$.
 By  (\ref{eq:3.24})
 and (\ref{eq:3.26}), we get \eqref{eqn_gssn}.
\cqfd
 
 ~
 
\noindent
\textit{Proof of Corollary \ref{crl_fp}. --- }
 The first part of the corollary follows at once from
 Proposition \ref{prop_gssn}. 
 The second part is an immediate consequence of the estimate 
  \begin{equation}\label{eq:3.27}
   \sup_{x\in (0,1)}     
   \mathbb{G}_p(x) 
   = 1+\mathcal{O}(p^{-1/2}). 
  \end{equation}
 To establish this estimate, 
 let us prove first that 
  \begin{equation}\label{eq:3.28}
 1\leq  \sup_{x\in (0,1)} |\!\log x|\sum_{\ell = 1}^{\infty}
   \mathcal{G}_0\big(\sqrt{p}[1+\log(x^{\ell})]\big) 
   = 1+\mathcal{O}(p^{-1/2}). 
  \end{equation}
 As putting $x=e^{-1}$ in $|\!\log x|\sum_{\ell = 1}^{\infty}
 \mathcal{G}_0\big(\sqrt{p}[1+\log(x^{\ell})]\big)$ 
 gives $1+\sum_{\ell = 2}^{\infty}\mathcal{G}_0
 \big(\sqrt{p}[1+\log(x^{\ell})]\big)\geq 1$, 
 we get already that the $\sup$ in \eqref{eq:3.28} 
 is bounded below by $1$. 

 Now we have
 \begin{align}\label{eq:3.29}
\sum_{\ell= \lfloor -(\log x)^{-1}\rfloor }^{ \lfloor -(\log 
x)^{-1}\rfloor +1} e^{-p(1+ \ell \log x)^2}
\leq 1+ e^{-p(\log x)^2/4} \quad \,
\text{if } |\log x|\leq 1,
\end{align}
and as a function of $s>0$,
$e^{-p(1+ s\log x)^2}$ increases when 
$s<\lfloor -(\log x)^{-1}\rfloor$ and decreases when
$s>\lfloor -(\log x)^{-1}\rfloor+1$, thus
\begin{align}\label{eq:3.30}\begin{split}
\Bigg(\sum_{\ell=1}^{\lfloor -(\log x)^{-1}\rfloor-1 }
&+ \sum_{\ell=\lfloor -(\log 
x)^{-1}\rfloor +2}^{\infty}\Bigg)
\mathcal{G}_0\big(\sqrt{p}  \big[1+\log(x^{\ell})\big]\big)\\
   & \leq  \int_{\R}\mathcal{G}_0\big(\sqrt{p}\big[1-s|\!\log(x)|\big]
	    \big)\,ds     \\
   &=  \int_{\R}\mathcal{G}_0\big(-\sqrt{p}s|\!\log(x)|\big)\,ds  
    =   \frac{C}{\sqrt{p}\,|\!\log x|},
\end{split}\end{align}
where we just omit the sum 
 $\sum_{\ell=1}^{\lfloor -(\log x)^{-1}\rfloor-1}$, 
 if $|\!\log x|\geq 1$;  the transition from the second
 to the third line simply comes from the translation 
 $s\leftarrow s+|\log x|^{-1}$. 

 Now $|\log x| e^{-p(1+\log x)^2} =
 (|\log x| -1) e^{-p(|\log x|-1)^2} 
 + e^{-p(|\log x|-1)^2} $. By using that the function
 $\eta\mapsto \eta e^{-\eta^2/2}$ is bounded on $\R$, we get
 from (\ref{eq:3.29}), (\ref{eq:3.30}) that for $x\in (0,1)$,
  \begin{equation}\label{eq:3.32}
   |\!\log x|\sum_{\ell = 1}^{\infty}
   \mathcal{G}_0\big(\sqrt{p}[1+\log(x^{\ell})]\big) 
   = \inf \{1, |\log x|\} + \mathcal{O}(p^{-1/2}),\:\: p\to +\infty. 
  \end{equation}
 With similar methods, one proves that 
  \begin{equation}\label{eq:3.33}
   \sup_{x\in (0,1)} |\!\log x|\sum_{\ell = 1}^{\infty}
   \mathcal{G}_1\big(\sqrt{p}[1+\log(x^{\ell})]\big) 
   = \mathcal{O}(1),\:\: p\to +\infty.  
  \end{equation}
  From  (\ref{eq:3.28}) and (\ref{eq:3.33}) we get 
  (\ref{eq:3.27}).
 \cqfd
 ~
 
\section{Elliptic Estimates for Kodaira Laplacians on 
\texorpdfstring{$\D^*$}{Delta*} 
and \texorpdfstring{$\Sigma$}{Sigma}}\label{S:EE}

In this section we establish a weighted elliptic estimate for
Kodaira Laplacians on $(\D^*, \omega_{\D^*})$ with weight 
$\big|\!\log(|z|^2)\big|^p$ such that the estimate is uniform on $p$,
and along $\D^*$.
This is a major 
analytic input in comparison with the 
compact situation.

Let $\dbar^{L^p*}$ be the adjoint of the Dolbeault operator
$\dbar^{L^p}$ on $(L^p, h^p)$ over $(\Sigma, \omega_{\Sigma})$.
Then the Kodaira Laplacian is defined as
\begin{align}\label{eq:4.0}
\square_{p}:= (\dbar^{L^p}+ \dbar^{L^p*})^2
= \dbar^{L^p}\dbar^{L^p*}+ \dbar^{L^p*}\dbar^{L^p}
: \Omega^{(0,\bullet)}(\Sigma, L^p)\to 
\Omega^{(0,\bullet)}(\Sigma, L^p).
\end{align}
We denote by $\square_{p}^{\D^*}$ the above operator when 
$(\Sigma,\omega_{\Sigma},L^p,h^p)$ is replaced by 
$(\D^*,\omega_{\D^*}, \C, \big|\!\log(|z|^2)\big|^p h_{0})$ 
(recall that $h_0$ denotes the standard flat Hermitian metric 
on the trivial line bundle $\C$).

\subsection{Estimate on the punctured disc 
\texorpdfstring{$\D^*$}{Delta*}: degree \texorpdfstring{$0$}{0}}
\label{s4.1}

Note that the Poincar\'e metric (\ref{eqn_omegaPcr}) 
on the punctured disc can be written as
 \begin{equation}  \label{eqn_PcrPttl}
  \omega_{\D^*} = - i\ddbar \log\big(\!-\log(|z|^2)\big). 
 \end{equation}

 Recall that the norm $\|\LargerCdot\|_{{\bL}^{2,2}_p(\D^*)}$
 was defined in Section \ref{ss:FS} {\bf(iii)}. 
 In what follows,  we adopt the notation $\mathsf{L}$ for 
the trivial line bundle $\C$ over the open unit disc $\D$,
thought of as endowed with the singular Hermitian metric 
$h_{\D^*} := \big|\!\log(|z|^2)\big| h_{0}$; 
similarly, for $p\geq 1$, $\mathsf{L}^p$ will implicitly refer to 
$\big(\C, \big|\!\log(|z|^2)\big|^p h_{0}\big)
= \big(\C, h_{\D^*}^p\big)$. 
Notice that with these conventions, \eqref{eqn_PcrPttl} 
can be interpreted as: 
\begin{align}\label{eq:4.1}
-i\omega_{\D^*} \text{ is the curvature of } (\mathsf{L},h_{\D^*})
  \text{ (and thus $-i\, p\omega_{\D^*}$ is that of 
$(\mathsf{L}^p,h_{\D^*}^p)$). } 
\end{align}
 
We prove in this section the following basic elliptic estimate on 
the Kodaira Laplacians $\square_p^{\D^*}$, 
associated to the data 
$\big(\mathbb{D}^*,\omega_{\mathbb{D}^*},
\mathsf{L}^p, h_{\D^*}^p \big)$.

\begin{prop}   \label{prop_TheEstimates}
    Let $s\geq 1$. Then there exists $C=C(s,h_{\D^*})$ such that
    for all $p\geq 1$, 
    and all $\sigma\in {\bL}^{2,2s}_p(\D^*)$, 
     \begin{equation}   \label{eqn_TheEstimates}
      \|\sigma\|_{{\bL}^{2,2s}_p(\D^*)}^2 
      \leq C\sum_{j=0}^{s} p^{4(s-j)}\big\|(\square_p^{\D^*})^j 
      \sigma\big\|_{{\bL}^{2}_p(\D^*)}^2.
     \end{equation}
   \end{prop}

Our strategy is as follows. 
We will write the detailed proof for $s=1$, 
and the proof follows then by induction on $s\in N^*$, 
we get it for $s\geq 2$.

For $s = 1$ we first establish an estimate analogous 
to \eqref{eqn_TheEstimates}
for the Laplace-Beltrami operator $\Delta_0$ of
$\omega_{\D^*}$, 
instead of the Kodaira Laplacian $\square_p^{\D^*}$
associated to $(\D^*, \omega_{\D^*}, \mathsf{L}^p, h_{\D^*}^p\big)$.
Then we deduce  \eqref{eqn_TheEstimates} by 
 K\"ahler identities. 

 To facilitate the computation, we introduce first
 new coordinates on $\D^*$ and explain some basic 
 geometric facts.
 
For $z\in \D^*$, we will use the coordinates 
$(t,\theta)\in \R\times (\R/2\pi \Z)$ with 
\begin{align}\label{eq:4.3}
t:=\log\big(-\log(|z|^2)\big), \quad z= |z| {e}^{i\theta}.
\end{align}
We denote also $\frac{\partial}{\partial t}$ by $\partial_{t}$,
and  $\frac{\partial}{\partial \theta}$ by $\partial_{\theta}$. 
Then we compute 
\begin{align}\label{eq:4.4}
\overline{z}\log(|z|^2)\frac{\partial}
    {\partial \overline{z}} 
  = \partial_{t} - \frac{i}{2}e^t\partial_{\theta},\quad 
  (\overline{z}\log(|z|^2))^{-1} d\overline{z}
  =\frac{1}{2}dt + i\, e^{-t}d\theta.
\end{align}
Thus we have
\begin{align}\label{eq:4.5}
\overline{\partial} = d\overline{z} \wedge 
\frac{\partial}{\partial \overline{z}} 
= (\frac{1}{2}dt + i\, e^{-t}d\theta)
( \partial_{t}- \frac{i}{2}e^t \partial_{\theta}),\quad
\partial t= \frac{1}{2}dt - i\, e^{-t}d\theta.
\end{align}
{}From (\ref{eq:4.5}) we obtain the following useful relation
\begin{align}\label{eq:4.6}
\omega_{\D^*}= - e^{-t}dt\wedge d\theta,\quad
 \big|\!\log(|z|^2)\big|^p \omega_{\D^*}  
 = -  e^{(p-1)t}dt\wedge d\theta,
\end{align}
and the metric associated with $\omega_{\D^*}$
in the coordinates $(t,\theta)$ is
\begin{align}\label{eq:4.7}
\frac{1}{2}(dt)^2 + 2 e^{-2t}( d\theta)^2,
\end{align}
thus $(\sqrt{2} \partial_{t}, \frac{1}{\sqrt{2}} e^t 
\partial_{\theta})$ is an orthonormal frame of $\omega_{\D^*}$.

Let $\nabla^{\D^*}$ be the Levi-Civita connection on
$(\D^*, \omega_{\D^*})$. Using \eqref{eq:4.7}
and the equality 
\[\nabla^{\D^*}_{\partial_{t}}\partial_{\theta}
- \nabla^{\D^*}_{\partial_{\theta}}\partial_{t}=
[\partial_{t}, \partial_{\theta}]=0,\] 
we compute that 
\begin{align}\label{eq:4.8}
\left\langle \nabla^{\D^*}_{\partial_{\theta}}\partial_{\theta} ,
\partial_{t}\right\rangle = -
\left\langle \nabla^{\D^*}_{\partial_{t}}\partial_{\theta} ,
\partial_{\theta}\right\rangle
= -\frac{1}{2}\partial_{t} \left\langle \partial_{\theta} ,
\partial_{\theta}\right\rangle
= 2 e^{-2t}.
\end{align}
From (\ref{eq:4.8}) we get
\begin{align}\label{eq:4.9}\begin{split}
   & \nabla^{\D^*}_{\partial_{t}}\partial_{t}=0,\quad
    \nabla^{\D^*}_{\partial_{\theta}}\partial_{\theta}
    = 4 e^{-2t} \partial_{t},\\
  & \nabla^{\D^*}_{\partial_{t}}\partial_{\theta}
= \nabla^{\D^*}_{\partial_{\theta}}\partial_{t}= -\partial_{\theta}.
\end{split}\end{align}
{}From (\ref{eq:4.9}) we get
\begin{align}\label{eq:4.10}\begin{split}
&\nabla^{\D^*} d\theta = d\theta\otimes dt 
+ dt\otimes d\theta,\quad
\nabla^{\D^*} dt = -4 e^{-2t} d\theta \otimes d\theta,\\
&\Delta_{0}= -2 (\partial_{t}\partial_{t}-\partial_{t})
- \frac{1}{2} e^{2t} \partial_{\theta}\partial_{\theta}.
\end{split}\end{align}

Let $\overline{\partial}^*$ (resp.\ $\overline{\partial}^{\mathsf{L}^p*}$)
be the adjoint of $\overline{\partial}$ on the trivial line bundle
$(\C,h_{0})$ (resp.\ on $(\C, \big|\!\log(|z|^2)\big|^p h_{0})$) 
over $(\D^*, \omega_{\D^*})$. 
By \eqref{eq:4.4} and  \eqref{eq:4.6}, we have 
the following expressions in the coordinates $(t,\theta)$,
\begin{align}\label{eq:4.12}
\overline{\partial}^{\mathsf{L}^p*}= \overline{\partial}^{*}
- p (\overline{\partial} t\wedge)^*
\quad \text{and }\, 
(\overline{\partial} t\wedge)^* d\overline{z}
= \left\langle  d\overline{z}, \overline{\partial} t\right\rangle
= \overline{z}\log(|z|^2).
\end{align}
By (\ref{eq:4.4}) and (\ref{eq:4.6}), we get for $f\in 
C^\infty(\D^*)$,
\begin{align}\label{eq:4.11}
\overline{\partial}^{\mathsf{L}^p*}( f \overline{\partial} t)
= - z\log(|z|^2)\frac{\partial}{\partial {z}} f 
+ (1-p) f, \quad
\overline{\partial}^{*} (f d\overline{z}) =
- |z|^{2}\log^{2}(|z|^2)\frac{\partial}{\partial {z}} f .
\end{align}
Thus the Kodaira Laplacian associated with 
$(\C, \big|\!\log(|z|^2)\big|^p h_{0})$ has the form 
\begin{equation}\label{eq:4.13}
\begin{split}
\square_p^{\D^*}= 
\overline{\partial}^{\mathsf{L}^p*} \overline{\partial}
+\overline{\partial}\, \overline{\partial}^{\mathsf{L}^p*}
&= \overline{\partial}^{*} \overline{\partial}
+\overline{\partial}\, \overline{\partial}^{*}
-p (  \overline{\partial} (\overline{\partial} t\wedge)^*
+ (\overline{\partial} t\wedge)^*\overline{\partial} )\\
&= \frac{1}{2}\Delta_{0}
-p  \left(\,  \overline{\partial} (\overline{\partial} t\wedge)^*
+ (\overline{\partial} t\wedge)^*\overline{\partial} \,\right),
\end{split}
\end{equation}
where we used the K\"ahler identity 
$\overline{\partial}^{*} \overline{\partial}
+\overline{\partial}\, \overline{\partial}^{*}=\frac12 \Delta_0$ 
for the last equality.

\smallskip
\noindent
\textit{Proof of Proposition \ref{prop_TheEstimates}. ---}
Notice that since the Hermitian line bundles $\mathsf{L}^p$ 
we consider here are powers of 
the line bundle $\big( \C, \big|\!\log(|z|^2)\big|h_0\big)$, 
the Chern connections $\nabla^p$ acting on the sections 
of these bundles, which are functions, 
are given by 
 \begin{equation}   \label{eqn_nablap}
  \nabla^p f = df + p f \,\partial t\,,\:\: f\in C^{\infty}(\D^*,L^p).
 \end{equation}
Therefore, for these $f$ and $p>1$,
 \begin{align}\label{eq:4.16}\begin{split}
  \|f\|_{\bL^{2,1}_p(\D^*)}^2 
   &=\int_{\D^*}  \big(|f|^2 + |\nabla^p_{\sqrt{2}\partial_{t}} f|^2 
   +|\nabla^p_{(1/\sqrt{2})e^t\partial_{\theta}} f|^2\big)
               \big|\!\log(|z|^2)\big|^p \omega_{\D^*}        \\
    &\leq 2 \int_{\D^*}  \big((p^2+1)|f|^2 + 2|\partial_{t} f|^2 
    +\tfrac{1}{2}|e^t\partial_{\theta} f|^2\big)
             \big|\!\log(|z|^2)\big|^p \omega_{\D^*},
 \end{split}\end{align}
and, similarly
 \begin{align}\label{eq:4.17}\begin{split}
  \|f\|_{{\bL}^{2,2}_p(\D^*)}^2  
  \lesssim \int_{\D^*}  \big(&p^4|f|^2 + p^2(|\partial_{t} f|^2 
  + |e^t\partial_{\theta} f|^2)                \\
& + |\partial_{t}^2 f|^2 + |e^t\partial_{t}\partial_{\theta} f|^2
  + |e^{2t}\partial_{\theta}^2 f|^2\big)
  \big|\!\log(|z|^2)\big|^p \omega_{\D^*}, 
\end{split} \end{align}
with the constants understood in $\lesssim$  
\textit{independent of} $p$.  

We will compute everything by using the coordinate $(t,\theta)$,
then $\int$ means $\int_{\R\times (\R/2\pi \Z)}$
and sometimes we identify $S^1$ to $ \R/2\pi \Z$.
We will work separately the real and imagine part of $f$, 
thus we assume now $f$ is real. 
Thus by (\ref{eq:4.6}) and simple integrations by parts, we get
 \begin{equation} \label{eq:4.18}
  \int_{\D^*} |e^t\partial_{\theta} f|^2 \big|\!\log(|z|^2)\big|^p 
  \omega_{\D^*} 
  = \int  |e^t\partial_{\theta} f|^2 e^{(p-1)t}dt d\theta
= \int (e^{2t}\partial_{\theta}^2 f)f 
e^{(p-1)t}dt d\theta,
 \end{equation}
and 
  \begin{equation}\label{eq:4.19}
  \int_{\D^*} |\partial_{t} f|^2 \big|\!\log(|z|^2)\big|^p 
  \omega_{\D^*} 
     = - \int (\partial_{t}^2 f) f     e^{(p-1)t}dt d\theta
       + \frac{(p-1)^2}{2} \int f^2  e^{(p-1)t}dt d\theta.
 \end{equation}
This way, by Young inequality, 
we obtain for every $\vareps>0$,
 \begin{equation}  \label{eqn_eps}
  \begin{aligned}
   \int_{\D^*} & (|\partial_{t} f|^2+|e^t\partial_{\theta} f|^2) 
   \big|\!\log(|z|^2)\big|^p \omega_{\D^*}          \\
              \leq &\Big(\vareps^{-1}+\frac{(p-1)^2}{2}\Big)
\int f^2  e^{(p-1)t}dt d\theta
      +\frac{\vareps}{2}\int \big(|\partial_{t}^2 f|^2
        +|e^{2t}\partial_{\theta}^2f|^2\big) e^{(p-1)t}dt d\theta.
  \end{aligned}
 \end{equation}
Taking $\vareps=p^{-2}$ we get
 \begin{multline}\label{eq:4.20}
  \int_{\D^*} p^2(|\partial_{t} f|^2+|e^t\partial_{\theta} f|^2)
  \big|\!\log(|z|^2)\big|^p \omega_{\D^*} \\
     \leq  2\int \big(p^4 |f|^2 + (|\partial_{t}^2 f|^2
     +|e^{2t}\partial_{\theta}^2 f|^2)\big) e^{(p-1)t}dt d\theta.
 \end{multline}
Thus, from (\ref{eq:4.17}) we have for $f\in {\bL}^{2,2}_p(\D^*) $,
 \begin{equation}\label{eq:4.21}
  \|f\|_{{\bL}^{2,2}_p(\D^*)}^2  \lesssim \int_{\D^*}  \big(p^4|f|^2 
     + |\partial_{t}^2 f|^2 + |e^t\partial_{t}\partial_{\theta} f|^2 
     + |e^{2t}\partial_{\theta}^2 f|^2\big)
     e^{(p-1)t}dt d\theta, 
 \end{equation}
with the implied constant independent of $p$. 

By (\ref{eq:4.10}),
 \begin{equation}    \label{eqn_L2Lapl}
  \begin{aligned}
 &  \int_{\D^*} (\Delta_0f)^2  \big|\!\log(|z|^2)\big|^p 
   \omega_{\D^*}                            \\        
                  &= 4 \int \big((\partial_t^2f)^2+(\partial_t f)^2
 +(\tfrac{e^{2t}}{4}\partial_{\theta}^2f)^2\big)e^{(p-1)t}\,
 dt d\theta  
 -8 \int(\partial_t^2f)(\partial_t f)e^{(p-1)t}\,dt d\theta     \\
& \quad+8\int(\partial_t^2f)(\tfrac{e^{2t}}{4}\partial_{\theta}^2 f)
   e^{(p-1)t}\,dt d\theta 
      -8\int(\partial_tf)(\tfrac{e^{2t}}{4}\partial_{\theta}^2 f)
      e^{(p-1)t}\,dt d\theta  .
  \end{aligned}
 \end{equation}
We deal with the mixed terms as follows: 
 \begin{itemize}
  \item $-8\int(\partial_t^2f)(\partial_t f)e^{(p-1)t}\,dt d\theta$: 
   an integration by parts yields: 
    \begin{equation}\label{eq:4.22}
     -8\int(\partial_t^2f)(\partial_t f)e^{(p-1)t}\,dt d\theta 
     =4 (p-1)\int (\partial_t f)^2 e^{(p-1)t}\,dt d\theta,  
    \end{equation}
   and we do not provide more efforts, as this quantity has 
   the favorable sign already 
   -- remember we want a bound below on 
   the $\bL^{2}_p(\D^*)$-norm of $\Delta_0f$; 
  \item $8\int(\partial_t^2f)(\tfrac{e^{2t}}{4}\partial_{\theta}^2 f)
  e^{(p-1)t}\,dt d\theta$: 
   exchanging $\partial_t$ and $\partial_{\theta}$ via integrations 
   by parts,  we get: 
    \begin{multline}\label{eq:4.23}
8\int (\partial_t^2f)(\tfrac{e^{2t}}{4}
\partial_{\theta}^2 f)e^{(p-1)t}\,dt d\theta    \\
   =  8\int(\tfrac{e^t}{2}\partial_t\partial_{\theta}f)^2 e^{(p-1)t}
   \,dt d\theta 
     - 8(p+1)\int(\partial_tf)(\tfrac{e^{2t}}{4}\partial_{\theta}^2f) 
     e^{(p-1)t}\,dt d\theta, 
    \end{multline}
   and collect the extra term 
   $-8(p+1)\int(\partial_tf)(\tfrac{e^{2t}}{4}\partial_{\theta}^2 f) 
   e^{(p-1)t}\,dt d\theta$ 
   together with the left over right-hand-side mixed term in 
   \eqref{eqn_L2Lapl}, 
   i.e. we deal with:
  \item $-8(p+2)\int(\partial_tf)(\tfrac{e^{2t}}{4}
  \partial_{\theta}^2f) e^{(p-1)t}\,dt d\theta$: 
    \begin{multline}\label{eq:4.24}
-8(p+2)\int (\partial_tf)(\tfrac{e^{2t}}{4}\partial_{\theta}^2 f) 
      e^{(p-1)t}\,dt d\theta                              \\
\geq  - 2\int (\tfrac{e^{2t}}{4}\partial_{\theta}^2 f)^2 
e^{(p-1)t}\,dt d\theta
       -8(p+2)^2 \int ( \partial_{t} f)^2 e^{(p-1)t}\,dt d\theta, 
    \end{multline}
By Cauchy-Schwarz inequalities and  (\ref{eq:4.19}),
    we get 
    \begin{align*}
    - 8(&p+2)^2\int ( \partial_{t} f)^2 e^{(p-1)t}\,dt d\theta \\
&\geq - 2\int ( \partial_{t}^2 f)^2 e^{(p-1)t}\,dt d\theta 
  -\big(8(p+2)^4+4 (p-1)^2(p+2)^2 \big) \int f^2 e^{(p-1)t}\,
  dt d\theta. 
    \end{align*}
  \end{itemize}
  We sum up what precedes as: 
   \begin{equation}\label{eq:4.25}
    \begin{aligned}
\int_{\D^*} (\Delta_0f)^2  \big|\!\log(|z|^2)\big|^p \omega_{\D^*} 
   &\geq  2 \int \big( (\partial_t^2f)^2 
   + (\tfrac{e^{t}}{2}\partial_t\partial_{\theta}^2f)^2
 +(\tfrac{e^{2t}}{4}\partial_{\theta}^2f)^2\big)e^{(p-1)t}\,
 dt d\theta                            \\
  &\quad- \big(8(p+2)^4+ 4(p-1)^2(p+2)^2 \big) 
  \int f^2 e^{(p-1)t}\, dt d\theta, 
    \end{aligned}
   \end{equation}
   By (\ref{eq:4.21}) and (\ref{eq:4.25}),  we get
   \begin{equation}   \label{eqn_TheEstimateLapl}
    \|f\|_{{\bL}^{2,2}_p(\D^*)}^2 \leq C\left(\|\Delta_0 f\|_{\bL^{2}_p
    (\D^*)}^2 + p^4\|f\|_{\bL^{2}_p(\D^*)}^2 \right),
   \end{equation}
  for some $C>0$ independent of both $p\geq 1$ and
  $f\in C^{\infty}_{0}(\D^*)$ with real values;  
  by density, this readily generalizes to $f \in {\bL}^{2,2}_p(\D^*)$ 
  with complex values, as $\Delta_0$ is a real operator. 
  
  We now carry out the replacement of $\Delta_0$ by 
  $\square_p^{\D^*}$ in \eqref{eqn_TheEstimateLapl}, 
  to get the desired estimate \eqref{eqn_TheEstimates}. 
  By (\ref{eq:4.12}) and (\ref{eq:4.13})
  we have the following identities for operators action 
  on functions defined on $\D^*$, 
   \begin{equation}\label{eq:4.27}
\square_p^{\D^*}=\frac{1}{2}\Delta_0 - p\widetilde{\partial}
\quad \text{with }
\widetilde{\partial} = \overline{z}\log(|z|^2)\frac{\partial}
    {\partial \overline{z}} 
  = \frac{\partial}{\partial t} - \frac{i}{2}e^t\frac{\partial}
  {\partial \theta}\,\cdot
   \end{equation}
  Let $f\in {\bL}^{2,2}_p(\D^*)$. 
Using inequality \eqref{eqn_eps} with $\vareps>0$ to be adjusted, 
(\ref{eq:4.21}) and (\ref{eqn_TheEstimateLapl}), we have:
   \begin{align}\label{eq:4.28}\begin{split}
    \int \big|\widetilde{\partial}f\big|^2 e^{(p-1)t}dt d\theta 
 &\leq 2\int (|\partial_tf|^2+|\tfrac{e^t}{2}\partial_{\theta}f|^2) 
	e^{(p-1)t}dt d\theta       \\
  &\leq \big(2\vareps^{-1}+p^2+\vareps C p^4 \big)
\|f\|_{\bL^{2}_p(\D^*)}^2 
+ \vareps C \|\Delta_0 f\|_{\bL^{2}_p(\D^*)}^2. 
\end{split} \end{align}
 From (\ref{eq:4.27}) and  (\ref{eq:4.28}) we are led to:
   \begin{align}\label{eq:4.29}\begin{split}
    \|\Delta_0 f &\|_{\bL^{2}_p(\D^*)}^2 
    =   \big\|2\big(\square_p^{\D^*} 
    + p \widetilde{\partial}\big) f\big\|_{\bL^{2}_p(\D^*)}^2 
      \leq 8 \|\square_p^{\D^*} f\|_{\bL^{2}_p(\D^*)}^2 
 + 8p^2\big\|\widetilde{\partial} f\big\|_{\bL^{2}_p(\D^*)}^2    \\
     &\leq 8 \|\square_p^{\D^*} f\|_{\bL^{2}_p(\D^*)}^2 
     + 8p^2\big(2\vareps^{-1}+p^2
     +\vareps C p^4 \big)\|f\|_{\bL^{2}_p(\D^*)}^2 
   + 8p^2\vareps C \|\Delta_0 f\|_{\bL^{2}_p(\D^*)}^2. 
  \end{split}  \end{align}
  Take $\vareps= \dfrac{1}{16Cp^2}$ to conclude that: 
   \begin{equation}   \label{eqn_estDelta0}
    \|\Delta_0 f\|_{\bL^{2}_p(\D^*)}^2 
    \leq 16\|\square_p f\|_{\bL^{2}_p(\D^*)}^2 
    + 2  Ap^4\|f\|_{\bL^{2}_p(\D^*)}^2 
   \end{equation}
  with  $A=256C+9$, say. 
  Plugged back into \eqref{eqn_TheEstimateLapl}, 
  this estimate gives exactly \eqref{eqn_TheEstimates},
  with a (new) constant $C>0$, uniform for $p\geq1$ and 
  $f\in {\bL}^{2,2}_p(\D^*)$. 

 The proof of Proposition \ref{prop_TheEstimates}
for $s=1$ is completed. Continuing by induction on $s\in \N^{*}$
we obtain the assertion for all $s\in \N^*$. 
\cqfd

\subsection{Estimate on the punctured Riemann surface 
\texorpdfstring{$\Sigma$}{Sigma} :  degree 0}\label{s4.2}

 We now consider the geometric situation of a punctured polarized 
 Riemann surface 
$(\Sigma,\omega_\Sigma,$ $L,h)$, satisfying moreover conditions 
$(\alpha)$ and $(\beta)$ of the Introduction.
Let $a\in D$. By assumption the following holds: 
 \textit{there exists a trivialization of $L$ around $a$ such that in 
 the associated local complex coordinate $z\in \D$, we have
  $h=\big|\!\log(|z|^2)\big|h_0$ on some coordinate disc 
  $\D_{r}$ centered at a and of radius $r\in(0,e^{-1})$.} 
 This way, the curvature $R^L$  
 of $h$ coincides with 
 $-i \omega_{\D^*}$ on $\D^*_{r}:=\D_{r}\smallsetminus\{0\}$. 

 \begin{prop}  \label{prop_EstimateSigma}
    For every $s\in\N^*$ there exists $C=C(s,h)$ such that for 
    all $p\gg 1$, and all $\sigma\in {\bL}^{2,2s}_p(h)
    ={\bL}^{2,2s}\big(\Sigma,\omega_\Sigma, L^p, h^p\big)$, 
     \begin{equation}   \label{eqn_EstimateSigma}
 \|\sigma\|_{{\bL}^{2,2s}_p(h)}^2 \leq C\sum_{j=0}^{s} p^{4(s-j)}\|
 (\square_p)^j \sigma\|_{{\bL}^2_p(h)}^2, 
     \end{equation}
    where $\square_p$ is the Kodaira Laplacian on $\Sigma$ 
    associated to $\omega_{\Sigma}$ and $h^p$. 
   \end{prop}
 \prf.\textit{ --- }
  Again, we do it for $s=1$. 
    
In the situation of the Proposition, we denote by $\overline{h}$ 
a smooth Hermitian metric on $L$ on the whole 
$\overline\Sigma$ such that  it coincides with $h$ on 
$\Sigma\smallsetminus \D_{r/2}$. 
It is an easy exercise to construct $\overline{h}$ so that 
$i \times$ its curvature, $\overline{\omega}$ say, is K\"ahler over 
the whole (compact) $\overline{\Sigma}$, 
which we take for granted until the end of this proof. 
Notice that $\omega_\Sigma$ and $\overline{\omega}$ coincide on 
$\Sigma\smallsetminus\D_{r/2}$. 
  
  Now the principle of the proof is to glue estimate 
  \eqref{eqn_TheEstimates} to the analogous 
  estimate for $\big(\overline\Sigma, L, \overline{h}\big)$, 
  that states the existence of $C>0$ such that for all 
  $p\gg 1$,  and all 
  $\sigma\in {\bL}^{2,2}_p\big(\overline{h}\big):=
 {\bL}^{2,2}\big(\overline\Sigma,\overline{\omega}, L^p, 
 \overline{h}^p\big)$, 
  \begin{align}\label{eq:4.30}
\|\sigma\|_{{\bL}^{2,2}_p(\overline{h})}^2 \leq 
C \big(\big\|\square_p^{\overline{\Sigma}} 
\sigma\big\|_{{\bL}^2_p(\overline{h})}^2
+ p^4\| \sigma\|_{{\bL}^2_p(\overline{h})}^2\big).
\end{align} 
  This estimate, as well as its generalization for 
  $\sigma\in \bL^{2, 2s}_p\big(\overline{h}\big)$, $s\geq1$, 
  can be found for instance in 
  \cite[(4.14)]{DLM06} or \cite[\S1.6.2]{mm}. 
  
 We denote $\nabla^{p,\Sigma*}$ the formal adjoint 
 of $\nabla^{p,\Sigma}$ acting on 
 $\Lambda (T^{*(0,1)}\Sigma) \otimes L^p$.
 By Lichnerowicz formula \cite[Remark 1.4.8]{mm}, 
 \begin{align}\label{eq:4.31}
2\square_{p}= \nabla^{p,\Sigma*}\nabla^{p,\Sigma}
-p R^L(w,\overline{w})
+(2p R^L + R^{T^{(1,0)}\Sigma}) (w,\overline{w}) \overline{w}^*
\wedge i_{\overline{w}},
\end{align}
and $w$ is an orthonormal frame of $T^{(1,0)}\Sigma$.
By \eqref{eq:4.1}, 
\begin{align}\label{eq:4.32}
R^L(w,\overline{w})=1, \quad
R^{T^{(1,0)}\Sigma}(w,\overline{w})= -2 
\quad \text{on } V= V_{1}\cup \ldots \cup V_{N}.
\end{align}
  From \eqref{eq:4.31} and  \eqref{eq:4.32}, 
 we have for any $\sigma\in {\bL}^{2,2}_p(h)$,
  \begin{align}\label{eq:4.33}
\Big| \|\nabla^{p,\Sigma} \sigma\|_{\bL^{2}_p(h)}^2
- 2 \left\langle  \square_{p}\sigma,\sigma\right\rangle_{p}
\Big| \leq C p \|\sigma\|_{\bL^{2}_p(h)}^2.
\end{align}

  Let $\chi$ be a cut-off function supported near $a$; 
  assume, more precisely, that  
  \begin{equation}\label{e:chi}
  \chi\in C^\infty\big(\overline\Sigma\big),\:\: 0\leq \chi\leq1, \:\:
  \chi\mid_{\,\overline{\D}_{r/2}}\equiv 1,\:\:
  \chi\mid_{\,\Sigma\smallsetminus \D_{2r/3}}\equiv 0.
  \end{equation}
  Let $p\geq 1$, and $\sigma\in {\bL}^{2,2}_p\big(h\big)
       ={\bL}^{2,2}\big(\Sigma,\omega_\Sigma, L^p, h^p\big)$. 
  Then $(1-\chi)\sigma \in {\bL}^{2,2}_p(\overline{h})$ and on 
  its support, $h$ coincides with $\overline{h}$; 
  likewise, $\chi\sigma$ can be interpreted as an element 
  of ${\bL}^{2,2}_p(\D^*)$ and on its support, 
  $h$ can be regarded as $h_{\D^*}$. 
  Therefore, 
   \begin{equation}   \label{eqn_split}
    \begin{aligned}
     \|\sigma\|_{{\bL}^{2,2}_p(h)}^2 &=\|\chi \sigma 
     + (1-\chi)\sigma\|_{{\bL}^{2,2}_p(h)}^2 
    \leq  2\big(\|\chi\sigma\|_{{\bL}^{2,2}_p(\D^*)}^2
     +\|(1-\chi)\sigma \|_{{\bL}^{2,2}_p(\overline{h})}^2\big)   \\
&\leq 2C\big(\big\|\square_p^{\D^*}(\chi\sigma) 
\big\|_{{\bL}^2_p(\D^*)}^2
+\big\|\square_p^{\overline{\Sigma}}[(1-\chi)\sigma]
\big\|_{{\bL}^2_p(\overline{h})}^2 \big)                \\
  &         \qquad \qquad +2Cp^4\big(\|\chi\sigma
  \|_{{\bL}^2_p(\D^*)}^2
  +\|(1-\chi)\sigma\|_{{\bL}^2_p(\overline{h})}^2\big)
    \end{aligned}
   \end{equation}
  where $C=\sup\big(C(h_{\D^*}),C(\overline{h})\big)$, with 
  $C(h_{\D^*})$,  resp. $C(\overline{h})$, 
  the constant from \eqref{eqn_TheEstimates}, resp. 
  from its analogue for $\big(\overline\Sigma, L, \overline{h}\big)$.  
  Thus defined, $C$ is independent of $\sigma$ and $p$. 
  
  Now $\|\chi\sigma\|_{{\bL}^2_p(h_{\D^*})}^2 = 
  \|\chi\sigma\|_{{\bL}^2_p(h)}^2\leq \|\sigma\|_{{\bL}^2_p(h)}^2$, 
  and $\|(1-\chi)\sigma\|_{{\bL}^2_p(\overline{h})}^2 
  \leq \|\sigma\|_{{\bL}^2_p(h)}^2$ as well. 
  The treatments of $\big\|\square_p^{\D^*}(\chi\sigma)
  \big\|_{{\bL}^2_p(h_{\D^*})}^2$ and 
  $\big\|\square_p^{\overline{\Sigma}}[(1-\chi)\sigma]
  \big\|_{{\bL}^2_p(\overline{h})}^2$ are done in the same spirit, 
  but require a little extra work.  For instance, we have on $\D_{2r/3}$, 
  by (\ref{eq:4.13}) and (\ref{eqn_nablap}), 
   \begin{equation} \label{eq:4.34}
    \begin{aligned}
    \square_p^{\D^*}(\chi\sigma) = \chi \square_p^{\D^*}\sigma
     -\big(\dbar\chi,\overline{(\nabla^{p})^{1,0}\sigma}
     \big)_{T^*\D^*}
      -\big(\dbar\sigma,\overline{\partial\chi}\big)_{T^*\D^*} 
  + \Big(\frac{1}{2}\Delta_0\chi\Big)\sigma  , 
  \end{aligned} \end{equation}
  hence 
   \begin{align} \label{eq:4.35}
    \begin{aligned}
    \big\|\square_p^{\D^*}(\chi\sigma)
    \big\|_{{\bL}^2_p(h_{\D^*})}^2                                                      
     &\leq 4\Big(\big\|\chi\square_p^{\D^*}\sigma
     \big\|_{{\bL}^2_p(h_{\D^*})}^2
     +\big\||\partial\chi|\dbar\sigma\big\|_{{\bL}^2_p(h_{\D^*})}^2    \\
     & \quad\qquad+\big\||\partial\chi|(\nabla^{p})^{1,0}
     \sigma\big\|_{{\bL}^2_p(h_{\D^*})}^2
   +\big\|(\tfrac{1}{2}\Delta_0\chi)\sigma
   \big\|_{{\bL}^2_p(h_{\D^*})}^2\Big)                           \\
     &=     4\Big(\big\|\chi\square_p\sigma
     \big\|_{{\bL}^2_p(h)}^2
  +\big\||\partial\chi|\dbar\sigma\big\|_{{\bL}^2_p(h)}^2    \\
& \quad\qquad+\big\||\partial\chi|(\nabla^{p,\Sigma})^{1,0}\sigma
     \big\|_{{\bL}^2_p(h)}^2
 +\big\|(\tfrac{1}{2}\Delta_0\chi)
 \sigma\big\|_{{\bL}^2_p(h)}^2\Big),                                     
   \end{aligned}\end{align}
  with the ${\bL}^2_p(h_{\D^*})$-norms, resp. 
  ${\bL}^2_p(h)$-norms,
  for 1-forms, resp. 1-forms with value in $L^p$, 
  computed with $\omega_{\D^*}\otimes h_{\D^*}^p$ on $\D^*$, 
  resp. with $\omega_\Sigma\otimes h^p$.  
  Consequently,  
   \begin{align}\label{eq:4.36}
    \begin{aligned}
\big\|\square_p^{\D^*}&(\chi\sigma)
\big\|_{{\bL}^2_p(h_{\D^*})}^2 
  \leq 4C\Big(\big\|\square_p\sigma\big\|_{{\bL}^2_p(h)}^2
  +\big\|\nabla^{p,\Sigma}\sigma\big\|_{{\bL}^2_p(h)}^2
    +\big\|\sigma\big\|_{{\bL}^2_p(h)}^2\Big)
    \end{aligned}\end{align}
  with $C=1+\max\big\{\big\|\tfrac{1}{2}
  \Delta_0\chi\big\|^2_{C^{0}(\D^*, \omega_{\D^*})}, 
\big\|\dbar\chi\big\|^2_{C^{0}(\D^*, \omega_{\D^*})}\big\}$, 
  that does not depend on $p$.  

  From (\ref{eq:4.33}) and (\ref{eq:4.35}) we get
  \begin{equation}\label{eq:4.40}
   \big\|\square_p^{\D^*}(\chi\sigma)
   \big\|_{{\bL}^2_p(h_{\D^*})}^2 
 \leq C\big(\big\|\square_p\sigma\big\|_{{\bL}^2_p(h)}^2
    +p \|\sigma\big\|_{{\bL}^2_p(h)}^2\big)
  \end{equation}
 for some $C$ independent of $p$ and $\sigma$.  
 Similarly, 
  \begin{equation}\label{eq:4.41}
   \big\|\square_p^{\overline{\Sigma}}\big((1-\chi)\sigma\big)
   \big\|_{{\bL}^2_p(\overline{h})}^2 
          \leq C\big(\big\|\square_p\sigma\big\|_{{\bL}^2_p(h)}^2
             +p \|\sigma\big\|_{{\bL}^2_p(h)}^2\big)
  \end{equation}
 with $C$ again independent of $p$ and $\sigma$. 
 
 In conclusion, it follows from (\ref{eqn_split}),
 (\ref{eq:4.40}) and (\ref{eq:4.41}), that there exists
 $C>0$ such that for any $p\gg1$ and 
 $\sigma\in {\bL}^{2,2}_p(h)$, we have
  \begin{equation}\label{eq:4.42}
     \|\sigma\|_{\bL^{2, 2}_p(h)}^2 
     \leq  C\big(\big\|\square_p\sigma\big\|_{{\bL}^2_p(h)}^2
  +p^4\|\sigma \|_{{\bL}^2_p(h)}^2\big).      \\
  \end{equation}
 
The proof of Proposition \ref{prop_EstimateSigma} for $s=1$ 
is completed.
The proof for general $s\in\N^*$ follows by induction
with the help of Proposition \ref{prop_TheEstimates}.
 \cqfd
    
 \begin{rmk}   \label{rmk:4.3}
  We made the assumption that $h$ polarizes 
  $\omega_{\Sigma}$ \emph{on the whole} $\Sigma$ 
  to fix ideas.  It is actually superfluous here, provided that 
  $h$ polarizes $\omega_{\Sigma}$ \emph{near the puncture(s)}, 
  and \emph{has curvature bounded below by} 
  $\vareps\omega_{\Sigma}$ \emph{on the whole} $\Sigma$;  
  both of these are automatically implied by conditions 
  $(\alpha)$ and $(\beta)$.  
 \end{rmk}

  \subsection{Bidegree \texorpdfstring{$(0,1)$}{(0,1)}}
  
This subsection will not be used in the rest of this paper, we 
include it here  only for completeness and its independent interest.
  
To prove that Propositions \ref{prop_TheEstimates} 
and \ref{prop_EstimateSigma} still hold in bidegree $(0,1)$, 
or, namely, for $\sigma$ a section of 
$T^{*(0,1)}\D^*\otimes \mathsf{L}^p$ or $T^{*(0,1)}
\Sigma\otimes L^p$, 
 an easy procedure is to observe that the following diagram: 
    \begin{equation}\label{eq:4.44}
    \xymatrix@C=3pc  {
    C_0^{\infty}(\D^*) 
    \ar[rr]^{\LargerCdot\otimes\frac{d\overline{z}}
	{\overline{z}\log(|z|^2)}\qquad\,\,} 
    \ar[d]_{\square_{p}- ie^t\partial_{\theta}} 
    && C_0^{\infty}(\D^*,T^{*(0,1)}\D^*) 
    \ar[d]^{\square_{p}} \\
    C_0^{\infty}(\D^*) \ar[rr]_{\LargerCdot\otimes\frac{d\overline{z}}
	{\overline{z}\log(|z|^2)}\qquad\,\,}
    && C_0^{\infty}(\D^*,T^{*(0,1)}\D^*)
  }
    \end{equation}
   commutes, where the horizontal arrows are isometries under
   $h_{\D^*}^p$ and $(h_{\D^*})^p\otimes\omega_{\D^*}$. 
  Indeed,  by (\ref{eq:4.4}) and (\ref{eq:4.5}),
   $\frac{d\overline{z}}{\overline{z}\log(|z|^2)}=\dbar t$,
   and by  (\ref{eq:4.5}) and (\ref{eq:4.11}), 
   for $g\in C^\infty(\D^*)$, we have
   \begin{equation}\label{eq:4.45}
   \begin{split}
\dbar\, \dbar^{L^p*} (g\dbar t)
&= \Big[\dbar^{L^p*}\dbar  g + \overline{z}\log(|z|^2)\frac{\partial}
    {\partial \overline{z}} g
    - {z}\log(|z|^2)\frac{\partial}{\partial {z}}g \Big] \dbar t\\
   & = (\square_{p}g - ie^t\partial_{\theta}g)\dbar t.
    \end{split}
\end{equation}
   \begin{prop}   \label{prop_TheEstimates2}
    Let $s\in \N^{*}$. 
    Then there exists $C=C(s,h_{\D^*})$ such that for all $p\geq 1$, 
    and all 
    $\sigma\in {\bL}^{2,2s}_p(\D^*) 
    = \bL^{2, 2s}\big(\D^*,\omega_{\D^*}, T^{*(0,1) }
    \D^*\otimes \mathsf{L}^p,\omega_{\D^*}\otimes h_{\D^*}^p\big)$, 
     \begin{equation}   \label{eqn_TheEstimates2}
      \|\sigma\|_{{\bL}^{2,2s}_p(\D^*)}^2 
      \leq C\sum_{j=0}^{s} p^{4(s-j)}\big\|(\square_p^{\D^*})^j 
      \sigma\big\|_{{\bL}^2_p(\D^*)}^2.
     \end{equation}
   \end{prop}
  \textit{Proof. --- }
   Indeed, take 
   $\sigma=f\frac{d\overline{z}}{\overline{z}\log(|z|^2)}
   = f \dbar t
   \in C^{\infty}(\D^*,T^{*(0,1)\D^*})$. 
   Then for instance, 
    \begin{equation*}
     (\nabla^p)^2\sigma = (\nabla^p)^2f\otimes\dbar t 
                           + 2\nabla^pf\otimes\nabla^{\D^*}\dbar t 
                           + f \otimes(\nabla^{\D^*})^2\dbar t
    \end{equation*}
   where $\nabla^{\D^*}$ is the Levi-Civita connection of 
   $\omega_{\D^*}$.  
  By (\ref{eq:4.5}),  (\ref{eq:4.7}) and (\ref{eq:4.9}),  
  $\dbar t$ is uniformly bounded at any order with respect to
   $\omega_{\D^*}$ on $\D^*$,  we get that 
    \begin{equation*}
     \int_{\D^*}\big|(\nabla^p)^2\sigma\big|^2_p\,\omega_{\D^*} 
\lesssim \int_{\D^*}\big(\big|(\nabla^p)^2f\big|^2_p 
	+ |\nabla^pf|^2_p + |f|^2\big)\,\omega_{\D^*} 
	= \| f\|_{{\bL}^{2,2}_p(\D^*)}^2,
    \end{equation*}
   independently of $p$. 
   By Proposition \ref{prop_TheEstimates}, 
   we thus have 
    \begin{equation}   \label{eqn_estmtesigma}
     \int_{\D^*}\big|(\nabla^p)^2\sigma\big|^2_p\,\omega_{\D^*} 
        \leq C \int_{\D^*}\big(\big|\square_p f\big|^2_p  
	+ p^4 |f|^2\big)\,\omega_{\D^*}
    \end{equation}
   for $p\gg1$, with $C$ independent of $p$. 
   By \eqref{eq:4.44}, we have
   \begin{equation}\label{eq:4.47}
   \begin{split}
     \int_{\D^*} \big|\square_p \sigma \big|^2_p\,\omega_{\D^*} 
         &= \int_{\D^*} \big|\square_p f
- ie^t\partial_{\theta}f\big|^2_p\,\omega_{\D^*}    \\   
 &\geq \frac{1}{2}\int_{\D^*}  |\square_p f|^2_p\,\omega_{\D^*} 
   - \int |e^t\partial_{\theta}f\big|^2 e^{(p-1)t}dt d\theta.
   \end{split}
   \end{equation}
   By \eqref{eqn_eps} with $\varepsilon=\dfrac{1}{2^9 p^2}$,
   \eqref{eq:4.25} and \eqref{eqn_estDelta0},
we obtain as in \eqref{eq:4.28},
   \begin{align}\label{eq:4.48}
\int \big|e^t\partial_{\theta}f\big|^2 e^{(p-1)t}dt d\theta
\leq C p^2 \|f\|_{{\bL}^{2}_p(\D^*)}^2
+ \frac{1}{4 p^2}\|\square_p f\|_{{\bL}^2_p(\D^*)}^2.
\end{align}
   As $\int_{\D^*} |\sigma|^2_p\,\omega_{\D^*}
   = \int_{\D^*} |f|^2_p\,\omega_{\D^*}$, we get from
   (\ref{eq:4.47}) and  (\ref{eq:4.48}), 
    \begin{equation*}
     \int_{\D^*} \big|\square_p f \big|^2_p\,\omega_{\D^*} 
         \leq 4\int_{\D^*}  |\square_p \sigma|^2_p\,\omega_{\D^*} 
     +C p^2  \int_{\D^*} |\sigma|^2_p\,\omega_{\D^*} .
    \end{equation*}
   This yields, coming back to \eqref{eqn_estmtesigma}, 
    \begin{equation}\label{eq:4.49}
     \int_{\D^*}\big|(\nabla^p)^2\sigma\big|^2_p\,\omega_{\D^*} 
        \leq C \int_{\D^*}\big(\big|\square_p \sigma\big|^2_p  
	+ p^4 |\sigma|^2\big)\,\omega_{\D^*}.  
    \end{equation}
  Consequently, we get for $p \gg 1$, 
   \begin{equation}\label{eq:4.50}
  \|\sigma\|_{{\bL}^{2,2}_p(\D^*)}^2
  \leq C\big(\|\square_p\sigma\|_{{\bL}^2_p(\D^*)}^2
  + p^4\|\sigma\|^2_{{\bL}^2_p(\D^*)}).
 \end{equation} 
 Now by induction on $s$, the assertion of 
 Proposition \ref{prop_TheEstimates2} follows for all $s\in \N^{*}$.
  \cqfd
  \smallskip

  Using moreover the same gluing procedure as in proving
  Proposition \ref{prop_EstimateSigma}, 
  we obtain the analogue of \eqref{eqn_TheEstimates2} on $\Sigma$, 
  when $\omega_{\Sigma}$ and $h$ 
  satisfy conditions $(\alpha)$ and $(\beta)$: 
   \begin{prop}   \label{prop_TheEstimateSigma2}
    Let $s\in\N^*$. 
    Then there exists $C=C(s,h)$ such that for all $p\geq 1$, 
    and all $\sigma\in {\bL}^{2,2s}_p(\Sigma)
    = \bL^{2,2s}\big(\Sigma,\omega_\Sigma, T^{*(0,1) }\Sigma
    \otimes L^p,\omega_\Sigma\otimes h^p\big)$, 
     \begin{equation}   \label{eqn_TheEstimateSigma2}
      \|\sigma\|_{{\bL}^{2,2s}_p(\Sigma)}^2 
      \leq C\sum_{j=0}^{s} p^{4(s-j)}\big\|(\square_p)^j 
      \sigma\big\|_{\bL^{2}_p(\Sigma)}^2.
     \end{equation}
   \end{prop}
 
 \section{Spectral Gap and Localization}\label{S:sgl}
We follow in this Section the localization scheme based on
the spectral gap and finite propagation speed \cite{mm} and 
show that the Bergman kernel localizes
near the singularities. As a consequence we obtain a first rough 
estimate, which will
be improved in the next Section.
   
Let $(M,\omega_M)$ be a complete K\"ahler manifold.
We will denote by $R^{\det}$ the curvature of the anticanonical 
line bundle
$(K^*_M, h^{K^*_M})$, where $h^{K^*_M}$ is induced by 
$\omega_M$. 

Let $(E,h^E)$ be a Hermitian holomorphic line bundle on $M$.
Let $\overline\partial^*$ be the formal adjoint of 
$\overline\partial$ with respect to 
$\langle\LargerCdot,\LargerCdot\rangle$ (cf. \eqref{eq:2.1}).
Let $\Box^E=\overline\partial^*\overline\partial$ be 
the Kodaira Laplace operator.
By \cite[Corollary\,3.3.4]{mm} the operator 
$\Box^E: C^\infty_{0}(M,E)\to C^\infty_{0}(M,E)$ 
is essentially
self-adjoint and we will denote its unique self-adjoint extension with
the same symbol
$\Box^E$. Note that the domain of this extension is 
$\operatorname{Dom}(\Box^E)=
\{\sigma\in {\bL}^2(M,E):\Box^E \sigma\in {\bL}^2(M,E)\}$.

Consider now a Hermitian holomorphic line bundle $(L,h)$ and
denote by $\Box_p:=\Box^{L^p}$ the Kodaira Laplace operator
corresponding to $(L^p,h^p)$.  
By \cite[Theorem\,6.1.1]{mm} and its proof we have the following.

   \begin{prop}[Spectral gap]  \label{prop_SpectralGap}
    Let $(M,\omega_M)$ be a complete K\"ahler manifold
and $(E,h^E)$ be a Hermitian holomorphic line bundle on $M$.
We assume there exist $\varepsilon>0$, $C>0$ such that 
$iR^E\geq\varepsilon\omega_M$ and $iR^{\det}\geq-C\omega_M$.
Then there exists $c=c(C,\varepsilon)>0$ such that for 
all $p\gg1$ we have
     \begin{equation}   \label{eqn_SpectralGap}
      {\rm Spec}(\square_p) \subset \{0\}\cup [cp, +\infty ) \,. 
     \end{equation}
   \end{prop}
 
 \begin{crl}\label{t5.2}
 The spectral gap \eqref{eqn_SpectralGap} holds for the Laplacian
 $\Box_p$ in the following situations:

\smallskip 
 (1) $(M,\omega_M)=(\D^*, \omega_{\D^*})$, 
 $(E,h^E)=(\C,\big|\!\log(|z|^2)\big| h_{0})$;
\smallskip

 (2)  $(M,\omega_M)=(\Sigma, \omega_\Sigma)$, $(E,h^E)=(L,h)$
as in Theorem \ref{thm_MainThm}. 
 \end{crl}   
 
 Indeed, by \eqref{eq:4.1} and  $iR^{\det}= -2 \omega_{\D^*}$, 
 the hypotheses of Proposition \ref{prop_SpectralGap} clearly 
 hold on $\D^*$. 
 Combining moreover with the condition
 ($\beta$) on $\Sigma$, we know that
 $iR^L\geq \varepsilon \omega_\Sigma$ and 
$iR^{\det}\geq  -C\omega_{\Sigma}$ on $\Sigma$,
for some $C>0$.  
We can thus apply Proposition \ref{prop_SpectralGap}
to get Corollary \ref{t5.2} in both situations.
 
\smallskip
  We assume here, without loss of generality, that the 
  puncture divisor $D$ in $\Sigma$ is reduced to one point $a$. 
  Let $\mathfrak{e}$ be the holomorphic frame of $L$ near $a$ 
  corresponding to the trivialization in the condition ($\alpha$).
  
   By the assumption ($\alpha$), ($\beta$), under our trivialization
  $\mathfrak{e}$ of $L$ on the coordinate $z$ on 
  $\D^*_{r}$ for some $0<r<e^{-1}$, we have the identification of
  the geometric data
  \begin{align}\label{eq:6.1}
(\Sigma,\omega_{\Sigma}, L,h)|_{\D^*_{r}}
= (\D^*,\omega_{\D^*}, \C, h_{\D^*})|_{\D^*_{r}}.
\end{align}

  We set: 
   \begin{itemize}
    \item  $F$ is the normalized Fourier transform of a smooth 
    cut-off function as in \cite[\S 4.1]{mm}, namely 
           \begin{equation}\label{eq:5.1}
F(u) = \Big(\int_{\R} f(v)\,dv\Big)^{-1}\int_{\R}e^{ivu}f(v)\,dv
           \end{equation}
with $f:\R\to[0,1]$ a smooth even function such that
$f(v)=1$ if $|v|\leq\epsilon/2$ and $f(v)=0$ if
$|v|\geq \epsilon$ for $\epsilon>0$. Thus $F$ is an even
function in the Schwartz space $\mathscr{S}(\R)$ with $F(0) = 1$.
Let $\widetilde{F}$ be the function satisfying 
$\widetilde{F}(u^2)=F(u)$ for all $u\in\R$.
 We consider the function
          \begin{equation}\label{eq:5.2}
          \phi_p:\R\longrightarrow\R, \:\:u\longmapsto
  \mathds{1}_{[cp,+\infty)}(|u|)\widetilde{F}(u)
          \end{equation}
          where $c>0$ is defined in \eqref{eqn_SpectralGap};
let  $K_p:=\phi_p(\Box_p)$ 
and let $K_p(\LargerCdot,\LargerCdot)$ be the associated kernel; 
we denote by $f_p(\LargerCdot,\LargerCdot)$ the function 
associated to $K_p(\LargerCdot,\LargerCdot)$ 
via the doubled trivialization around $a$
used above; for $x\in \D^*_r$, we set $f_{p,x}$ for 
the one-variable function $y\mapsto f_p(x, y)$; then 
    \begin{align}\label{eq:5.3}
    K_p(x,y)= f_p(x,y) \mathfrak{e}^p(x)\otimes
    (\mathfrak{e}^p(y))^*,
      \end{align}
and $(\mathfrak{e}^p(y))^*$ is the metric dual of 
 $\mathfrak{e}^p(y)$ with respect to $h^p$, that is,
 $\mathfrak{e}^p(y)^*\cdot\mathfrak{e}^p(y)
 =|\mathfrak{e}^p(y)|_{h^p}^2$.
    \item $\chi$ a cut-off function as in \eqref{e:chi}; 
    \item $\rho:\Sigma\to [1,+\infty)$ is a smooth function 
    such that $\rho(z)=\big|\!\log(|z|^2)\big|$  on $\D^*_r$.
   \end{itemize}

\begin{prop}\label{t5.3}  For any $\ell, m \geq 0$, 
   $\gamma>\frac{1}{2}$, there exists 
   $C_{\ell,m,\gamma}>0$ such that for any $p>1$, we have 
\begin{equation}\label{eq:5.4}
      \big\| \rho(x)^{-\gamma}\rho(y)^{-\gamma} K_p(x,y)
      \big \|_{C^m(h^p)}
      \leq C_{\ell,m,\gamma}p^{-\ell},
    \end{equation}
    in the sense of \eqref{eq:2.13}.
   \end{prop}
   \textit{Proof. --- } 
   We proceed as follows. 
   Take $p\gg1$, and pick ${\bf a}$ and $b$ two real parameters to 
   be determined later; then by definition of $f_p$ and $K_p$ :
    \begin{align}\label{eq:5.5}\begin{split}
     &\big\|\rho(x)^{{\bf a}}\rho(y)^{b}\chi(x)\chi(y)f_p(x,y)
     \big\|^2_{{\bL}^2(\Sigma\times\Sigma)}                           \\
 &= \int\limits_{x\in \D^*} \rho(x)^{2{\bf a}}\chi(x)^2  
 \omega_{\Sigma}(x)
 \int\limits_{y\in \D^*}\rho(y)^{2b} \chi(y)^2\big|f_p(x,y)\big|^2\,
 \omega_{\Sigma}(y)                              \\
 &= \int\limits_{x\in \D^*}\!\!\!  
 \rho(x)^{2{\bf a}-p}\chi(x)^2     
 \bigg\langle \int\limits_{y\in \D^*} 
 \!\!\!\big\langle K_p(x,y),\rho(y)^{2b-p}
 \chi(y)^2f_p(x,y)\mathfrak{e}^p(y)\big\rangle_{h_{\D^*}^p}\!
 \omega_{\D^*}(y), \mathfrak{e}^p(x)\bigg\rangle_{h_{\D^*}^p}
 \!\!\!\omega_{\D^*}(x) \\
 &= \int\limits_{x\in \D^*}  \rho(x)^{2{\bf a}-p}\chi(x)^2 
\Big\langle K_p(\rho^{2b-p}\chi^2 f_{p,x} \mathfrak{e}^p),
\mathfrak{e}^p\Big\rangle_{h_{\D^*}^p}\!(x)\, \omega_{\D^*}(x).
\end{split}    \end{align}
   As $|\mathfrak{e}^p|_{h_{\D^*}^p}=\rho^{p/2}$, 
   we have from (\ref{eq:5.5}), 
    \begin{align}\label{eq:5.6}\begin{split}
     \big\|\rho(x)^{{\bf a}}\rho(y)^{b}\chi(x)& \chi(y)f_{p}(x,y)
     \big\|^2_{{\bL}^2(\Sigma\times\Sigma)}    \\
 &= \int_{x\in \D^*}  \rho(x)^{2{\bf a}-\frac{p}{2}}\chi(x)^2 
\big|K_p(\rho^{2b-p}\chi^2 f_{p,x} \mathfrak{e}^p)(x)
\big|_{h_{\D^*}^p} \omega_{\D^*}(x).
 \end{split}   \end{align}
   Now for all $x\in \D^*_{2r/3}$, 
    \begin{align}\label{eq:5.7}\begin{split}
     \big|K_p(\rho^{2b-p}\chi^2 f_{p,x} \mathfrak{e}^p)(x)
     \big|_{h_{\D^*}^p}
&\leq \rho(x)^{1/2}\big\|\rho^{-1/2}K_p(\rho^{2b-p}\chi^2 
f_{p,x} \mathfrak{e}^p)\big\|_{C^0(\D^*_{2r/3})}                       \\
 &=   \rho(x)^{1/2} \big\|\rho^{-1}|K_p(\rho^{2b-p}\chi^2 f_{p,x}
 \mathfrak{e}^p)|_{h_{\D^*}^p}^2\big\|_{C^0(\D^*_{2r/3})}^{1/2}, 
  \end{split}     \end{align}
   and thus: 
    \begin{align}\label{eq:5.8}\begin{split}
\big\|\rho(x)^{{\bf a}}\rho(y)^{b}& \chi(x)\chi(y)f_p(x,y)
\big\|^2_{{\bL}^2(\Sigma\times\Sigma)}            \\
\leq &\int_{x\in \D^*}  \rho(x)^{2{\bf a}+\frac{1-p}{2}}\chi(x)^2        
\big\|\rho^{-1}|K_p(\rho^{2b-p}\chi^2 f_{p,x} 
\mathfrak{e}^p)|_{h_{\D^*}^p}^2
\big\|_{C^0(\D^*_{2r/3})}^{1/2}\omega_{\D^*}(x). 
  \end{split}  \end{align}
   We momentarily set $g(z)=|K_p(\rho^{2b-p}\chi^2 f_{p,x}
   \mathfrak{e}^p)|_{h_{\D^*}^p}^2(z)$, 
   and use the embedding
   ${\bL}^{1,3}_{\rm wtd}(\Sigma, \omega_{\Sigma})
   \longhookrightarrow C^0(\Sigma, \omega_{\Sigma})$
(cf. Lemma \ref{t2.4} a)),  that gives, taking supports into account: 
    \begin{align}\label{eq:5.10}\begin{split}
     \big\|\rho^{-1}g\big\|_{C^0(\D^*_{2r/3})} 
 &\leq c_0 \int_{\D^*_{r}} \rho\big(\big|\rho^{-1}g\big|
 +\ldots +\big|(\nabla^\Sigma)^3(\rho^{-1}g)\big|\big) 
 \omega_{\Sigma}  \\ 
         &\leq A^3 c_0\int_{\D^*_{r}} \big(\big|g\big|
+\ldots   +\big|(\nabla^\Sigma)^3g\big|\big)  \omega_{\Sigma}                                      
  \end{split}  \end{align}
with $A=1+\|d\log\rho\|_{C^2(\D^*_{r}, \omega_{\D^*})}$,  
which is finite by $\log\rho = \log\big(-\log(|z|^2)\big) = t$
and \eqref{eq:4.10}.
   Moreover, 
    \begin{align}\label{eq:5.11}\begin{split}
     \big|g\big|+\ldots+&\big|(\nabla^\Sigma)^3g\big|     \\
          &\leq C \big(|K_p(\rho^{2b-p}\chi^2 f_x 
	  \mathfrak{e}^p)|_{h_{\D^*}^p}^2
         +\ldots 
+|(\nabla^{p,\Sigma})^3K_p(\rho^{2b-p}\chi^2 f_{p,x} 
\mathfrak{e}^p)|_{h_{\D^*}^p}^2\big).
   \end{split} \end{align}
 We obtain from \eqref{eq:5.10}, \eqref{eq:5.11}, 
    \begin{align}\label{eq:5.12}
     \big\|\rho^{-1}g\big\|_{C^0(\D^*_{2r/3})}
\leq C \big\|K_p(\rho^{2b-p}\chi^2 f_{p,x} 
\mathfrak{e}^p)\big\|_{\bL^{2, 3}_p(h)}^2. 
    \end{align}
By Propositions \ref{prop_EstimateSigma} 
and \ref{prop_SpectralGap}, (\ref{eq:5.2})
and since $\widetilde{F}\in \mathscr{S}(\R)$,
we have that for any fixed $\ell$, there exists $C_{\ell}>0$ such that
for any $x\in \D^*_{r}, p\in \N^*$, 
 \begin{equation} \label{eq:5.12a} 
   \big\|K_p(\rho^{2b-p}\chi^2 f_{p,x} 
   \mathfrak{e}^p)\big\|_{\bL^{2, 3}_p(h)} 
  \leq C_{\ell}p^{-\ell} \big\| \rho^{2b-p}\chi^2 f_{p,x} 
\mathfrak{e}^p\big\|_{{\bL}^2_{p}(h)},
\end{equation} 
and thus \eqref{eq:5.8}, \eqref{eq:5.12}, \eqref{eq:5.12a} yield
    \begin{align}\label{eq:5.13}\begin{split}
     \big\|&\rho(x)^{\bf a}\rho(y)^{b}\chi(x)\chi(y)f_p(x,y)
     \big\|^2_{{\bL}^2(\Sigma\times\Sigma)}                  \\
      &\leq C_{\ell}\,p^{-\ell}\int_{x\in\D^*} \chi(x)^2 
	     \rho(x)^{2{\bf a}+\frac{1-p}{2}} 
 \big\| \rho^{2b-p}\chi^2 f_{p,x} \mathfrak{e}^p
 \big\|_{{\bL}^2_{p}(h)}  \omega_{\D^*}(x) \\
             &\leq C_{\ell}\,p^{-\ell}
 \bigg(\int_{\D^*} \chi^2\rho^{2{\bf a}-p+1}\omega_{\D^*}
 \bigg)^{\!\frac{1}{2}}
                   \bigg(\int_{x\in\D^*} \chi(x)^2 \rho(x)^{2{\bf a}}
\big\| \rho^{2b-p}\chi^2 f_{p,x} \mathfrak{e}^p
\big\|_{{\bL}^2_{p}(h)}^2
\omega_{\D^*}(x)\bigg)^{\!\frac{1}{2}}          \\
             &\leq C_{\ell}\,p^{-\ell}
 \bigg(\int_{\D^*} \chi^2\rho^{2{\bf a}-p+1}\omega_{\D^*}
 \bigg)^{\!\frac{1}{2}}
\,\big\|\rho(x)^{{\bf a}}\rho(y)^{2b-\frac{p}{2}}\chi(x)\chi(y)
f_p(x,y)\big\|_{{\bL}^2(\Sigma\times\Sigma)}, 
  \end{split}  \end{align}
   since 
    \begin{align*}
     \int_{x\in\D^*} \chi(x)^2 &\rho(x)^{2{\bf a}}\big\| \rho^{2b-p}
 \chi^2 f_{p,x} \mathfrak{e}^p\big\|_{{\bL}^2_{p}(h)}^2 
 \omega_{\D^*}(x)  \\ 
          =& \int_{x\in\D^*} \chi(x)^2 \rho(x)^{2{\bf a}} 
\omega_{\D^*}(x)
 \int_{y\in\D^*}\rho(y)^{4b-p} \chi(y)^4 |f_p(x,y)|^2 
 \omega_{\D^*}(y)                             \\
\leq& \int_{x\in\D^*} \chi(x)^2 \rho(x)^{2{\bf a}} 
\omega_{\D^*}(x)
\int_{y\in\D^*}\rho(y)^{4b-p} \chi(y)^2 |f_p(x,y)|^2 
\omega_{\D^*}(y)                             \\
          =& \big\|\rho(x)^{{\bf a}}\rho(y)^{2b-\frac{p}{2}}\chi(x)
 \chi(y)f_p(x,y)\big\|_{{\bL}^2(\Sigma\times\Sigma)}^2.
    \end{align*}
   Moreover, 
   $\int_{\D^*} \chi^2\rho^{2{\mathbf a}-p+1}\omega_{\D^*}$ 
   is finite as soon as ${\mathbf a}<\frac{p}{2}$. 
   Fixing ${\mathbf a}=\frac{p}{2}-\delta$, $\delta > 0$, 
   $\int_{\D^*} \chi^2\rho^{2{\bf a}-p+1}\omega_{\D^*}
         = \int_{\D^*} \chi^2\rho^{1-2\delta}\omega_{\D^*}
   =      \int_{\D^*} \chi^2 e^{-2\delta t}dt\,d\theta
         <\infty$ 
   is independent on $p$, 
   and consequently, the previous inequality reads: 
    \begin{equation}\label{eq:5.15}
     \big\|\rho(x)^{\frac{p}{2}-\delta}\rho(y)^{b}
     \chi(x)\chi(y)f_p(x,y)\big\|_{{\bL}^2(\Sigma\times\Sigma)} 
     \leq C_{\ell,\delta}\, p^{-\ell}
    \end{equation}
for all $b\leq \frac{p}{2}$ by using 
$\rho(y)^{2b-\frac{p}{2}}\leq \rho(y)^{b}$, 
and with $C_{\ell,\delta}$ independent 
of $b$, hence in particular for $p\geq1$,
    \begin{equation}\label{eq:5.16}
     \big\|\rho(x)^{\frac{p}{2}-\delta}\rho(y)^{\frac{p}{2}
 -\delta}\chi(x)\chi(y)f_p(x,y)\big\|_{{\bL}^2(\Sigma\times\Sigma)} 
   \leq C_{\ell,\delta}p^{-\ell}.
    \end{equation}
   With the same techniques, we extend this estimate to 
   higher orders, namely for any $k\geq 0$,
 \begin{equation}\label{eq:5.17}
 \big\|\big(\rho(x)\rho(y)\big)^{\frac{p}{2} - \delta}
 \chi(x)\chi(y)f_p(x,y)\big\|_{\bL^{2,k}(\Sigma\times\Sigma)}
 \leq C_{\ell,k, \delta}p^{-\ell},  
\end{equation}

Observe that by \eqref{eqn_Cmnorm} and  \eqref{eq:4.10},
$d\log \rho = dt$ is $C^m(\Sigma,\omega_{\Sigma})$ bounded,
 and the factor $\rho$ in the definition of 
 $\bL^{2,k}_{\rm wtd}$, thus \eqref{eq:5.17} implies    
    \begin{equation}\label{eq:5.18}
     \big\|\big(\rho(x)\rho(y)\big)^{\frac{p}{2}
     -\gamma}\chi(x)\chi(y)f_p(x,y)\big\|_{\bL^{2,k}_{\rm wtd}
     (\Sigma\times\Sigma)}
  \leq C_{\ell,k,\gamma}p^{-\ell},  
    \end{equation}
   with $\gamma=\delta+\frac{1}{2}>\frac{1}{2}$. 
   By Lemma \ref{t2.4} b) and \eqref{eq:5.18},
   for every $\ell, m\geq0$, $\gamma>\frac{1}{2}$
   there exists $C_{\ell,m,\gamma}>0$ such that for all $p\geq1$
   we have
        \begin{equation}\label{eq:5.19}
     \big\|\big(\rho(x)\rho(y)\big)^{\frac{p}{2}
     -\gamma}\chi(x)\chi(y)f_p(x,y)
     \big\|_{C^m(\Sigma\times\Sigma)}
   \leq C_{\ell,m,\gamma}p^{-\ell}, 
    \end{equation}
   which can be rewritten in the sense of  \eqref{eq:2.13} as
    \begin{equation}\label{eq:5.20}
     \big\|\rho(x)^{-\gamma}\rho(y)^{-\gamma}\chi(x)\chi(y)
     K_p(x,y)\big\|_{C^m(h^p)}\leq C_{\ell,m,\gamma}p^{-\ell}.  
    \end{equation} 
   Such an estimate is already well-known far from 
   $D_a=(\{a\}\times\Sigma)+(\Sigma\times\{a\})$ in 
   $\Sigma\times\Sigma$, 
   where the weights $\rho$ can be omitted,
   cf. \cite[\S 6.1]{mm}.  
   We moreover prove the analogous estimates on 
   \[\big(\D^*_{r/2}\times(\Sigma\smallsetminus\D_{r/3})\big)
   \cup \big((\Sigma\smallsetminus\D_{r/3})\times\D^*_{r/2}\big)\] 
   along the same lines. 
   We thus come to the conclusion that
    \eqref{eq:5.4} holds.
 \cqfd
 
 ~

 Estimates \eqref{eq:5.4} might look a bit disappointing, 
 as these are stated with 
 \textit{negative} weights;  so far, this does not even tell us that
 $K_p(x,y)$ is bounded \textit{near} the divisor 
 $(\{a\}\times\Sigma)+(\Sigma\times\{a\})$ 
 in $\overline{\Sigma}\times\overline{\Sigma}$
 for the product Hermitian norm, 
 which contrasts with our knowledge that the Bergman kernels of 
 $L^p$ do vanish \textit{along} this divisor.  
 
 We see in next part that this rough estimate, 
 together with the vanishing property, suffice however to estimate 
 sharply Bergman kernels on the whole $\Sigma\times\Sigma$.  
 
 \begin{rmk} \label{t5.4}
  The embedding ${\bL}^{1,3}_{\rm wtd}\longhookrightarrow C^0$ 
  does not hold on the whole $\D^*$. 
  Now Proposition \ref{t5.3} still holds on $\D^*\times \D^*$ 
  near $\{0\}+\{0\}$
  and, more generally, far from $\partial\D\times \partial\D$, 
  as functions with supports in $\D^*_{r}$ or $\D^*_{2r/3}$ 
  in $\D^*$ can be thought of as functions around the punctures
  in $\Sigma$, to which Lemma \ref{t2.4} applies. 
  More precisely, let $K^{\D^*}_p(x,y)$ be the kernel of 
$\phi_p(\Box_p^{\D^*})$ on $\D^*$ with respect to 
$\omega_{\D^*}$; then for all $0<r<1$ and all $\ell,m\geq 0$, 
and $\gamma>\frac{1}{2}$, there exists 
  $C_{\ell,m,\gamma}=C_{\ell,m,\gamma}(r)>0$ such that for 
  all $p\geq1$,
    \begin{equation}\label{eq:5.25}
\big\|\rho(x)^{-\gamma}\rho(y)^{-\gamma}\chi(x)\chi(y) 
K^{\D^*}_p(x,y)
\big\|_{C^m(\mathbb{D}^*_r\times\mathbb{D}^*_r)}
            \leq C_{\ell,m,\gamma}p^{-\ell}.
    \end{equation}
 \end{rmk}

 \section{Proofs of the main results
}
\label{S:ibke}
In this section, we will establish the results stated in the Introduction. 
We first complete the proofs of Theorems \ref{t6.1}, \ref{thm_MainThm} 
and Corollary \ref{crl_CrlIntro}, 
then that of Corollary \ref{crl_CrlIntro1}, 
and finally give the details needed to establish Theorems \ref{t2.5} 
and \ref{t0.6}. 

Theorems \ref{t6.1}, \ref{thm_MainThm}  
and Corollary \ref{crl_CrlIntro} will be a consequence of
Proposition \ref{t5.3}.
  The principal idea is to
  combine Proposition \ref{t5.3} and the \textit{holomorphicity} 
  of sections associated to Bergman kernels, 
  together with the fact that for $p\geq2$, 
  ${\bL}^2$ holomorphic sections of $(L^p,h^p)$ over 
  $(\Sigma,\omega_{\Sigma})$ 
  \textit{vanish at $D$}, 
  and similarly for the holomorphic sections in 
  $\bL^{2}\big(\D^*,\omega_{\D^*}, \C, h_{\D^*}^p\big)$, 
  which vanish at $0$, as already noticed. 

Let us fix a point $a\in D$ and work around this point,
in coordinates centered at $a$.
 \comment{
  By the assumption ($\alpha$), ($\beta$), under our trivialization
  $\mathfrak{e}$ of $L$ on the coordinate $z$ on 
  $\D^*_{r}$ for some $0<r<e^{-1}$, we have the identification of
  the geometric data
  \begin{align}\label{eq:6.1}
(\Sigma,\omega_{\Sigma}, L,h)\mid_{\,\D^*_{r}}\,
= (\D^*,\omega_{\D^*}, \C, h_{\D^*}^p)\mid_{\,\D^*_{r}}.
\end{align}
}
For $x,y\in \D^*_{r}$,  under our identification \eqref{eq:6.1} 
and convention after (\ref{eq:3.7}), we write 
   \begin{equation}\label{eq:6.3}\begin{split}
  & {B}^{\D^*}_p(x,y) 
    = \big|\!\log(|y|^2)\big|^{p} \beta^{\D^*}_p(x,y); \\
 &   {B}_p(x,y) 
    = \big|\!\log(|y|^2)\big|^{p} \beta^{\Sigma}_p(x,y),
   \end{split}\end{equation}
  where, by \eqref{eq:3.7}, $\beta^{\D^*}_p$ is 
  a holomorphic function of $x$ and $\overline{y}$, namely 
   \begin{equation}\label{eq:6.4}
    \beta^{\D^*}_p(x,y) = \frac{1}{2\pi(p-2)!} 
    \sum_{\ell=1}^{\infty} \ell^{p-1} x^{\ell}\overline{y}^{\ell}, 
   \end{equation}
  which vanishes along 
  $\{x=0\}+\{y=0\} = \{xy=0\}\subset \D\times\D$. 
By Remark \ref{t3.2}, \eqref{eq:2.2} and the convention 
after \eqref{eq:3.7}, $\beta^{\Sigma}_p$ is a holomorphic function
of $x$ and $\overline{y}$ vanishing along 
  $\{x=0\}+\{y=0\}\subset \D_{r}\times \D_{r}$, and 
  \begin{align}\label{eq:6.5}
\Big |( B_{p} - {B}^{\D^*}_p)(x,y) \Big |_{h^p} = 
 \big|\!\log(|x|^2)\big|^{p/2}\big|\!\log(|y|^2)\big|^{p/2}
\Big |(\beta^{\D^*}_p - \beta^{\Sigma}_p)(x,y) \Big |.
\end{align}
   \textit{Proof of Theorem \ref{t6.1}. --- }
   By \eqref{eq:5.2}, we have
   \begin{align}\label{eq:6.9}\begin{split}
    \widetilde{F}\big(\Box_p\big)(x,y) - {B}_p(x, y) 
    &= \phi_p\big(\Box_p\big)(x,y) = K_p(x,y), \\
    \widetilde{F}\big(\Box_p^{\D^*}\big)(x,y) - {B}_p^{\D^*}(x, y)    
    &= \phi_p\big(\Box_p^{\D^*}\big)(x,y) = K_p^{\D^*}(x,y).
  \end{split} \end{align}
By the finite propagation speed for the wave operators
\cite[Theorem D.2.1]{mm}, we have 
\begin{equation}\label{eq:6.4a}
\begin{split}
&\supp \widetilde{F}(\square_{p})(x,\LargerCdot)
\subset B\Big(x, \frac{\epsilon}{\sqrt{2}}\Big)
\text{ and } \widetilde{F}(\square_{p})(x,\LargerCdot)\\
&\text{depends only on the restriction of 
$\square_{p}$ to $B\Big(x, \frac{\epsilon}{\sqrt{2}}\Big)$}.
\end{split}
\end{equation}
Here $B(x, \frac{\epsilon}{\sqrt{2}})$ is the geodesic ball 
with center at $x$ and radius of $\frac{\epsilon}{\sqrt{2}}$
for $\omega_{\mathbb{D}^*}$.

Thus from \eqref{eq:6.1} and (\ref{eq:6.4a}) (here we fix
$\epsilon>0$ such that $\epsilon\leq d_{\D^*}(\partial D_{r},
\partial D_{r/2})$), we have
\begin{equation}\label{eqn_FisF}
\widetilde{F}\big(\Box_p\big)(x,y)
=\widetilde{F}\big(\Box_p^{\D^*}\big)(x,y)
\:\:\text{ for all $x,y\in\D^*_{r/2}$}\,.
\end{equation}
 We have from \eqref{eq:6.9} and \eqref{eqn_FisF},
\begin{equation}\label{eq:6.5a}
{B}_p(x, y)-{B}^{\D^*}_p(x, y)=K_p^{\D^*}(x,y)-K_p(x,y)
\:\:\text{ for all $x,y\in\D^*_{r/2}$}\,.
\end{equation}

   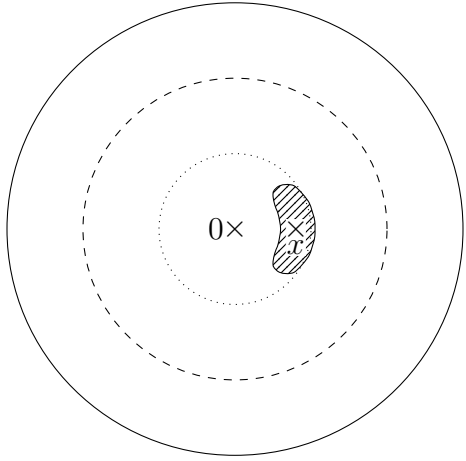
\begin{figure}
   \begin{centering}
   \begin{tikzpicture}[scale=.5]
    \draw (0,0) circle (6) ;
    \draw (0,0) node {$\times$} node [left] {$0$};
    \draw[black,pattern=north east lines] (1.2,0) ..controls +(0,-.1) 
    and +(.05,.1).. (1.1,-.5) ..controls +(-.05,-.1) and +(0,.1).. (1,-.9) 
  ..controls +(0,-.1) and +(-.1,.1).. (1.1, -1.1) ..controls +(.1,-.1) 
  and +(-.1,-.05).. (1.55,-1.15) 
    ..controls +(.1,.05) and +(-.05,-.1).. (1.95,-.7) ..controls +(.05, .1) 
    and +(-.01,-.1).. (2.1, -.2)
      ..controls +(.01,.1) and +(.01,-.1).. (2.1, .2) ..controls +(-.01,.1) 
      and +(.05, -.1).. (1.95,.7) 
..controls +(-.05,.1) and +(.1,-.05).. (1.55,1.15) ..controls +(-.1,.05)
and +(.1,.1).. (1.1, 1.1) 
   ..controls +(-.1,-.1) and +(0,.1).. (1,.9) ..controls +(0,-.1) 
   and +(-.05,.1).. (1.1,.5) 
  ..controls +(.05,-.1) and +(0,.1).. (1.2,0);
    \draw[dashed] (0,0) circle (4) ;
    \draw[dotted] (0,0) circle (2) ;
    \draw[white,fill=white]  (1.35,-.7) -- (1.85,-.7) -- (1.85,.25) -- 
    (1.35,.25) -- (1.35,-.7);
    \draw (1.6,0) node {$\times$} node [below] {$x$};  
   \end{tikzpicture}
   \caption{In \eqref{eqn_FisF}, $\tilde{F}(\Box_p)(x,\point)$ 
   depends only 
   on the restriction of $\Box_p$ \\ to a \textit{geodesic}
   ball of radius $\frac{\epsilon}{\sqrt{2}}$ 
            around $x$, and vanishes outside this ball \\ 
            (the circles $\{|z|=r/2\}$ and $\{|z|=r\}$
	    are at Poincar\'{e} distance \\
            $\frac{1}{\sqrt{2}}\log\big(1-\frac{\log 2}{\log r}\big)$ 
             from each other).}
   \end{centering}
   \end{figure}
Note that by (\ref{eq:4.5}) and (\ref{eq:4.10}), 
\begin{align}\label{eq:6.8}
\Big ||\log |z|^2|^p\Big|_{C^k} = \Big |e^{tp}\Big|_{C^k} 
\leq C_k p^k e^{tp}.
\end{align}
By (\ref{eqn_Cmnorm}), (\ref{eq:2.13}), (\ref{eq:4.5}), 
(\ref{eq:4.10}), (\ref{eqn_nablap}),
(\ref{eq:6.1}) and (\ref{eq:6.8}),
we infer that for any $k\in \N$ there exists $C>0$
such that for any $p\in \N^*$, $x,y\in\D^*_{r/2}$, we have
\begin{align}\label{eq:6.10}\begin{split}
&\big|B^{\D^*}_p-B_p\big|_{C^k(h^p)} (x,y)
\leq C p^k  \big|\!\log(|x|^2)
    \big|^{p/2}\big|\!\log(|y|^2)\big|^{p/2}
   \big|\beta^{\D^*}_p-\beta^{\Sigma}_p\big|_{C^k} (x,y),\\
  & \big|\beta^{\D^*}_p-\beta^{\Sigma}_p\big|_{C^k} (x,y)
\leq C p^k  \big|\!\log(|x|^2)
    \big|^{-p/2}\big|\!\log(|y|^2)\big|^{-p/2}
   \big|B^{\D^*}_p- B_p\big|_{C^k(h^p)} (x,y).
\end{split}\end{align}
By Proposition \ref{t5.3}, (\ref{eq:6.5a}) and \eqref{eq:6.10},
we get for $x,y\in \D^*_{r/2}$, 
   \begin{equation}   \label{eqn_estmteb}
    \big|\beta^{\D^*}_p-\beta^{\Sigma}_p\big| _{C^k} (x,y)
    \leq C_{\ell,\gamma}p^{-\ell} \big|\!\log(|x|^2)
    \big|^{\gamma-p/2}\big|\!\log(|y|^2)\big|^{\gamma-p/2}.
   \end{equation}
\comment{  
 but we are interested in deriving estimates on quantities such as
   \begin{equation*}
    \beta^{\Sigma}_p(x,y)\big|\!\log(|x|^2)\big|^{p/2}
    \big|\!\log(|y|^2)\big|^{p/2}
   \end{equation*}
  (for instance, 
   \begin{equation*}
    \beta^{\Sigma}_p(x,x)\big|\!\log(|x|^2)\big|^{p}
    =\big|{B}_p^{\Sigma}(x,x)\big|_{h^{\Sigma}} 
    = \sup_{\|\sigma\|_{H^p_{(2)}(\Sigma)}=1} |\sigma(x)|^2_{h^p} 
   \end{equation*}
   for $x$ near $a$ in $\Sigma$) 
  from the knowledge of the model ${B}^{\D^*}_p(x,x)$, 
  in which the previous estimate does not allow to conclude. 
  }
  Our task is thus to refine estimate \eqref{eqn_estmteb} 
  by working directly on $\beta^{\D^*}_p-\beta^{\Sigma}_p$.
  We set: 
   \begin{equation}\label{eq:6.6a}
    \beta^{\D^*,\Sigma}_p = \beta^{\D^*}_p-\beta^{\Sigma}_p.   
   \end{equation}
  As $\beta^{\D^*,\Sigma}_p$ is a holomorphic function of $x$ 
  and $\overline{y}$ vanishing along $\{xy=0\}$, 
  one can write 
   \begin{equation}\label{eq:6.7a}
    \beta^{\D^*,\Sigma}_p (x,y)= x\overline{y}g_p(x,y) 
   \end{equation}
  for some smooth $g_p$, holomorphic in $x$ and $\overline{y}$. 
  Now for $\vareps>0$ small, 
   \begin{equation}\label{eq:6.8a}
    \sup_{|x|,|y|\leq\vareps^{1/2}} |g_p(x,y)| 
    =  \sup_{|x|=|y|=\vareps^{1/2}} |g_p(x,y)| 
    = \vareps^{-1} \sup_{|x|=|y|=\vareps^{1/2}} 
    \big|\beta^{\D^*,\Sigma}_p (x,y)\big|
   \end{equation}
  by holomorphicity of $g_p$. 
 From (\ref{eqn_estmteb}) and (\ref{eq:6.8a}), we get: 
   \begin{equation}\label{eq:6.9a}
    \sup_{|x|,|y|\leq\vareps^{1/2}} |g_p(x,y)|
    \leq C_{\ell,\gamma}p^{-\ell}\vareps^{-1}
    |\log\vareps|^{-p+2\gamma}. 
   \end{equation}
 From (\ref{eq:6.7a}), (\ref{eq:6.9a}),
for $|x|, |y|\leq e^{-p}$, we have
   \begin{equation}   \label{eqn_ineqb1}\begin{split}
    \big|\beta^{\D^*,\Sigma}_p (x,y)\big| 
    & \leq  |x||y| \sup_{|x|,|y|\leq e^{-p}  } |g_p(x,y)| \\
    & \leq C_{\ell,\gamma}p^{-\ell}e^{2p}(2p)^{-p+2\gamma}|x||y|.
  \end{split} \end{equation}
  Pick $\alpha\geq 0$. 
  As the function 
  $u\mapsto u\big(-\log(u^2)\big)^{\frac{p}{2}+\alpha}$
  is increasing on $(0,e^{-\frac{p}{2}-\alpha}]$, 
  hence on $(0,e^{-p}]$, thus
   \begin{equation}\label{eq:6.14}
    |x| \leq e^{-p}(2p)^{\frac{p}{2}
    +\alpha}\big|\!\log(|x|^2)\big|^{-\frac{p}{2}-\alpha}
 \quad \text{if   } |x|\leq e^{-p}.
   \end{equation}
  We thus convert \eqref{eqn_ineqb1} into: 
   \begin{equation}\label{eq:6.15}
    \big|\beta^{\D^*,\Sigma}_p (x,y)\big| 
                 \leq 2^{2(\gamma+\alpha)}C_{\ell,\gamma}
 p^{-\ell+2(\gamma+\alpha)}
  \big|\!\log(|x|^2)\big|^{-\frac{p}{2}-\alpha}\big|\!\log(|y|^2)
  \big|^{-\frac{p}{2}-\alpha}
   \end{equation}
  on $|x|,|y|\leq e^{-p}$; 
  up to increasing $\ell$ and adjusting the constant, we thus have 
   \begin{equation}\label{eq:6.16}
    \big|\beta^{\D^*,\Sigma}_p (x,y)\big| 
                 \leq C_{\ell,\alpha}p^{-\ell}
 \big|\!\log(|x|^2)\big|^{-\frac{p}{2}-\alpha}
 \big|\!\log(|y|^2)\big|^{-\frac{p}{2}-\alpha}
   \end{equation}
  on $\{|x|,|y|\leq e^{-p}\}$. 
  Notice that such an estimate holds on 
  $\{e^{-p}\leq |x|, |y|\leq r\leq e^{-1}\}$ as well, 
  since then $\big|\log(|x|^2)\big|,\big|\log(|y|^2)\big|\leq 2p$, 
  and thus \eqref{eqn_estmteb} gives for any $\beta\geq 0$,
   \begin{equation}\label{eq:6.17}
    \big|\beta^{\D^*,\Sigma}_p (x,y)\big| 
 \leq C_{\ell,\gamma}p^{-\ell}(2p)^{2(\gamma+\beta)}
 \big|\!\log(|x|^2)\big|^{-\beta-p/2}
 \big|\!\log(|y|^2)\big|^{-\beta-p/2}. 
   \end{equation}
  We conclude by the case 
  $|x|\leq e^{-p},\, e^{-p}\leq |y|\leq r\leq e^{-1}$;  
  fixing $y$, by the holomorphicity of $x$ 
  and \eqref{eqn_estmteb},
   \begin{align}\label{eq:6.18}\begin{split}
\sup_{|x|\leq e^{-p}}\Big|\frac{1}{x}
\beta_p^{\D^*,\Sigma}(x,y)\Big| 
 &= \sup_{|x| = e^{-p}}\Big|\frac{1}{x}\beta_p^{\D^*,\Sigma}(x,y)
 \Big|                         \\
 &= e^{p}\sup_{|x| = e^{-p}}\Big|\beta_p^{\D^*,\Sigma}(x,y)\Big|                               \\
&\leq e^pC_{\ell,\gamma}p^{-\ell}(2p)^{-\frac{p}{2}
+\gamma}\big|\!\log(|y|^2)\big|^{-\frac{p}{2}+\gamma}
  \end{split} \end{align}
  so that 
   \begin{equation}\label{eq:6.19}
    \big|\beta_p^{\D^*,\Sigma}(x,y)\big| 
    \leq e^pC_{\ell,\gamma}p^{-\ell}(2p)^{-\frac{p}{2}
    +\gamma}|x|\big|\!\log(|y|^2)\big|^{-\frac{p}{2}+\gamma}
   \end{equation}
  on $\{|x|\leq e^{-p},\, e^{-p}\leq |y|\leq r\leq e^{-1}\}$. 
  Now on this set, 
  $\big|\!\log(|y|^2)\big|^{\frac{p}{2}
  +\gamma}\leq (2p)^{\alpha+\gamma}\big|\!\log(|y|^2)
  \big|^{-\frac{p}{2}-\alpha}$ 
 by (\ref{eq:6.14}) and (\ref{eq:6.19}), we get
   \begin{equation}\label{eq:6.20}
\big|\beta_p^{\D^*,\Sigma}(x,y)\big| \leq 
2^{2(\alpha+\gamma)}C_{\ell,\gamma}p^{-\ell+2(\alpha+\gamma)}
\big|\!\log(|x|^2)\big|^{-\frac{p}{2}-\alpha}\big|\!\log(|y|^2)
\big|^{-\frac{p}{2}-\alpha}
   \end{equation}
  on $\{|x|\leq e^{-p},\, e^{-p}\leq |y|\leq r\leq e^{-1}\}$.  
  This holds on $\{|y|\leq e^{-p},\, e^{-p}\leq |x|\leq r\leq e^{-1}\}$ 
  as well by symmetry. 

  We get from (\ref{eq:6.16})---(\ref{eq:6.20}): 
    For all $\ell>0$, $\alpha\geq 0$, there exists 
    $C_{\ell,\alpha}>0$ such that for $p>1$,
    on $\D_{r/2}^* \times \D_{r/2}^*$,
     \begin{equation}\label{eq:6.21}
      \big|(\beta_p^{\D^*}-\beta_p^{\Sigma})(x,y)\big| 
            \leq C_{\ell,\alpha}p^{-\ell}
                 \big|\!\log(|x|^2)\big|^{-\frac{p}{2}-\alpha}
\big|\!\log(|y|^2)\big|^{-\frac{p}{2}-\alpha}.
     \end{equation}
   From (\ref{eq:6.5}) and (\ref{eq:6.21}) we get 
   (\ref{eq:6.6}) for $k=0$. 
   
   Observe that by (\ref{eq:4.3}), 
   (\ref{eq:4.11}), for any $k\in \N$, there exists 
   $C_{k}>0$ such that on $\D^*_{r}$,
   \begin{align}\label{eq:6.17a}
\Big| \log |z|^2\Big|_{C^k}\leq C_{k} \Big| \log |z|^2\Big|,\quad
|z|_{C^k} \leq C_{k} \Big| \log |z|^2\Big|^k |z|.
\end{align}
 From  (\ref{eq:6.7a}) and (\ref{eq:6.17a}) we have 
\begin{align}\label{eq:6.18a}
 \big|\beta_p^{\D^*, \Sigma}(x,y)\big|_{C^k}
            \leq C_{k} |x| |y| \sum_{k_{1}+k_{2}+k_{3}\leq k}
  \big|\!\log(|x|^2)\big|^{k_{1}}
\big|\!\log(|y|^2)\big|^{k_{2}} \Big| g_{p}(x,y)\Big|_{C^{k_{3}}}.
\end{align}
By (\ref{eq:4.4}) and (\ref{eq:4.10}) we get
\begin{multline}\label{eq:6.19a}
\Big| g_{p}(x,y)\Big|_{C^j}\\
\leq C \sum_{j_{1}+j_{2}\leq j} 
 \big|\!\log(|x|^2)\big|^{j_{1}}\big|\!\log(|y|^2)\big|^{j_{2}} 
 |x|^{\min (1,j_{1})} |y|^{\min (1,j_{2})}
\Big| \frac{\partial^{j_{1}+j_{2}}}{\partial x^{j_{1}} \partial 
\overline{y}^{j_{2}}}g_{p}(x,y)\Big|.
\end{multline}
Thus from (\ref{eq:6.18a}) and (\ref{eq:6.19a}), 
for any $k\in \N$, there exists 
   $C_{k}>0$ such that  for $x,y\in \D^*_{r}$, 
\begin{align}\label{eq:6.20a}
 \big|\beta_p^{\D^*, \Sigma}(x,y)\big|_{C^k}
 \leq C_{k} |x| |y|  \big|\!\log(|x|^2)\big|^{k}
\big|\!\log(|y|^2)\big|^{k}
\sum_{j_{1}+j_{2}\leq k} \Big| \frac{\partial^{j_{1}+j_{2}}}
{\partial x^{j_{1}} \partial 
\overline{y}^{j_{2}}}g_{p}(x,y)\Big|.
\end{align}
Now we can combine the argument for (\ref{eq:6.21}) and 
(\ref{eq:6.20a}) to get (\ref{eq:6.6}) for $k\geq 1$. 
The proof of Theorem \ref{t6.1} is complete.
   \cqfd
  \begin{crl}\label{crl6.2} For any $\ell, m\in \N$, $\varepsilon>0$, 
	  and  $\delta>0$, 
  there exists $C=C(\ell, m, \vareps, \delta)>0$ such that
  for all $p\in\N^*$, and 
  $x,y\in V_1\cup\ldots\cup V_N$, $d(x,y)>\varepsilon$,  
  with the local coordinate $x_{i}, y_{j}$,
  in the sense of (\ref{eq:2.13}), we have
   \begin{equation}    \label{eq6.20b}
     \Big | B_p(x_{i}, y_{j})\Big |_{C^m(h^{p})} \leq Cp^{-\ell}
     \, \big|\!\log(|x_{i}|^2)\big|^{-\delta}
	 \big|\!\log(|y_{j}|^2)\big|^{-\delta}.   
   \end{equation}
If $x\in V_1\cup\ldots\cup V_N$,
  with the local coordinate $x_{i}$, 
  $y\in \Sigma\smallsetminus (V_1\cup\ldots\cup V_N)$, 
and  $d(x,y)>\varepsilon$,  then
 \begin{equation}    \label{eq6.21b}
     \Big | B_p(x_{i}, y)\Big |_{C^m(h^{p})} 
	 + \Big | B_p(y, x_{i})\Big |_{C^m(h^{p})} \leq Cp^{-\ell}
     \, \big|\!\log(|x_{i}|^2)\big|^{-\delta}.   
   \end{equation}
  If $x, y\in \Sigma\smallsetminus (V_1\cup\ldots\cup V_N)$, 
  and $d(x,y)>\varepsilon$, then 
  $ | B_p(x, y) |_{C^m(h^{p})} \leq Cp^{-\ell}.$
\end{crl}
\textit{Proof. --- }   Pick $\vareps>0$ (small), and, if needed, replace 
the $\epsilon$ fixed in Section \ref{S:sgl} by this new $\vareps$. 
As $d(x,y)>\varepsilon$, by 
\eqref{eq:6.9} and \eqref{eq:6.4a}, we get 
   \begin{align}\label{eq:6.22b}
- {B}_p(x, y)  = K_p(x,y).
\end{align}
We explain now in detail how to establish \eqref{eq6.20b}.
The argument of \eqref{eq:6.10} implies
\begin{align}\label{eq:6.23b}\begin{split}
&\big|B_p\big|_{C^k(h^p)} (x,y)
\leq C p^k  \big|\!\log(|x|^2)
    \big|^{p/2}\big|\!\log(|y|^2)\big|^{p/2}
   \big|\beta^{\Sigma}_p\big|_{C^k} (x,y),\\
  & \big|\beta^{\Sigma}_p\big|_{C^k} (x,y)
\leq C p^k  \big|\!\log(|x|^2)
    \big|^{-p/2}\big|\!\log(|y|^2)\big|^{-p/2}
   \big|B_p\big|_{C^k(h^p)} (x,y).
\end{split}\end{align}
By Proposition \ref{t5.3}, and \eqref{eq:6.22b},
we get the analogue of \eqref{eqn_estmteb} in the situation 
of \eqref{eq6.20b}, that is, for any $\gamma>\frac{1}{2}$,
\begin{align}\label{eq:6.24b}
\big|\beta^{\Sigma}_p\big| _{C^k} (x,y)
    \leq C_{\ell,\gamma}p^{-\ell} \big|\!\log(|x|^2)
    \big|^{\gamma-p/2}\big|\!\log(|y|^2)\big|^{\gamma-p/2}.
\end{align}
By the argument of \eqref{eq:6.17}--\eqref{eq:6.20a} 
and \eqref{eq:6.23b}, we improve then \eqref{eq:6.24b}
to \eqref{eq6.20b}.
The same argument as above works for 
the situation of \eqref{eq6.21b}.
The last part of Corollary is from 
\eqref{eq:5.4}, \eqref{eq:6.9} and \eqref{eq:6.4a}.
\cqfd
\\

\smallskip
\noindent 
    \textit{Proof of Theorem \ref{thm_MainThm}. --- }
  This follows immediately by taking $x=y$ in Theorem \ref{t6.1}.
    \cqfd

\smallskip
\noindent 
\textit{Proof of Corollary \ref{crl_CrlIntro}. --- } 
   By  Theorem \ref{thm_MainThm}, Proposition \ref{t3.3} and 
   (\ref{eq:3.20}), we get
    \begin{equation}\label{eq:6.23}
      \sup_{|x|\leq r\leq e^{-1}} \big|{B}_p(x,x)
      \big|_{h^p} = \frac{p^{3/2}}{(2\pi)^{3/2}}+\mathcal{O}(p)
      \qquad\text{as }p\to +\infty.   
    \end{equation}
    Combining (\ref{e:bke}), (\ref{e:Pvar}) and (\ref{eq:6.23}), 
    we obtain the desired conclusion. 
    \cqfd
    \smallskip
    
We turn now to the proof of Corollary \ref{crl_CrlIntro1}.
Let $\overline\Sigma$ be a compact Riemann surface of genus 
$g$, let $D=\{a_1,\ldots,a_N\}\subset\overline\Sigma$ be 
a finite set and $\Sigma=\overline\Sigma\smallsetminus D$. 

The Uniformization Theorem, see 
\cite[Theorems IV.5.6, IV.6.3, IV.6.4, IV.8.6]{fk}, 
readily implies that conditions (i)-(iv) from the introduction are 
equivalent, taking into account that $\chi(\Sigma)=2-2g-N$ 
and the degree of $L$ equals $-\chi(\Sigma)$:
{}from \cite[Theorem IV.8.6]{fk} we see first that (i) and (iii) 
are equivalent, and combining \cite[Theorems IV.6.3, IV.6.4]{fk} we 
know that (iii) implies (ii), since the Riemann surfaces with universal 
covering the sphere or $\C$, are the sphere, $\C$, $\C^*$
or the tori. This means that (ii) implies (iii).
Finally, by \cite[p.\ 214]{GH}, (ii) and (iv) are equivalent.
  \begin{lem}\label{L:ab}
  Let $\overline\Sigma$ be a compact Riemann surface of genus 
  $g$ and
 $D=\{a_1,\ldots,a_N\}\subset\overline\Sigma$ 
 a finite set such that $2g-2+N>0$.
 Denote $\Sigma=\overline\Sigma\smallsetminus D$ and 
 $L=K_{\overline\Sigma}\otimes
 \mathscr{O}_{\overline{\Sigma}}(D)$.
 There exists a metric $\omega_\Sigma$ on $\Sigma$ 
 and a singular Hermitian metric $h$ on $L$, 
 such that $(\Sigma,\omega_{\Sigma})$
 and the formal square root of $(L,h)$ satisfy
 the conditions $(\alpha)$ and $(\beta)$.
  \end{lem}
 \prf.\textit{ --- } Since $\chi(\Sigma)=2-2g-N<0$, the universal 
 covering of $\Sigma$ is $\mathbb{H}$ and $\Sigma$ admits 
 a K\"ahler-Einstein metric
 $\omega_{\Sigma}$ of constant negative curvature $-4$, 
 induced by the Poincar\'e metric 
 $\omega_\mathbb{H}= \frac{idz\wedge d\overline{z}}
 {4 |{\rm Im} z|^2}$ on $\mathbb{H}$.
It is a classical fact that every $a\in D$ has a coordinate 
neighborhood
$\big(\overline{U}_a,z\big)$ in $\overline\Sigma$ such that in 
this coordinate 
$\omega_{\Sigma}$ is exactly given by $\omega_{\D^*}(z)$ 
on $U_a = \overline{U}_a\smallsetminus\{a\}$,
see e.g. \cite[p.\,79, (6.7)]{bor}.

Note that  $\omega_{\Sigma}$ extends to a closed 
strictly positive $(1,1)$-current 
$\omega_{\overline\Sigma}$ on $\overline\Sigma$.
\comment{
By \cite[Lemma 6.8]{CM11}, there exists a singular metric $h$
on $L$ over $\overline\Sigma$, which is smooth on $\Sigma$ 
and whose curvature current satisfies 
$iR^L=-\ric_{\omega_{\Sigma}} = \omega_{\overline\Sigma}$ 
on $\Sigma$.
If $\omega_{\Sigma}=\frac{i}{2}gdz\wedge d\overline{z}$ 
on $\overline{U}_a\smallsetminus\{a\}$ 
in the coordinate $(\overline{U}_a,z)$ 
near $a\in D$, the weight of $h$ on $\overline{U}_a$ has the 
form $\varphi=\frac12\log(|z|^2g)$.
As one can choose $z$ so that $g=2|z|^{-2}\log^{-2}(|z|^2)$ 
}
Let $h^{K_{\Sigma}}$ be the metric on $K_{\Sigma}$ induced by 
$\omega_{\Sigma}$. Then we have  
\begin{align}\label{eq:6.22}
|dz|_{h^{K_{\Sigma}}}^2= |z|^{2}\log^{2}(|z|^2)
\, \,  \text{in}\, \, \big(\overline{U}_a,z\big).
\end{align}
Let $\sigma$ be the canonical section of 
$\mathscr{O}_{\overline{\Sigma}}(D)$.
The singular metric $h^{\mathscr{O}_{\overline{\Sigma}}(D)}$
on $\mathscr{O}_{\overline{\Sigma}}(D)$ is defined by
$|\sigma|^2_{h^{\mathscr{O}_{\overline{\Sigma}}(D)}}=1$.
The isomorphism 
$$K_{\Sigma}\longrightarrow K_{\Sigma}\otimes 
\mathscr{O}_{\overline{\Sigma}}(D)|_{\Sigma}= L|_{\Sigma},
\quad s\longmapsto s\otimes \sigma$$
over $\Sigma$ and the metrics 
$h^{K_{\Sigma}}$ and 
$h^{\mathscr{O}_{\overline{\Sigma}}(D)}$ induce 
 the metric $h$ on $L|_{\Sigma}$. The 
curvature of the line bundle $(L|_{\Sigma}, h)$
is given by $-2 i\omega_{\Sigma}$. 
Since $\frac{\sigma}{z}$ is a holomorphic frame of 
$\mathscr{O}_{\overline{\Sigma}}(D)$ on $\overline{U}_a$,
$dz\otimes \frac{\sigma}{z}$ is a holomorphic frame of $L$ on
$\overline{U}_a$.
Then $|dz\otimes \frac{\sigma}{z}|_{h}^2 = (\log(|z|^2))^{2}$,
and thus $(\Sigma,\omega_{\Sigma})$ and 
the (formal) square root of $(L,h)$ satisfy 
conditions $(\alpha)$ and $(\beta)$.
\cqfd
\\
\smallskip

Let $\Gamma\cong\pi_1(\Sigma)$ be the group of 
deck transformations of the 
covering $\mathbb{H}\to\Sigma$.
Then $\Gamma\backslash\mathbb{H}\cong\Sigma$ has 
finite hyperbolic volume and 
$\Gamma$ is a Fuchsian group of the first kind
without elliptic elements. We denote by 
$\pi:\mathbb{H}\to\Gamma\backslash\mathbb{H}$
the canonical projection.

The space $\mathcal{M}_{2p}^\Gamma$ of 
$\Gamma$-modular forms 
of weight $2p$ is by definition the space 
of holomorphic functions $f\in\mathscr{O}(\mathbb{H})$
satisfying the functional equation
\begin{equation}\label{eq:6.24}
f(\gamma z)=(cz+d)^{2p}f(z),\quad z\in\mathbb{H},\:\: 
\gamma=\begin{pmatrix}a&b\\c&d\end{pmatrix}\in\Gamma,
\end{equation}
and which extend holomorphically to the cusps of $\Gamma$ (fixed
points of the parabolic elements).
If $f\in\mathscr{O}(\mathbb{H})$ satisfies
\eqref{eq:6.24}, then 
$fdz^{\otimes p}\in H^0(\mathbb{H},K_\mathbb{H}^p)$
descends to a holomorphic section $\Phi(f)$ of 
$H^0(\Sigma,K^p_\Sigma)\cong H^0(\Sigma,L^p)$. 
By \cite[Propositions\,3.3,\,3.4(b)]{Mum77}, $\Phi$ induces 
an isomorphism
$\Phi:\mathcal{M}_{2p}^\Gamma
\to H^0\big(\overline\Sigma,L^p\big)$.

The subspace of $\mathcal{M}_{2p}^\Gamma$ consisting of
modular forms vanishing at the cusps is called the 
space of \emph{cusp forms} (Spitzenformen) of weight $2p$ 
of $\Gamma$, denoted by $\mathcal{S}_{2p}^\Gamma$.
The space of cusps forms is endowed with 
the Petersson scalar product
\[
\langle f,g\rangle:=
\int_U f(z)\overline{g(z)}(2y)^{2p}\,dv_{\mathbb{H}}(z),
\]
where $U$ is a fundamental domain for $\Gamma$ and
$dv_{\mathbb{H}}=\frac{1}{2} y^{-2}dx\wedge dy$ is the 
hyperbolic volume form. The Bergman density function of 
$\big(\mathcal{S}_{2p}^\Gamma,\langle
\LargerCdot,\LargerCdot\rangle\big)$
is defined by taking any orthonormal basis $(f_j)$ and setting
\[S^\Gamma_p(z)=\sum_j|f_j(z)|^2(2y)^{2p}\,,\quad z\in U.\]
Under the above isomorphism, $\mathcal{S}_{2p}^\Gamma$
is identified to 
the space
$H^0\big(\overline\Sigma,L^p\otimes
\mathscr{O}_{\overline\Sigma}(D)^{-1}\big)=
H^0\big({\overline\Sigma},K_{\overline\Sigma}^p
\otimes\mathscr{O}_{\overline\Sigma}(D)^{p-1}\big)$
of holomorphic sections of $L^p$ over $\overline\Sigma$
vanishing on $D$.

If we endow $K_\mathbb{H}$ with the Hermitian metric
induced by the Poincar\'e metric on $\mathbb{H}$, 
the scalar product of two elements 
$udz^{\otimes p}, vdz^{\otimes p}\in K^p_{\mathbb{H},z}$
is $\langle udz^{\otimes p}, vdz^{\otimes p}\rangle
=u\overline{v}(2y)^{2p}$.
Hence, the Petersson scalar product corresponds to the
${\bL}^2$ scalar product of pluricanonical forms on $\Sigma$,
\[\langle f,g\rangle=
\int_{\Sigma} \langle\Phi(f),\Phi(g)\rangle\, \omega_{\Sigma}\,,
\quad f,g\in\mathcal{S}_{2p}^\Gamma\,. 
\] 
The isomorphism $\Phi$ gives thus an isometry 
(see also \cite[Section\,6.4]{CM11})
\begin{equation}\label{eq:6.25}
\mathcal{S}_{2p}^\Gamma\cong
H^0\big(\overline\Sigma,L^p\otimes\mathscr{O}_{\overline\Sigma}
(D)^{-1}\big)\cong
H^0_{(2)}(\Sigma,K_\Sigma^p)\cong H^0_{(2)}(\Sigma, L^p),
\end{equation}
where $H^0_{(2)}(\Sigma, L^p)$ is the space of holomorphic
sections of $L^p$ that are square-integrable with respect to 
the volume form $\omega_\Sigma$
and the metric $h^p$ on $L^p$, with $h$ introduced in
Lemma \ref{L:ab}. 
Moreover, $H^0_{(2)}(\Sigma,K_\Sigma^p)$ is the space
of ${\bL}^2$-pluricanonical
sections with respect to the metric $h^{K^p_\Sigma}$
and the volume form $\omega_\Sigma$, where we denote by 
$h^{K_\Sigma}$ the Hermitian metric induced by 
$\omega_\Sigma$ on $K_\Sigma$. 
We let now $B^\Gamma_p$ be the Bergman density function
of $H^0_{(2)}(\Sigma, L^p)$, defined as in \eqref{e:BFS1}. 
We have
\begin{equation}\label{eq:6.25a}
S^\Gamma_p(z)=B^\Gamma_p(\pi(z)),\quad z\in U.
\end{equation}
We thus identify the space of cusp forms 
$\mathcal{S}_{2p}^\Gamma$ to a subspace of holomorphic sections
of $L^p$ by \eqref{eq:6.25}
and its Bergman density function $S^\Gamma_p$ to $B^\Gamma_p$
by \eqref{eq:6.25a}.

\smallskip
\noindent 
     \textit{Proof of Corollary \ref{crl_CrlIntro1}. --- }
In view of Lemma \ref{L:ab}, this follows immediately from 
Corollary \ref{crl_CrlIntro}
applied for the even powers of the square root of 
$L=K_{\overline\Sigma}\otimes\mathscr{O}_{\overline{\Sigma}}(D)$,
and from \eqref{eq:6.25a}.
Even if the square root of $L$ is a formal line bundle
we can apply Corollary 1.3 to its even powers, i.\,e., to $L^p$,
which has the effect of scaling $p$ to $2p$. This explains
the occurrence of $p/\pi$ in the leading term of \eqref{eqn_CrlIntro1}
and \eqref{eqn_CrlIntro2}, as well as in Theorems \ref{t0.5} and \ref{t0.6}.
\cqfd
   
\smallskip
\noindent 
     \textit{Proof of Theorem \ref{t2.5}. --- }
    Since we have the same Sobolev constants
     for $\widetilde{M}/\Gamma$ 
   and $M$ (cf. \cite[Theorem A.1.6, (A.1.15)]{mm}),
   the proof of \cite[Proposition 4.1]{DLM06}
   or of \cite[Proposition 4.1.5]{mm} shows that
   for any $l,m\in \N$, there exists 
   $C_{l,m}>0$ (independent of $\Gamma\subset \pi_{1}(M)$)
   such that for any
   $p\in \N^*, x,y\in \widetilde{M}/\Gamma$ 
   \begin{align}\label{eq:6.31}
\Big|K_{p}^\Gamma (x,y)- B_{p}^\Gamma(x,y)
\Big|_{C^m(\widetilde{M}/\Gamma)} \leq C_{l,m} p^{-l}.
\end{align}
 By the finite propagation speed for the wave operator,
 for $d(x,y)\leq \epsilon$, we have
 \begin{align}\label{eq:6.32}
K_{p}^\Gamma (x,y)=
\pi_{\Gamma}^* K_{p}(\pi_{\Gamma}(x),\pi_{\Gamma}(y)).
\end{align}
 From (\ref{eq:6.31}) and (\ref{eq:6.32}) we
conclude Theorem \ref{t2.5}. \cqfd
    
\smallskip
\noindent 
\textit{Proof of Theorem \ref{t0.5}. --- }
For a subgroup $\Gamma\subset \Gamma_{0}$ of finite index
let $\pi_{\Gamma}: \Gamma\backslash\mathbb{H}
     \to \Gamma_{0}\backslash\mathbb{H}$
     be the canonical finite covering.
We will add a superscript $\Gamma$ for the various
 objects living on $\Gamma\backslash\mathbb{H}$.
  We fix $0<r<e^{-1}$ such that the ends 
$V_{1}, \ldots, V_{N}$ of $\Gamma_{0}\backslash\mathbb{H}$
are of the form $(\D^*_{r},\omega_{\D^*}, \C, h_{\D^*})$.
Now for any finite index subgroup $\Gamma\subset \Gamma_{0}$, 
the ends of 
$\Gamma\backslash\mathbb{H}$ are the connected components
$\{V^\Gamma_{i}\}_{i}$ of $\pi_{\Gamma}^{-1}(V_{j})$,
$j=1,\ldots, N$. Moreover, if $V^\Gamma_{i}$ is 
a connected component of $\pi_{\Gamma}^{-1}(V_{j})$,
there exists $n\in \N^*$ such that the map
$\pi_{\Gamma}: V^\Gamma_{i}\to V_{j}$ is given by
$\D^*_{r^{1/n}}\to \D^*_{r} : z\mapsto z^n$
(cf. \cite[Theorem 5.10]{Fo}). Let 
\begin{align}\label{eq:6.39}
V_{i,r}^\Gamma= V^\Gamma_{i}\cap 
\pi_{\Gamma}^{-1}(\D^*_{r^n}) (= \D^*_{r}).
\end{align}
On $(\Gamma\backslash\mathbb{H}) 
\smallsetminus 
\bigcup_{i}V_{i,r/4}^\Gamma$,
the Sobolev constants are the same as the ones on 
$(\Gamma_{0}\backslash\mathbb{H})
\smallsetminus \bigcup_{j}V_{j,r/4}$,
where $V_{j,r/4}= \D^*_{r/4}$ under the identification
of $V_{j}$ and $\D^*_{r}$.

As the curvatures on $\Gamma\backslash\mathbb{H}$ 
are pull-backs of the corresponding curvatures on 
$\Gamma_{0}\backslash\mathbb{H}$, 
we see from the proof of \cite[Theorem 6.1.1]{mm} that the
spectral gap property, Proposition \ref{prop_SpectralGap},
holds uniformly on the set of subgroups 
$\Gamma\subset \Gamma_{0}$,
i.\,e., there exists $c_{\Gamma_{0}}>0$, $p_{0}>0$
such that for all $p\geq p_{0}$,
and all subgroups $\Gamma\subset \Gamma_{0}$, we have
\begin{align}\label{eq:6.43}
 {\rm Spec}(\square^\Gamma_p) \subset \{0\}\cup
 [c_{\Gamma_{0}}p, +\infty ) \,.
\end{align}

Using \eqref{eq:6.43} and arguing as in the proof of Theorem \ref{t2.5},
we see that for any $k\in \N^*$, there exists $C_{\Gamma_{0},k}>0$
(depending only on $\Gamma_{0}$ and $k$) such that
for any $p\geq 1$, $x\in (\Gamma\backslash\mathbb{H}) 
\smallsetminus \bigcup_{i}V_{i,r/4}^\Gamma$,
\begin{align}\label{eq:6.40}
\Big|B^\Gamma_p(x)- \frac{1}{\pi}p + \frac{1}{2\pi}
\Big|_{C^m}\leq C_{\Gamma_{0},k} p^{-k}.
\end{align}
 Now on $V_{i,r}^\Gamma$ we have
 \begin{align}\label{eq:6.41}
\big(V_{i,r}^\Gamma, \omega_{\Gamma\backslash\mathbb{H}}, 
L, h\big)\simeq (\D^*_{r}, \omega_{\D^*}, \C, h_{\D^*}).
\end{align}
{}From the proof of Proposition \ref{t5.3} and (\ref{eq:6.41}),
by using (\ref{eq:6.43}) in (\ref{eq:5.12a}), we get
that for any $\ell, m \geq 0$, 
   $\gamma>\frac{1}{2}$, there exists 
   $C_{\ell,m,\gamma}>0$ such that for any $p\geq p_{0}$ and any 
   $\Gamma\subset \Gamma_{0}$ with finite index, we have 
\begin{equation}\label{eq:6.45}
      \big\| \rho(x)^{-\gamma}\rho(y)^{-\gamma} K_p^{\Gamma}(x,y)
      \big \|_{C^m(h^p)}
      \leq C_{\ell,m,\gamma}p^{-\ell}, \quad\text{for any }
      x,y\in V_{i,r}^\Gamma.
    \end{equation}
Finally, from (\ref{eq:6.45}) and the proof of Theorem \ref{t6.1}, 
we obtain that for any $\ell, m\geq 0$, and every $\delta >0$, 
  there exists $C_{\Gamma_{0}}>0$ such that
  for all $p\in\N^*$ and any $\Gamma\subset \Gamma_{0}$
  with finite index, 
   \begin{equation}    \label{eq:6.46}
     \Big\|B^\Gamma_p - B_p^{\D^*}
\Big\|_{C^m} (z)\leq C_{\Gamma_{0}}p^{-\ell} 
\big|\!\log(|z|^2)\big|^{-\delta}\quad
\text{ for } \, \, z\in \cup_{i}V_{i,r}^\Gamma.
   \end{equation}
  Now it is clear from Proposition \ref{t3.3}, 
  (\ref{eq:3.20}), (\ref{eq:6.40}) and (\ref{eq:6.46}),
  that (\ref{eqn_CrlIntro2a}) holds.
  \cqfd
  
\smallskip
\noindent  
 \textit{Proof of Theorem \ref{t0.6}. ---}
As in the proof of (\ref{eq:6.31}), there exists $C_{l,m}>0$
(independent of $\Gamma\subset \Gamma_{0}$)
   such that for any  $p\in \N^*$, 
   $x,y\in (\Gamma\backslash\mathbb{H}) \smallsetminus 
\cup_{i}V_{i,r/4}^\Gamma$,
   \begin{align}\label{eq:6.50}
\big|K_{p}^\Gamma (x,y)- B_{p}^\Gamma(x,y)
\big|_{C^m} \leq C_{l,m} p^{-l}.
\end{align}

 Note that the property of the finite propagation speed
of solutions of hyperbolic equations still holds on orbifolds,
as shown in \cite[\S 6.6]{M05}.  
Hence \eqref{eq:6.50} implies that \eqref{eq:6.40} holds for 
$x\in (\Gamma\backslash\mathbb{H}) \smallsetminus 
\big(\bigcup_{i}V_{i,r/4}^\Gamma\cup 
\bigcup_{j}\pi_{\Gamma}^{-1}(U_{x_{j}})\big)$.
Moreover, on each end the argument 
\eqref{eq:6.43}--\eqref{eq:6.46} goes through, 
thus we get the uniform estimate (\ref{eq:6.46})
with the constant $C_{\Gamma_{0}}$ depending only
on $\Gamma_{0}$ and independent of 
$\Gamma\subset \Gamma_{0}$.

Now on each component of $\pi_{\Gamma}^{-1}(U_{x_{j}})$,
observe that the stabilizer group $\Gamma_{x_{j}^\Gamma}$
of $x_{j}^\Gamma$
acts in the normal coordinate around $x_{j}^\Gamma$
in $\mathbb{H}$ by rotation, with $x_{j}^\Gamma$
being the unique fixed point.  We denote by $Z$ the (real) 
normal coordinates around a point $x\in\mathbb{H}$.
By \cite[(4.114), (5.21), (5.23)]{DLM06} 
(cf.\ \cite[Theorem 5.4.11]{mm}, \cite[Theorem 0.2]{DLM12})
 for the $2p$-th tensor power of 
 $(L|_{\Gamma\backslash\mathbb{H}})^{1/2}
 = K_{\Gamma\backslash\mathbb{H}}^{1/2}$, 
there exists $C_{0}>0$  such that
for any $k,l>0$,  there exist $N>0$, $C_{k,l}>0$,
depending only on $\Gamma_{0}, k,l$, such that in 
the normal coordinates around $x_{j}^\Gamma$
in $\mathbb{H}$, we have when $p\to \infty$,
\begin{equation}\label{apb4.20}
\begin{split}
\Bigg |\frac{1}{2p} B_p^\Gamma(Z,Z)
&- \sum_{\nu=0}^k \bb_\nu(Z)(2p)^{-\nu}\\
& - \sum_{\nu=0}^{2k}  (2p)^{-\frac{\nu}{2}}
\sum_{1\neq \gamma\in\Gamma_{\!\!x_j^\Gamma}} 
e^{i \theta_\gamma p}
\mathcal{K}_{\nu}(\sqrt{2p}Z) e^{-
p(1-e^{i\theta_\gamma })|z|^2}
\Bigg |_{C^{l}}\\
&\leqslant C_{k,l} \left(p^{-k-1} + p^{-k+\frac{l-1}{2}}
\left(1+\sqrt{p} |Z|\right)^N
e^{- \sqrt{C_0\, p}\, |Z|}\right),
\end{split}
\end{equation}
where $\mathcal{K}_{\nu}(Z)$ are polynomials 
in the real coordinates $Z$.
Note also that $\bb_\nu$ are given by  (\ref{eq:2.10}).
By \cite[(4.107), (4.105), (4.117), (5.4)]{DLM06}
(or \cite[Remark 4.1.26, (4.1.84), (4.1.92)]{mm}),  we have 
$\mathcal{K}_{0}=\frac{1}{2\pi}, \mathcal{K}_{1}=0$.
This implies in particular that \eqref{eq:0.8} holds on 
$\pi_{\Gamma}^{-1}(U_{x_{j}})$.

The above arguments yield
(\ref{eq:0.7}), (\ref{eq:0.8}) and (\ref{eq:0.10}).
To obtain (\ref{eq:0.9}), notice that 
by our choice of $q_{0}$,
all factors $e^{i\theta_\gamma q_{0}p}$
in (\ref{apb4.20}) are $1$. Thus (\ref{eq:0.8}) implies
(\ref{eq:0.9}). 
  \cqfd

 \appendix

 \section{Proof of Lemma 
 \texorpdfstring{\ref{lem_approx_psip}}{3.4}} \label{app}
  We prove in this appendix the existence of a constant $C$ 
  such that for all $\zeta>0$ and all $p\geq 1$, 
   \begin{equation}   \label{eqn_deltap}
    p\big(1+p(1-\zeta)^2\big)|\delta_p(\zeta)| \leq C, 
   \end{equation}
  where we recall the notation from (\ref{eq:3.21}):
   \begin{equation*}
    \delta_p(\zeta) = e^{p(1-\zeta+\log \zeta)}
    - e^{-\frac{p}{2}(1-\zeta)^2}\Big(1-\frac{p}{3}(1-\zeta)^3\Big). 
   \end{equation*}
  Clearly, \eqref{eqn_deltap} holds for any fixed $p$, that is: 
  for any $p\geq1$, there exists $C_p$ such that for all $\zeta>0$, 
  $p\big(1+p(1-\zeta)^2\big)|\delta_p(\zeta)| \leq C_p$. 
  We thus want to show that the $C_p$ can be chosen 
  independent of $p$; this we do arguing by contradiction: 
  we hence assume that there exists a positive sequence
  $(\zeta_p)_{p\geq1}$ such that, 
  up to passing to a subsequence, 
   \begin{equation}   \label{eqn_deltap2}
    p\big(1+p(1-\zeta_p)^2\big)|\delta_p(\zeta_p)| 
    \xrightarrow{p\to \infty} \infty.  
   \end{equation}
  Up to passing to a subsequence, $(\zeta_p)$ converges in 
  $[0,\infty]$; 
  we first distinguish the three cases $(\zeta_p)\to 0$, 
  $(\zeta_p)\to\infty$, 
  and $(\zeta_p)\to\ell\in\R_{>0}$. 
  \begin{enumerate}
   \item $(\zeta_p)\to 0$: here, 
   $p\big(1+p(1-\zeta_p)^2\big)\sim p^2$,  whereas
   $e^{p(1-\zeta_p+\log \zeta_p)}\leq e^{-p}$ for $p$ large, 
  hence 
  \begin{equation*}
 p\big(1+p(1-\zeta_p)^2\big)e^{p(1-\zeta_p+\log \zeta_p)}
 \longrightarrow 0; 
\end{equation*}
 likewise, $e^{-\frac{p}{2}(1-\zeta)^2}\leq e^{-\frac{p}{4}}$
 for $p$ large, 
and $\big(1-\frac{p}{3}(1-\zeta)^3\big)\sim -\frac{p}{3}$, 
    hence 
\begin{equation}\label{eq:A4}
 p\big(1+p(1-\zeta_p)^2\big)e^{-\frac{p}{2}(1-\zeta_p)^2}
 \big(1-\frac{p}{3}(1-\zeta_p)^3\big) \longrightarrow 0.
\end{equation}
 This way, $p\big(1+p(1-\zeta_p)^2\big)|\delta_p(\zeta_p)|\to 0$, 
 which contradicts \eqref{eqn_deltap2}. 
   \item $(\zeta_p)\to \infty$: one has 
   $p\big(1+p(1-\zeta_p)^2\big)\sim p^2\zeta_p^2$ 
and $e^{p(1-\zeta_p+\log \zeta_p)}\leq e^{-\frac{p\zeta_p}{2}}$ 
for $p$ large, hence 
\begin{equation}\label{eq:A5}
p\big(1+p(1-\zeta_p)^2\big)e^{p(1-\zeta_p+\log \zeta_p)}
\lesssim p^2\zeta_p^2e^{-\frac{p\zeta_p}{2}} \longrightarrow 0; 
\end{equation}
moreover $e^{-\frac{p}{2}(1-\zeta_{p})^2}
\leq e^{-\frac{p\zeta_p}{4}}$ for $p$ large, and 
$\big(1-\frac{p}{3}(1-\zeta_{p})^3\big)\sim \frac{p}{3}\zeta_p^3$,  
hence 
\begin{equation*}
p\big(1+p(1-\zeta_p)^2\big)e^{-\frac{p}{2}(1-\zeta_p)^2}
\big(1-\frac{p}{3}(1-\zeta_p)^3\big) 
\lesssim \frac{p}{3}\zeta_p^3e^{-\frac{p\zeta_p}{4}} 
 = \frac{1}{3p^2} (p\zeta_p)^3e^{-\frac{p\zeta_p}{4}} 
 \longrightarrow 0.
 \end{equation*}
Here again, 
$p\big(1+p(1-\zeta_p)^2\big)|\delta_p(\zeta_p)|\to 0$,
and \eqref{eqn_deltap2} is contradicted. 
   \item $(\zeta_p)\to \ell$: one must deal here with the 
   dichotomy $\ell \neq 1\Big{/}\ell =1$. 
    \begin{enumerate}
     \item $\ell\neq 1$: $p\big(1+p(1-\zeta_p)^2\big) 
     \sim (1-\ell)^2p^2$,  
and as $\zeta \mapsto 1-\zeta+\log\zeta$ is strictly convex 
and attains 0 at $\zeta=1$, 
$e^{p(1-\zeta_p+\log \zeta_p)}\leq e^{-\vareps p}$ 
for $p$ large, with some $\vareps>0$, hence 
 \begin{equation*}
p\big(1+p(1-\zeta_p)^2\big)e^{p(1-\zeta_p+\log \zeta_p)}
\lesssim (1-\ell)^2p^2e^{-\vareps p} \longrightarrow 0; 
\end{equation*}
furthermore, 
$e^{-\frac{p}{2}(1-\zeta_{p})^2}\leq 
e^{-\frac{p(1-\ell)^2}{4}}$
for $p$ large, and 
$\big(1-\frac{p}{3}(1-\zeta_{p})^3\big)\sim \frac{p}{3}(1-\ell)^3$, 
hence 
\begin{equation}\label{eq:A6}
 p\big(1+p(1-\zeta_p)^2\big)e^{-\frac{p}{2}(1-\zeta_p)^2}
 \big|1-\frac{p}{3}(1-\zeta_p)^3\big| 
 \lesssim \frac{p^{3}}{3}|1-\ell|^5e^{-\frac{p(1-\ell)^2}{4}} 
 \longrightarrow 0.
 \end{equation}
Once more, this yields 
$p\big(1+p(1-\zeta_p)^2\big)|\delta_p(\zeta_p)|\to 0$, 
and a contradiction to \eqref{eqn_deltap2}. 
     \item $\ell = 1$: setting $z_p=\zeta_p-1\to 0$, we must 
     one last time distinguish between three different cases: 
$|pz_p^3|\to \infty$, $pz_p^3\to \mu=\lambda^3\in \R^*$ 
and $pz_p^3\to 0$ (up to passing to a subsequence). 
      \begin{enumerate}
       \item $|pz_p^3|\to \infty$: in particular, 
       $pz_p^2=|pz_p^3|^{2/3}p^{1/3}\to \infty$. 
Here $p\big(1+p(1-\zeta_p)^2\big)\sim p^2z_p^2$ 
and $1-\zeta_p+\log\zeta_p = -z_p +\log(1+z_p) 
\leq -\frac{z_p^2}{3}$ for $p$ large, 
thus 
\begin{multline}\label{eq:A8}
p\big(1+p(1-\zeta_p)^2\big)e^{p(1-\zeta_p+\log \zeta_p)}
\lesssim p^2z_p^2e^{-\frac{pz_p^2}{3}} 
= p^2z_p^2e^{-\frac{|pz_p^3|^{2/3}p^{1/3}}{3}}   \\
 \leq p^2z_p^2e^{-(|pz_p^3|^{2/3}+p^{1/3})}
=p^{4/3}e^{-p^{1/3}} |pz_p^3|^{2/3}e^{-|pz_p^3|^{2/3}}  
\longrightarrow 0, 
\end{multline}
  where the last inequality holds as soon as $p^{1/3}$ and
  $|pz_p^3|^{2/3}\geq 6$. As for the other summand, 
 \begin{align}
    p&\big(1+p(1-\zeta_p)^2\big) e^{-\frac{pz_p^2}{2}}
      \Big|1-\frac{p}{3}z_p^3\Big|                           \label{eq:A9} \\
     &\sim  |p^3z_p^5|e^{-\frac{pz_p^2}{2}} 
= p^{4/3}|pz_p^3|^{5/3}e^{-\frac{|pz_p^3|^{2/3}p^{1/3}}{2}} 
\nonumber\\
     &\leq p^{4/3}|pz_p^3|^{5/3}e^{-(|pz_p^3|^{2/3}+p^{1/3})}
      =p^{4/3}e^{-p^{1/3}} |pz_p^3|^{5/3}e^{-|pz_p^3|^{2/3}}    
 \longrightarrow 0,                                               \nonumber
 \end{align}
 (the last inequality holds as soon as $p^{1/3}$ and
 $|pz_p^3|^{2/3}\geq 6$). In conclusion,
$p\big(1+p(1-\zeta_p)^2\big)|\delta_p(\zeta_p)|\to 0$, 
in contradiction with \eqref{eqn_deltap2}. 
 \item $pz_p^3\to \mu=\lambda^3 \neq 0$, that is: 
       $z_p\sim \lambda p^{-1/3}$. 
  First, $p\big(1+p(1-\zeta_p)^2\big) \sim p^2z_p^2 
  \sim \lambda^2p^{4/3}$; 
 moreover $1-\zeta_p+\log\zeta_p = -z_p +\log(1+z_p) 
 = - \frac{z_p^2}{2} + \frac{z_p^3}{3} + \mathcal{O}(p^{-4/3})$, 
   hence
\begin{align}\label{eq:A10}\begin{split}
        p\big(1+p(1-\zeta_p)^2\big)|\delta_p(\zeta_p)| 
        &\sim \lambda^2p^{4/3}e^{-\frac{pz_p^2}{2}}
\Big|e^{\frac{pz_p^3}{3} + \mathcal{O}(p^{-1/3})}-1
+\frac{pz_p^3}{3}\Big|  \\
&\lesssim \lambda^2p^{4/3}e^{-\frac{\lambda^2p^{1/3}}{3}}
\Big|e^{\mu/3}-1-\frac{\mu}{3}\Big|
         \longrightarrow 0,  
\end{split}\end{align}
       a contradiction with \eqref{eqn_deltap2}.
   \item $pz_p^3\to 0$: In this very last case, where again
       $pz_p^4=pz_p^3\cdot z_p\to0$, 
      \begin{align}\label{eq:A11}\begin{split}
      \delta_p(\zeta_p) &= e^{-\frac{pz_p^2}{2}}
      \Big(e^{p(\frac{z_p^3}{3}+\mathcal{O}(z_p^4))}
      - 1-\frac{p}{3}z_p^3\Big)      \\
     &= e^{-\frac{pz_p^2}{2}}\Big(\big(1+\frac{pz_p^3}{3}
     +\mathcal{O}(p^2z_p^6)\big)\big(1+\mathcal{O}(pz_p^4)\big)      
      - 1-\frac{p}{3}z_p^3\Big)      \\
     &= e^{-\frac{pz_p^2}{2}}\big(\mathcal{O}(p^2z_p^6)
     +\mathcal{O}(pz_p^4)\big), 
  \end{split} \end{align}
 hence
 \begin{align}\label{eq:A12}\begin{split}
 (p+p^2z_p^2)\delta_p(\zeta_p) 
 =& \mathcal{O}\big(p^3z_p^6e^{-\frac{pz_p^2}{2}}\big) 
    + \mathcal{O}\big(p^2z_p^4e^{-\frac{pz_p^2}{2}}\big)                       
   + \mathcal{O}\big(p^4z_p^8e^{-\frac{pz_p^2}{2}}\big)     \\
     =&\mathcal{O}(1)
    \end{split}  \end{align} 
 independently of the behavior of $p^2z_p^2$, as 
 $z\mapsto z^ke^{-\frac{z^2}{2}}$, $k=2,3,4$, 
 are bounded on $\R$. In other words, 
 $p\big(1+p(1-\zeta_p)^2\big)\delta_p(\zeta_p)
 =\mathcal{O}(1)$,
 and this final contradiction of \eqref{eqn_deltap2} 
 ends the proof of Lemma \ref{lem_approx_psip}. 
                             \cqfd
      \end{enumerate}
    \end{enumerate}    
  \end{enumerate}
  
 ~
 
 \begin{small}

 \renewcommand{\refname}
 
\end{small}
 

\end{document}